\documentclass{amsart}

\usepackage{a4wide}
\usepackage{verbatim}
\usepackage{amssymb}
\numberwithin{equation}{section}
\newtheorem{thm}{Theorem}[section]
\newtheorem{prop}[thm]{Proposition}
\newtheorem{lemma}[thm]{Lemma}
\newtheorem{cor}[thm]{Corollary}
\theoremstyle{remark}
\newtheorem{remark}[thm]{Remark}

\newcommand{\bul}{\makebox(0,0){$\bullet$}}
\newcommand{\bulo}{\makebox(0,0){$\circ$}}

\title{Properties of generalized univariate hypergeometric functions}
\author{F.J. van de Bult}
\address{F.J. van de Bult, KdV Institute for Mathematics, Universiteit van Amsterdam,
Plantage Muidergracht 24, 1018 TV Amsterdam, The Netherlands.}
\email{fjvdbult@science.uva.nl}
\author{E.M. Rains}
\address{E.M. Rains, Department of Mathematics, University of California, Davis, USA.}
\email{rains@math.ucdavis.edu}
\author{J.V. Stokman}
\address{J.V. Stokman,  KdV Institute for Mathematics, Universiteit van Amsterdam,
Plantage Muidergracht 24, 1018 TV Amsterdam, The Netherlands.}
\email{jstokman@science.uva.nl}

\begin{document}

\begin{abstract}
Based on Spiridonov's analysis of elliptic generalizations of the Gauss hypergeometric
function, we develop a common framework for $7$-parameter families of generalized elliptic, hyperbolic
and trigonometric univariate hypergeometric functions. In each case we derive the symmetries of the generalized
hypergeometric function under the Weyl group of type $E_7$ (elliptic, hyperbolic)
and of type $E_6$ (trigonometric) using the appropriate versions of the Nassrallah-Rahman beta integral, and
we derive contiguous relations using fundamental addition formulas for theta and sine functions. The top level
degenerations of the hyperbolic and trigonometric hypergeometric functions are identified with Ruijsenaars' relativistic hypergeometric function
and the Askey-Wilson function, respectively. We show that the degeneration process yields various
new and known identities for hyperbolic and trigonometric special functions.
We also describe an intimate connection between the hyperbolic and trigonometric
theory, which yields an expression of the hyperbolic hypergeometric function as an explicit bilinear sum
in trigonometric hypergeometric functions.
\end{abstract}

\maketitle

\section{Introduction}

The Gauss hypergeometric function, one of the cornerstones in the theory of classical
univariate special functions, has been generalized in various fundamental directions.
A theory on multivariate root system analogues of the Gauss hypergeometric function, due to Heckman and Opdam,
has emerged, forming the basic tools to solve trigonometric and
hyperbolic quantum many particle systems of Calogero-Moser type and generalizing the Harish-Chandra theory of spherical functions
on Riemannian symmetric spaces (see \cite{HS} and references therein). A further important development has been the
generalization to $q$-special functions, leading to the theory of Macdonald polynomials \cite{Macd}, which play a fundamental role in
the theory of relativistic analogues of the trigonometric quantum Calogero-Moser systems (see e.g. \cite{Ruijsproced})
and in harmonic analysis on quantum compact symmetric spaces (see e.g. \cite{N}, \cite{Let}).
In this paper, we focus on far-reaching generalizations of the
Gauss hypergeometric function within the classes of elliptic, hyperbolic and trigonometric {\it univariate} special functions.

Inspired by results on integrable systems, Ruijsenaars \cite{Ruijs1} defined gamma functions of rational,
trigonometric, hyperbolic and elliptic type. Correspondingly there are four types of special function
theories, with the rational (resp. trigonometric) theory being the
standard theory on hypergeometric (resp. $q$-hypergeometric) special functions, while the hyperbolic
theory is well suited to deal with unimodular base $q$.
The theory of elliptic special functions, initiated by Frenkel and Turaev in \cite{FrenkelTuraev},
is currently in rapid development. The starting point of our analysis is the definition of the various generalized hypergeometric functions
as an explicit hypergeometric integral of elliptic, hyperbolic and trigonometric type depending on seven auxiliary parameters (besides the bases).
The elliptic and hyperbolic analogue of the hypergeometric function are due to Spiridonov \cite{Spir2}, while
the trigonometric analogue of the hypergeometric function is essentially an integral representation of the function $\Phi$ introduced
and studied extensively by Gupta and Masson in \cite{GuptaMasson}. Under a suitable parameter discretization,
the three classes of generalized hypergeometric functions reduce to Rahman's \cite{Ra}
(trigonometric), Spiridonov's  \cite{Spir2} (hyperbolic), and Spiridonov's and Zhedanov's \cite{SZ}, \cite{Spir2} (elliptic)
families of biorthogonal rational functions.

Spiridonov \cite{Spir2} gave an elementary derivation of the symmetry of the elliptic hypergeometric
function with respect to a twisted action of the Weyl group of type $E_7$ on the parameters using the
elliptic analogue \cite{Spir} of the Nassrallah-Rahman \cite{NR} beta integral.
In this paper we follow the same approach to establish the $E_6$-symmetry (respectively $E_7$-symmetry)
of the trigonometric (respectively hyperbolic) hypergeometric function, using now the Nassrallah-Rahman beta integral
(respectively its hyperbolic analogue from \cite{Stokman}).
The $E_6$-symmetry of $\Phi$ has recently been established in \cite{LiJ} by different methods.
Spiridonov \cite{Spir2} also gave elementary derivations of contiguous relations for the
elliptic hypergeometric function using the fundamental addition formula for theta functions (see \eqref{eqellfund}), entailing
a natural elliptic analogue of the Gauss hypergeometric differential equation.
Following the same approach we establish contiguous relations and generalized Gauss hypergeometric equations
for the hyperbolic and trigonometric hypergeometric function.
For $\Phi$ it leads to simple proofs of various results from \cite{GuptaMasson}.

Although the elliptic hypergeometric function is the most general amongst the generalized hypergeometric functions
under consideration (rigorous limits between the different classes of special functions have
been obtained in the recent paper \cite{Rainslimits} of the second author), it is also the most
rigid in its class, in the sense that it does not admit natural degenerations within the class of elliptic
special functions itself (there is no preferred limit point on an elliptic curve). On the other hand, for the hyperbolic and
trigonometric hypergeometric functions various interesting degenerations within their classes are possible, as we point out in
this paper. It leads to many nontrivial identities and results,
some of which are new and some are well known. In any case, it provides new insight in identities, e.g.
as being natural consequences of symmetry breaking in the degeneration process, and it places many identities
and classes of univariate special functions in a larger framework.
For instance, viewing the trigonometric hypergeometric function as a degeneration of the elliptic hypergeometric
function, we show that the breaking of symmetry (from $E_7$ to $E_6$) leads to a second important integral representation of $\Phi$.

Moreover we show that Ruijsenaars' \cite{Ruijs2}
relativistic analogue $R$ of the hypergeometric function is a degeneration of the hyperbolic hypergeometric function,
and that the $D_4$-symmetry \cite{Ruijs3} of $R$ and the four Askey-Wilson second-order difference
equations \cite{Ruijs2} satisfied by $R$ are direct consequences of the $E_7$-symmetry and the contiguous relations of the hyperbolic hypergeometric
function. Similarly, the Askey-Wilson function \cite{KenS} is shown to be a degeneration of the trigonometric hypergeometric
function. In this paper we aim at deriving the symmetries of (degenerate) hyperbolic and trigonometric
hypergeometric functions directly from appropriate hyperbolic and trigonometric beta integral evaluations
using the above mentioned techniques of Spiridonov \cite{Spir2}. The rational level, in which case the Wilson function \cite{G}
appears as a degeneration, will be discussed in a subsequent paper of the first author.

We hope that the general framework proposed in this paper will shed light on the
fundamental, common structures underlying various quantum relativistic Calogero-Moser systems and various quantum
noncompact homogeneous spaces. In the present univariate setting,
degenerations and specializations of the generalized hypergeometric functions play a key role in solving rank one cases of quantum
relativistic integrable Calogero-Moser systems and in harmonic analysis on various quantum $\hbox{SL}_2$ groups.
On the elliptic level, the elliptic hypergeometric function provides solutions of particular cases of
van Diejen's \cite{vD} very general quantum relativistic Calogero-Moser systems of elliptic type (see e.g. \cite{Spir2}),
while elliptic biorthogonal rational functions have been
identified with matrix coefficients of the elliptic quantum $\hbox{SL}_2$ group in \cite{KNR}.  On the hyperbolic level, the Ruijsenaars' R-function
solves the rank one case of a quantum relativistic Calogero-Moser system of hyperbolic type (see \cite{Ruijs4})
and arises as a matrix coefficient of the modular double of the quantum $\hbox{SL}_2$ group (see \cite{vdBult}).
On the trigonometric level, similar results are known for the Askey-Wilson function, which is a degeneration
of the trigonometric hypergeometric function (see \cite{KenS} and \cite{KenS0}).
For higher rank only partial results are known, see e.g. \cite{KH}, \cite{Rainstrafo} (elliptic) and \cite{Stok} (trigonometric).

The outline of the paper is as follows.
In Section 2 we discuss the general pattern of symmetry breaking when integrals with $E_7$-symmetry are degenerated.
In Section 3 we introduce Spiridonov's \cite{Spir2} elliptic hypergeometric function. We shortly recall Spiridonov's \cite{Spir2}
techniques to derive the $E_7$-symmetry and the contiguous relations for the elliptic hypergeometric function.
In Section 4 these techniques are applied for the hyperbolic hypergeometric function and its top level degenerations.
We show that a reparametrization of the top level degeneration of the hyperbolic hypergeometric function is Ruijsenaars' \cite{Ruijs2}
relativistic hypergeometric function $R$. Key properties of $R$, such as a new integral representation,
follow from the symmetries and contiguous relations of the hyperbolic hypergeometric function.
In Section 5 these techniques are considered on the trigonometric level. We link the top level degeneration of the trigonometric
hypergeometric function to the Askey-Wilson function. Moreover, we show that the techniques lead to elementary derivations
of series representations and three term recurrence relations of the various trigonometric integrals. The trigonometric integrals
are contour integrals over indented unit circles in the complex plane, which can be re-expressed as integrals over the real line with indentations
by ``unfolding'' the trigonometric integral.
We show that this provides a link with Agarwal type integral representations of basic hypergeometric series (see \cite[Chapter 4]{GenR}).
Finally, in Section 6 we extend the techniques from \cite{Stokman} to connect the hyperbolic and trigonometric theory.
It leads to an explicit expression of the hyperbolic hypergeometric function as a bilinear sum of trigonometric hypergeometric functions.
For the top level degeneration,
it provides an explicit link between Ruijsenaars' relativistic hypergeometric function and the Askey-Wilson
function.\\

\noindent
\textbf{Acknowledgements:} Rains was supported in part by NSF Grant No. DMS-0401387.
Stokman was supported by the Netherlands Organization for Scientific Research (NWO) in the
VIDI-project ``Symmetry and modularity in exactly solvable models''.

\subsection{Notation}

We denote $\sqrt{\cdot}$ for the branch of the square root $z\mapsto z^{\frac{1}{2}}$
on $\mathbb{C}\setminus\mathbb{R}_{<0}$ with positive values on $\mathbb{R}_{>0}$.

\section{Weyl groups and symmetry breaking}\label{notations}
The root system of type $E_7$ and its parabolic root sub-systems plays an important role in this article. In this section we describe our
specific choice of realization of the root systems and Weyl groups, and we explain the general pattern of symmetry breaking
which arises from degenerating integrals with Weyl group symmetries.

Degeneration of integrals with Weyl group symmetries in general causes symmetry breaking since the direction of degeneration in parameter space is not
invariant under the symmetry group. All degenerations we consider are of the following form.
For a basis $\Delta$ of a given irreducible, finite root system $R$ in Euclidean space $(V,\langle\cdot,\cdot\rangle)$ with
associated Weyl group $W$ we denote
\[V^+(\Delta)=\{v\in V \, | \, \langle v,\alpha\rangle\geq 0 \quad \forall\,\alpha\in\Delta\}
\]
for the associated positive Weyl chamber. We will study integrals $I(u)$ meromorphically depending on a parameter $u\in\mathcal{G}$.
The parameter space will be some complex hyperplane $\mathcal{G}$
canonically isomorphic to the complexification $V_{\mathbb{C}}$ of $V$,
from which it inherits a $W$-action. The integrals under consideration will be $W$-invariant under an associated twisted $W$-action.
We degenerate such integrals by taking limits in parameter space along distinguished directions $v\in V^+(\Delta)$.
The resulting degenerate integrals will thus inherit symmetries with respect to the isotropy subgroup
\[W_v=\{\sigma\in W \,\,\, | \,\,\, \sigma v=v\},
\]
which is a standard parabolic subgroup of $W$ with respect to the given basis $\Delta$,
generated by the simple reflections $s_\alpha$, $\alpha\in\Delta\cap v^\perp$  (since $v \in V^+(\Delta)$).

All symmetry groups we will encounter are parabolic subgroups of the Weyl group $W$ of type $E_8$.
We use in this article the following explicit realization of the root system $R(E_8)$ of type $E_8$.
Let $\epsilon_k$ be the $k$th element of the standard orthonormal basis of $V=\mathbb{R}^8$, with corresponding scalar
product denoted by $\langle\cdot,\cdot\rangle$. We also denote $\langle\cdot,\cdot\rangle$ for its complex bilinear extension
to $\mathbb{C}^8$. We write $\delta=\frac{1}{2}(\epsilon_1+\epsilon_2+\cdots +\epsilon_8)$.
We realize the root system $R(E_8)$ of type $E_8$ in $\mathbb{R}^8$ as
\[R(E_8)=\{v=\sum_{j=1}^8c_j\epsilon_j+c\delta \,\,\, | \,\,\,  \langle v,v\rangle=2,\,\,
c_j,c\in\mathbb{Z} \,\,\hbox{and} \,\, \sum_{j=1}^8c_j \,\, \hbox{even} \}.
\]
For later purposes, it is convenient to have explicit notations for the roots in $R(E_8)$. The roots
are $\pm\alpha_{jk}^+$ ($1\leq j<k\leq 8$), $\alpha_{jk}^-$ ($1\leq j\not=k\leq 8$),
$\beta_{jklm}$ ($1\leq j<k<l<m\leq 8$),  $\pm\gamma_{jk}$  ($1\leq j<k\leq 8$) and $\pm\delta$, where
\begin{align*}
\alpha_{jk}^+ &=\epsilon_j+\epsilon_k,\\
\alpha_{jk}^-&=\epsilon_j-\epsilon_k,\\
\beta_{jklm} & = \frac12(\epsilon_j+\epsilon_k+\epsilon_l+\epsilon_m-\epsilon_n-\epsilon_p-\epsilon_q-\epsilon_r),\\
\gamma_{jk} &=\frac12(-\epsilon_j-\epsilon_k+\epsilon_l+\epsilon_m+\epsilon_n+\epsilon_p+\epsilon_q+\epsilon_r)\\
\end{align*}
and with $(j,k,l,m,n,p,q,r)$ a permutation of $(1,2,3,4,5,6,7,8)$.

The canonical action of the associated Weyl group $W$ on $\mathbb{C}^8$ is determined by
the reflections $s_{\gamma}u = u- \langle u,\gamma \rangle \gamma$ for $u\in \mathbb{C}^8$ and $\gamma\in R(E_8)$.
It is convenient to work with two different choices $\overline{\Delta}_1$, $\overline{\Delta}_2$ of bases for $R(E_8)$, namely
\begin{equation*}
\begin{split}
\overline{\Delta}_1&=\{\alpha_{76}^-, \beta_{1234}, \alpha_{65}^-, \alpha_{54}^-, \alpha_{43}^-, \alpha_{32}^-, \alpha_{21}^-,\alpha_{18}^+\},\\
\overline{\Delta}_2&=\{\alpha_{23}^-,\alpha_{56}^-,\alpha_{34}^-, \alpha_{45}^-, \beta_{5678}, \alpha_{18}^-, \alpha_{87}^-,\gamma_{18}\},
\end{split}
\end{equation*}
with corresponding (affine) Dynkin diagrams
\begin{equation}\label{dyndiae81}
\begin{picture}(210,46)(0,0)  
\put(0,30){\makebox(0,0){$\circ$}}
\put(0,37){\makebox(0,0)[b]{$\delta$}}
\multiput(3,30)(6,0){5}{\line(1,0){3}}
\put(30,30){\bul}
\put(30,37){\makebox(0,0)[b]{$\alpha_{18}^+$}}
\put(60,30){\bul}
\put(60,37){\makebox(0,0)[b]{$\alpha_{21}^-$}}
\put(90,30){\bul}
\put(90,37){\makebox(0,0)[b]{$\alpha_{32}^-$}}
\put(120,30){\bul}
\put(120,37){\makebox(0,0)[b]{$\alpha_{43}^-$}}
\put(150,30){\bul}
\put(150,37){\makebox(0,0)[b]{$\alpha_{54}^-$}}
\put(180,30){\bul}
\put(180,37){\makebox(0,0)[b]{$\alpha_{65}^-$}}
\put(210,30){\bul}
\put(210,37){\makebox(0,0)[b]{$\alpha_{76}^-$}}
\put(150,0){\bul}
\put(30,30){\line(1,0){180}}
\put(150,30){\line(0,-1){30}}
\put(157,0){\makebox(0,0)[l]{$\beta_{1234}$}}
\end{picture}
\end{equation}
and
\begin{equation}\label{dyndiae82}
\begin{picture}(210,46)  
\put(0,30){\makebox(0,0){$\circ$}}
\put(0,37){\makebox(0,0)[b]{$\delta$}}
\multiput(3,30)(6,0){5}{\line(1,0){3}}
\put(30,30){\bul}
\put(30,37){\makebox(0,0)[b]{$\gamma_{18}$}}
\put(60,30){\bul}
\put(60,37){\makebox(0,0)[b]{$\alpha_{87}^-$}}
\put(90,30){\bul}
\put(90,37){\makebox(0,0)[b]{$\alpha_{18}^-$}}
\put(120,30){\bul}
\put(120,37){\makebox(0,0)[b]{$\beta_{5678}$}}
\put(150,30){\bul}
\put(150,37){\makebox(0,0)[b]{$\alpha_{45}^-$}}
\put(180,30){\bul}
\put(180,37){\makebox(0,0)[b]{$\alpha_{34}^-$}}
\put(210,30){\bul}
\put(210,37){\makebox(0,0)[b]{$\alpha_{23}^-$}}
\put(150,0){\bul}
\put(30,30){\line(1,0){180}}
\put(150,30){\line(0,-1){30}}
\put(157,0){\makebox(0,0)[l]{$\alpha_{56}^-$}}
\end{picture}
\end{equation}
respectively, where the open node corresponds to the simple affine root, which we have labeled by the highest root of $R(E_8)$ with respect to the
given basis (which in both cases is given by $\delta\in V^+(\Delta_j)$).
The reason for considering two different basis is the
following: we will see that degenerating an elliptic hypergeometric integral with $W(E_7)$-symmetry
to the trigonometric level in the direction of the basis element $\alpha_{18}^+\in\overline{\Delta}_1$, respectively the basis element
$\gamma_{18}\in\overline{\Delta}_2$, leads to two essentially different trigonometric hypergeometric integrals
with $W(E_6)$-symmetry. The two integrals can be easily related since they arise
as degeneration of the same elliptic hypergeometric integral. This leads directly to highly nontrivial
trigonometric identities, see Section \ref{trigsection} for details.

This remark in fact touches on the basic philosophy of this paper: it is the symmetry breaking in the degeneration of
hypergeometric integrals which lead to various nontrivial identities. It forms an explanation why there
are so many more nontrivial identities on the hyperbolic, trigonometric and rational level when compared
to the elliptic level.

Returning to the precise description of the relevant symmetry groups, we will mainly encounter stabilizer subgroups of the isotropy subgroup $W_{-\delta}$.
Observe that $W_{-\delta}$ is a standard parabolic subgroup
of $W$ with respect to both bases $\overline{\Delta}_j$ since $-\delta\in V^+(\overline{\Delta}_j)$ ($j=1,2$), with associated simple
reflections $s_\alpha$, $\alpha\in\Delta_1:=\overline{\Delta}_1\setminus \{\alpha_{18}^+\}$, respectively
$s_\alpha$, $\alpha\in\Delta_2:=\overline{\Delta}_2\setminus \{\gamma_{18}\}$. Hence $W_{-\delta}$ is isomorphic to the Weyl group
of type $E_7$, and we accordingly write
\[W(E_7):=W_{-\delta}.
\]
We realize the corresponding standard parabolic root system $R(E_7)\subset R(E_8)$ as
\[R(E_7)= R(E_8)\cap \delta^{\perp}\subseteq \delta^\perp\subset \mathbb{R}^8.
\]
Both $\Delta_1$ and $\Delta_2$ form a basis of $R(E_7)$, and the associated (affine) Dynkin diagrams
are given by
\begin{equation}\label{dyndiae71}
\begin{picture}(180,46)  
\put(0,30){\bul}
\put(0,37){\makebox(0,0)[b]{$\alpha_{21}^-$}}
\put(30,30){\bul}
\put(30,37){\makebox(0,0)[b]{$\alpha_{32}^-$}}
\put(60,30){\bul}
\put(60,37){\makebox(0,0)[b]{$\alpha_{43}^-$}}
\put(90,30){\bul}
\put(90,37){\makebox(0,0)[b]{$\alpha_{54}^-$}}
\put(120,30){\bul}
\put(120,37){\makebox(0,0)[b]{$\alpha_{65}^-$}}
\put(150,30){\bul}
\put(150,37){\makebox(0,0)[b]{$\alpha_{76}^-$}}
\put(180,30){\bulo}
\put(180,37){\makebox(0,0)[b]{$\alpha_{78}^-$}}
\put(90,0){\bul}
\put(97,0){\makebox(0,0)[l]{$\beta_{1234}$}}
\put(0,30){\line(1,0){150}}
\put(90,30){\line(0,-1){30}}
\multiput(153,30)(7,0){4}{\line(1,0){3}}
\end{picture}
\end{equation}
and
\begin{equation}\label{dyndiae72}
\begin{picture}(180,46)  
\put(0,30){\bul}
\put(0,37){\makebox(0,0)[b]{$\alpha_{87}^-$}}
\put(30,30){\bul}
\put(30,37){\makebox(0,0)[b]{$\alpha_{18}^-$}}
\put(60,30){\bul}
\put(60,37){\makebox(0,0)[b]{$\beta_{5678}$}}
\put(90,30){\bul}
\put(90,37){\makebox(0,0)[b]{$\alpha_{45}^-$}}
\put(120,30){\bul}
\put(120,37){\makebox(0,0)[b]{$\alpha_{34}^-$}}
\put(150,30){\bul}
\put(150,37){\makebox(0,0)[b]{$\alpha_{23}^-$}}
\put(180,30){\bulo}
\put(180,37){\makebox(0,0)[b]{$\beta_{1278}$}}
\put(90,0){\bul}
\put(97,0){\makebox(0,0)[l]{$\alpha_{56}^-$}}
\put(0,30){\line(1,0){150}}
\put(90,30){\line(0,-1){30}}
\multiput(153,30)(7,0){4}{\line(1,0){3}}
\end{picture}
\end{equation}
respectively (where we have used that $\alpha_{78}^-$, respectively $\beta_{1278}$, is the highest root of $R(E_7)$ with respect to
the basis $\Delta_1$, respectively $\Delta_2$). Note that the root system $R(E_7)$ consists of the roots of the form $\alpha_{jk}^-$ and
$\beta_{jklm}$.

The top level univariate hypergeometric integrals which we will consider in this article depend meromorphically on a parameter $u\in\mathcal{G}_c$
with $\mathcal{G}_c\subset V_{\mathbb{C}}=\mathbb{C}^8$ ($c\in \mathbb{C}$) the complex
hyperplane
\[\mathcal{G}_c=\frac{c}{2}\delta+\delta^\perp=\{u=(u_1,u_2,\ldots,u_8)\in\mathbb{C}^8\,\, |\,\,  \sum_{j=1}^8 u_j = 2 c\}.
\]
The action on $\mathbb{C}^8$ of the isotropy subgroup $W(E_7)=W_{-\delta}\subset W$ preserves the hyperplane $\delta^\perp$ and fixes $\delta$,
hence it canonically acts on $\mathcal{G}_c$. We extend it to an action of the associated affine Weyl group $W_a(E_7)$ of $R(E_7)$ as follows.
Denote $L$ for the ($W(E_7)$-invariant) root lattice $L\subset \delta^\perp$ of $R(E_7)$, defined as the $\mathbb{Z}$-span of all $R(E_7)$-roots.
The affine Weyl group $W_a(E_7)$ is the semi-direct product group
$W_a(E_7) = W(E_7) \ltimes L$.
The action of $W(E_7)$ on $\mathcal{G}_c$ can then be extended to an action of the affine Weyl-group
$W_a(E_7)$ depending on an extra parameter $z\in\mathbb{C}$ by
letting $\gamma\in L$ act as the shift
\[
\tau_\gamma^z u = u - z \gamma,\qquad u\in\mathcal{G}_c.
\]
We suppress the dependence on $z$ whenever its value is implicitly clear from context.

We also use a multiplicative version of the $W(E_7)$-action on $\mathcal{G}_c$.
Consider the action of the group $C_2$ of order two on $\mathbb{C}^8$, with the non-unit element of $C_2$ acting by multiplication by $-1$ of each
coordinate. We define the parameter space $\mathcal{H}_c$ for a parameter $c\in\mathbb{C}^\times=\mathbb{C}\setminus \{0\}$ as
\[
\mathcal{H}_c = \{t=(t_1,\ldots,t_8) \in \mathbb{C}^8 ~|~ \prod_{j=1}^8 t_j = c^2\} / C_2.
\]
Note that this is well defined because if $t$ satisfies $\prod t_i=c^2$, then so does $-t$.
We sometimes abuse notation by simply writing $t=(t_1,\ldots,t_8)$ for the element $\pm t$ in $\mathcal{H}_c$
if no confusion can arise.

We view the parameters $t\in\mathcal{H}_{\exp(c)}$ as the exponential parameters associated to $u\in\mathcal{G}_c$.
Modding out by the action of the $2$-group $C_2$ allows us to put a $W_a(E_7)$-action on $\mathcal{H}_{\exp(c)}$,
which is compatible to the $W_a(E_7)$-action on $\mathcal{G}_c$ as defined above. Concretely, consider the surjective map
$\psi_c:\mathcal{G}_c \to \mathcal{H}_{\exp(c)}$ defined by
\[
\psi_c(u)=\pm(\exp(u_1),\ldots,\exp(u_8)),\qquad u\in\mathcal{G}_c.
\]
For $u\in \mathcal{G}_{c}$ we have $\psi^{-1}_c(\psi_c(u)) = u +2\pi i L$, where $L$ is the root lattice of $R(E_7)$ as defined above.
Since $L$ is $W(E_7)$-invariant, we can now define the action of $W_a(E_7)$ on $\mathcal{H}_{\exp(c)}$ by
$\sigma \psi_c(u) = \psi_c(\sigma u)$, $\sigma\in W_a$ (for any auxiliary parameter $z\in\mathbb{C}$).

Regardless of whether we view the action of the affine Weyl group additively or
multiplicatively, we will use the abbreviated notations $s_{jk} = s_{\alpha_{jk}^-}$, $w=s_{\beta_{1234}}$ and
$\tau_{jk}^z = \tau_{\alpha_{jk}^-}^z$ throughout the article. Note
that $s_{jk}$ ($j\not=k$) acts by interchanging the $j$th and $k$th coordinate. Furthermore, $W(E_7)$ is generated
by the simple reflections $s_\alpha$ ($\alpha\in\Delta_1$), which are the simple permutations  $s_{j,j+1}$ ($j=1,\ldots,6$) and $w$.
The multiplicative action of $w$ on $\mathcal{H}_c$ is explicitly given by
$w(\pm t)=\pm(st_1,st_2,st_3,st_4,s^{-1}t_5,s^{-1}t_6,s^{-1}t_7,s^{-1}t_8)$ where
$s^2=c/t_1t_2t_3t_4=t_5t_6t_7t_8/c$. Finally, note that
the longest element $v$ of the Weyl group $W(E_7)$ acts by multiplication with -1 on the root system $R(E_7)$,
and hence it acts by $v u= c/2-u$ on $\mathcal{G}_c$ and by $v(\pm t)=\pm(c^{\frac{1}{2}}/t_1,\ldots,c^{\frac{1}{2}}/t_8)$
on $\mathcal{H}_c$.

\section{The univariate elliptic hypergeometric function}


\subsection{The elliptic gamma function}

We will use notations which are consistent with \cite{GenR}.
We fix throughout this section two bases $p,q \in \mathbb{C}$ satisfying $|p|,|q|<1$.
The $q$-shifted factorial is defined by
\[
\bigl(a;q\bigr)_{\infty} = \prod_{j=0}^\infty (1-aq^j).
\]
We write $\bigl(a_1,\ldots,a_m;q\bigr)_{\infty}=\prod_{j=1}^m\bigl(a_j;q\bigr)_{\infty}$,
$(az^{\pm 1};q\bigr)_{\infty}=(az,az^{-1};q\bigr)_{\infty}$ etc.
as shorthand notations for products of $q$-shifted factorials.
The renormalized Jacobi theta-function is defined by
\[
\theta(a;q) = \bigl(a,q/a;q\bigr)_\infty.
\]
The elliptic gamma function \cite{Ruijs1}, defined by the infinite
product
\[
\Gamma_e (z;p,q) = \prod_{j,k=0}^\infty
\frac{1-z^{-1}p^{j+1}q^{k+1}}{1-zp^jq^k},
\]
is a meromorphic function in $z\in\mathbb{C}^\times=\mathbb{C}\setminus \{0\}$
which satisfies the difference equation
\begin{align}\label{diffGe}
\Gamma_e(qz;p,q) & = \theta(z;p) \Gamma_e(z;p,q),
\end{align}
satisfies the reflection equation
\begin{align*}
\Gamma_e(z;p,q) & = 1/\Gamma_e(pq/z;p,q),
\end{align*}
and is symmetric in $p$ and $q$,
\begin{equation*}
\Gamma_e(z;p,q) = \Gamma_e(z;q,p).
\end{equation*}
For products of theta-functions and elliptic gamma functions we use
the same shorthand notations as for the $q$-shifted factorial, e.g.
\[
\Gamma_e(a_1,\ldots,a_m;p,q) = \prod_{j=1}^m\Gamma_e(a_j;p,q).
\]

In this section we call a sequence of points a downward (respectively upward) sequence of points if it is of the form
$ap^jq^k$ (respectively $ap^{-j}q^{-k}$) with $j,k\in \mathbb{Z}_{\geq 0}$ for some $a\in \mathbb{C}$.
Observe that the elliptic gamma function $\Gamma_e(az;p,q)$, considered as a meromorphic function in $z$, has poles at
the upward sequence $a^{-1}p^{-j}q^{-k} $ ($j,k \in \mathbb{Z}_{\geq 0}$) of points and has zeros at
the downward sequence $a^{-1}p^{j+1}q^{k+1}$ ($j,k\in \mathbb{Z}_{\geq 0}$) of points.

\subsection{Symmetries of the elliptic hypergeometric function}
The fundamental starting point of our investigations is Spiridonov's \cite{Spir}
elliptic analogue of the classical beta integral,
\begin{equation}\label{eqellbeta}
\frac{(q;q)_\infty (p;p)_\infty}{2} \int_{\mathcal{C}} \frac{\prod_{j=1}^6 \Gamma_e(t_j z^{\pm 1};p,q)}
{\Gamma_e(z^{\pm 2};p,q)} \frac{dz}{2\pi i z} =
\prod_{1\leq j<k\leq 6} \Gamma_e(t_jt_k;p,q)
\end{equation}
for generic parameters $t\in\mathbb{C}^6$ satisfying the balancing condition $\prod_{j=1}^6 t_j = pq$, where
the contour $\mathcal{C}$ is chosen as a deformation
of the positively oriented unit circle $\mathbb{T}$ separating the downward sequences
$t_jp^{\mathbb{Z}_{\geq 0}}q^{\mathbb{Z}_{\geq 0}}$ ($j=1,\ldots,6$) of poles of the integrand
from the upward sequences $t_j^{-1}p^{\mathbb{Z}_{\leq 0}}q^{\mathbb{Z}_{\leq 0}}$ ($j=1,\ldots,6$).
Note here that the factor $1/\Gamma_e(z^{\pm 2};p,q)$ of the integrand is analytic on
$\mathbb{C}^\times$. Moreover, observe that we can take the positively oriented
unit circle $\mathbb{T}$ as contour if the parameters satisfy $|t_j|<1$ ($j=1,\ldots,6$). Several
elementary proofs of \eqref{eqellbeta} are now known, see e.g. \cite{Spir}, \cite{Spiri} and \cite{Rainstrafo}.

We define the integrand $I_e(t;z)=I_e(t;z;p,q)$ for the univariate elliptic hypergeometric function as
\[
I_e(t;z;p,q) = \frac{\prod_{j=1}^8 \Gamma_e(t_jz^{\pm 1};p,q)}{\Gamma_e(z^{\pm 2};p,q)},
\]
where $t=(t_1,t_2,\ldots, t_8) \in \bigl(\mathbb{C}^\times\bigr)^8$.
For parameters $t\in\mathbb{C}^8$ with
$\prod_{j=1}^8 t_j=p^2q^2$ and $t_it_j \not\in p^{\mathbb{Z}_{\leq 0}}q^{\mathbb{Z}_{\leq 0}}$
for $1\leq i,j\leq 8$ (possibly equal), we can define the elliptic hypergeometric
function $S_e(t)=S_e(t;p,q)$ by
\[
S_e(t;p,q) = \int_{\mathcal{C}} I_e(t;z;p,q) \frac{dz}{2\pi i z},
\]
where the contour $\mathcal{C}$ is a deformation of $\mathbb{T}$ which separates the downward
sequences $t_jp^{\mathbb{Z}_{\geq 0}}q^{\mathbb{Z}_{\geq 0}}$ ($j=1,\ldots,8$) of poles of
$I_e(t;\cdot)$ from the upward sequences $t_j^{-1}p^{\mathbb{Z}_{\leq 0}}q^{\mathbb{Z}_{\leq 0}}$
($j=1,\ldots,8$). If the parameters satisfy $|t_j|<1$ this contour can again be taken as the positively
oriented unit circle $\mathbb{T}$.

The elliptic hypergeometric function
$S_e$ extends uniquely to a meromorphic function on
$\{t\in \mathbb{C}^8:\; \prod t_j=p^2q^2\}$. In fact, for a particular value $\tau$ of the parameters
for which the integral is not defined,
we first deform for $t$ in a small open neighborhood of $\tau$ the contour $\mathcal{C}$ to include
those upward poles which collide at $t=\tau$ with downward poles. The resulting expression
is an integral which is analytic at an open neighborhood of $\tau$ plus a sum of residues depending
meromorphically on the parameters $t$. This expression yields the desired meromorphic extension of $S_e(t)$
at $\tau$.
For further detailed analysis of meromorphic dependencies of integrals like $S_e$, see e.g.
\cite{Ruijs2} and \cite{Rainstrafo}.

Since $I_e(t;-z)=I_e(-t;z)$ where $-t=(-t_1,\ldots,-t_8)$, we have
$S_e(t)=S_e(-t)$, hence we can and will view $S_e$ as a meromorphic function
$S_e:\mathcal{H}_{pq} \to \mathbb{C}$. Furthermore, $S_e(t)$
is the special case $II_{\textup{BC}_1}^1$ of Rains' \cite{Rainstrafo}
multivariate elliptic hypergeometric integrals $II_{BC_n}^m$, and it coincides with Spiridonov's
\cite[\S 5]{Spir2} elliptic analogue $V(\cdot)$ of the Gauss hypergeometric function.

\begin{remark}\label{biorthell}
Note that $S_e(t;p,q)$ reduces to the elliptic beta integral \eqref{eqellbeta}
when e.g. $t_1t_6=pq$. More generally, for e.g. $t_1t_6=p^{m+1}q^{n+1}$ ($m,n\in\mathbb{Z}_{\geq 0}$)
it follows from \cite[Thm. 11]{Spi4} that $S_e(t;p,q)$ essentially coincides with
the two-index elliptic biorthogonal rational function $R_{nm}$ of
Spiridonov \cite[App. A]{Spi4}, which is the product of two very-well-poised terminating elliptic
hypergeometric ${}_{12}V_{11}$ series (the second one with the role of the bases $p$ and $q$ reversed).
\end{remark}

Next we determine the explicit $W(E_7)$-symmetries of $S_e(t)$
in terms of the $W(E_7)$ action on $t\in\mathcal{H}_{pq}$ from
Section \ref{notations}. This result was previously obtained by Rains \cite{Rainstrafo}
and by Spiridonov \cite{Spir2}. We give here a proof which is similar to
Spiridonov's \cite[\S 5]{Spir2} proof.
\begin{thm}\label{lemellsym}
The elliptic hypergeometric function $S_e(t)$ \textup{(}$t\in\mathcal{H}_{pq}$\textup{)}
is invariant under permutations
of $(t_1,\ldots,t_8)$ and it satisfies
\begin{equation}\label{eqellwsym}
S_e(t;p,q) = S_e(wt;p,q) \prod_{1\leq j <k\leq 4} \Gamma_e(t_jt_k;p,q) \prod_{5\leq j<k \leq 8} \Gamma_e(t_jt_k;p,q)
\end{equation}
as meromorphic functions in $t\in \mathcal{H}_{pq}$, where \textup{(}recall\textup{)} $w=s_{\beta_{1234}}$.
\end{thm}
\begin{proof}
The permutation symmetry is trivial. To prove \eqref{eqellwsym} we first prove it for parameters
$t \in \mathbb{C}^8$ satisfying $\prod_{j=1}^8t_j=p^2q^2$ and satisfying the additional restraints
$|t_j|<1$ ($j=1,\ldots,8$),
$|t_j|>|pq|^{\frac{1}{3}}$ ($j=5,\ldots,8$) and $|\prod_{j=5}^8t_j|<|pq|$
(which defines a non-empty open subset of parameters of $\{t\in\mathbb{C} \, | \, \prod_{j=1}^8t_j=p^2q^2 \}$
since $|p|,|q|<1$). For these special values of the parameters we consider the double integrand
\begin{align*}
\int_{\mathbb{T}^2}
\frac
  {\prod_{j=1}^4 \Gamma_e(t_jz^{\pm 1};p,q) \Gamma_e(s x^{\pm 1} z^{\pm 1};p,q)
      \prod_{j=5}^8 \Gamma_e(t_j s^{-1} x^{\pm 1})}
  {\Gamma_e(z^{\pm 2},x^{\pm 2};p,q)}
\frac{dz}{2 \pi i z} \frac{dx}{2\pi i x},
\end{align*}
where $s$ is chosen to balance both the $z$ as the $x$ integral,
so $s^2\prod_{j=1}^4 t_j = pq = s^{-2} \prod_{j=5}^8 t_j$. By the additional parameter restraints
we have $|s|<1$ and $|t_j/s|<1$ for $j=5,\ldots,8$, hence the integration contour $\mathbb{T}$
separates the downward pole sequences of the integrand from the upward ones for both integration
variables. Using the elliptic beta integral \eqref{eqellbeta} to integrate this double integral
either first over the variable $z$, or first over the variable $x$, now yields \eqref{eqellwsym}.
Analytic continuation then implies the identity \eqref{eqellwsym} as meromorphic functions on
$\mathcal{H}_{pq}$.
\end{proof}
An interesting equation for $S_e(t)$ arises from Theorem \ref{lemellsym} by considering
the action of the longest element $v$ of $W(E_7)$, using its decomposition
\begin{equation}\label{valternative}
v=s_{45}s_{36}s_{48}s_{37}s_{34}s_{12}ws_{37}s_{48}ws_{35}s_{46}w
\end{equation}
as products of permutations and $w$.
\begin{cor}\label{vell}
We have
\begin{equation}\label{eqellrefl}
S_e(t;p,q) = S_e(vt;p,q) \prod_{1\leq j<k\leq 8} \Gamma_e(t_jt_k;p,q)
\end{equation}
as meromorphic functions in $t\in \mathcal{H}_{pq}$.
\end{cor}
\begin{remark}
Corollary \ref{vell} is the special case $n=m=1$ of \cite[Thm. 3.1]{Rainstrafo},
see also \cite[\S 5, (iii)]{Spir2} for a proof close to our present derivation.
\end{remark}


\subsection{Contiguous relations}\label{contE}

For sake of completeness we recall here Spiridonov's \cite[\S 6]{Spir2} derivation of certain contiguous relations
cq. difference equations for the elliptic hypergeometric function $S_e(t)$
(most notably, Spiridonov's elliptic hypergeometric equation).
The starting point is the fundamental theta function identity
\cite[Exercise 2.16]{GenR},
\begin{equation}\label{eqellfund}
\frac1y\theta(ux^{\pm 1},yz^{\pm 1};p) + \frac1z\theta(uy^{\pm1},zx^{\pm1};p) + \frac1x\theta(uz^{\pm1},xy^{\pm1};p)=0,
\end{equation}
which holds for arbitrary $u,x,y,z\in \mathbb{C}^\times$.
For the $W_a(E_7)$-action on $\mathcal{H}_{pq}$ we take in this subsection
$\tau_{ij}=\tau^{-\log(q)}_{ij}$, which multiplies $t_i$ by $q$ and divides $t_j$ by $q$.
Note that the $q$-difference operators $\tau_{ij}$
are already well defined on $\{t\in \mathbb{C}^8 \, | \, \prod_{j=1}^8 t_j=p^2q^2\}$.

Using the difference equation \eqref{diffGe} of the elliptic gamma function and using \eqref{eqellfund},
we have
\[
\frac{\theta(q^{-1}t_8t_7^{\pm1};p)}{\theta(t_6t_7^{\pm1};p)}I_e(\tau_{68}t;z) +
(t_6 \leftrightarrow t_7) = I_e(t;z),
\]
where $(t_6\leftrightarrow t_7)$ means the same term with $t_6$ and $t_7$ interchanged.
For generic $t\in\mathbb{C}^8$ with $\prod_{j=1}^8t_j=p^2q^2$ we
integrate this equality over $z\in\mathcal{C}$, with $\mathcal{C}$ a deformation of $\mathbb{T}$
which separates the upward and downward pole sequences of all
three integrands at the same time. We obtain
\begin{equation}\label{eqellcont1}
\frac{\theta(q^{-1}t_8t_7^{\pm1};p)} {\theta(t_6t_7^{\pm1};p)}S_e(\tau_{68} t) +
(t_6 \leftrightarrow t_7) = S_e(t)
\end{equation}
as meromorphic functions in $t\in\mathcal{H}_{pq}$.
This equation is also the $n=1$ instance of \cite[Thm. 3.1]{Rainsrecur}.
Note that in both terms on the left hand side the same parameter $t_8$ is divided by $q$,
while two different parameters ($t_6$ and $t_7$) are multiplied by $q$.
We can obtain a different equation (i.e. not obtainable by applying an $S_8$ symmetry to
\eqref{eqellcont1})
by substituting the parameters $vt$ in \eqref{eqellcont1}, where $v\in W(E_7)$ is the longest Weyl group
element, and by using \eqref{eqellrefl}.
The crux is that $\tau_{68}vt = v\tau_{86}t$. We obtain
\begin{equation}\label{eqellcont2}
\frac{\theta(t_7/qt_8;p)}{\theta(t_7/t_6;p)} \prod_{j=1}^5 \theta(t_jt_6/q;p) S_e(\tau_{86} t)
+ (t_6 \leftrightarrow t_7) = \prod_{j=1}^5 \theta(t_jt_8;p) S_e(t)
\end{equation}
for $t\in \mathcal{H}_{pq}$. We arrive at Spiridonov's \cite[\S 6]{Spir2} elliptic hypergeometric equation for $S_e(t)$.
\begin{thm}[\cite{Spir}]
We have
\begin{equation}\label{eqellcont}
A(t) S_e(\tau_{87} t;p,q) + (t_7 \leftrightarrow t_8) = B(t)S_e(t;p,q)
\end{equation}
as meromorphic functions in $t\in\mathcal{H}_{pq}$, where $A$ and $B$ are defined by
\begin{align*}
A(t) &= \frac{1}{t_8\theta(t_7/qt_8,t_8/t_7;p)}  \prod_{j=1}^6 \theta(t_j t_7/q;p) \\
B(t) &= \frac{\theta(t_7t_8/q;p)}{t_6\theta(t_7/qt_6,t_8/qt_6;p)}\prod_{j=1}^5 \theta(t_jt_6;p)
-\frac{\theta(t_6/t_8,t_6t_8;p)}{t_6\theta(t_7/qt_6,t_7/qt_8,t_8/t_7;p)} \prod_{j=1}^5 \theta(t_jt_7/q;p) \\
& -\frac{\theta(t_6/t_7,t_6t_7;p)}{t_6\theta(t_7/t_8,t_8/qt_6,t_8/qt_7;p)} \prod_{j=1}^5 \theta(t_jt_8/q;p).
\end{align*}
\end{thm}
\begin{remark}
Note that $B$ has an $S_6$-symmetry in $(t_1,t_2,\ldots,t_6)$ even though it is not directly apparent from
its explicit representation.
\end{remark}

\begin{proof}
This follows by taking an appropriate combination of three contiguous relations for $S_e(t)$.
Specifically, the three contiguous relations
are \eqref{eqellcont1} and \eqref{eqellcont2} with $t_6$ and $t_8$ interchanged,
and \eqref{eqellcont1} with $t_7$ and $t_8$ interchanged.
\end{proof}

By combining contiguous relations for $S_e(t)$ and exploring the $W(E_7)$-symmetry of $S_e(t)$, one can obtain
various other contiguous relations involving $S_e(\tau_x t)$, $S_e(\tau_y t)$, and $S_e(\tau_z t)$
for suitable root lattice vectors $x, y, z\in L$. A detailed analysis of such procedures is undertaken for
three term transformation formulas on the trigonometric setting by Lievens and Van der Jeugt
\cite{LiJ} (see also Section \ref{trigsection}).

\begin{remark}
Interchanging the role of the bases $p$ and $q$ and using the
symmetry of $S_e(t;p,q)$ in $p$ and $q$, we obtain contiguous relations for $S_e(t;p,q)$ with respect to
multiplicative $p$-shifts in the parameters.
\end{remark}


\section{Hyperbolic hypergeometric integrals}


\subsection{The hyperbolic gamma function}

We fix throughout this section $\omega_1,\omega_2 \in \mathbb{C}$ satisfying $\Re(\omega_1),\Re(\omega_2)>0$,
and we write
\[
\omega = \frac{\omega_1+\omega_2}{2}.
\]
Ruijsenaars' \cite{Ruijs1} hyperbolic gamma function \cite{Ruijs1} is defined by
\[
G(z;\omega_1,\omega_2) = \exp\left( i \int_{0}^\infty \left( \frac{\sin(2zt)}{2\sinh(\omega_1 t)\sinh(\omega_2 t)} - \frac{z}{\omega_1\omega_2 t}
\right) \frac{dt}{t} \right)
\]
for $z\in\mathbb{C}$ satisfying $|\Im(z)|<\Re(\omega)$.
There exists a unique meromorphic extension of $G(\omega_1,\omega_2;z)$ to $z\in \mathbb{C}$ satisfying
\begin{align*}
G(z;\omega_1,\omega_2) &= G(z;\omega_2,\omega_1), \\
G(z;\omega_1,\omega_2) &= G(-z;\omega_1,\omega_2)^{-1}, \\
G(z+i\omega_1;\omega_1,\omega_2;) &= -2 i s((z+i\omega)/\omega_2)
G(z;\omega_1,\omega_2),
\end{align*}
where we use the shorthand notation $s(z)=\sinh(\pi z)$.  In this section we will omit the
$\omega_1,\omega_2$ dependence of $G$ if no confusion is possible, and we formulate all results only with respect to
$i\omega_1$-shifts. We use similar notations for products of hyperbolic gamma functions as for $q$-shifted factorials
and elliptic gamma functions, e.g.
\[G(z_1,\ldots,z_n;\omega_1,\omega_2)=\prod_{j=1}^nG(z_j;\omega_1,\omega_2).
\]

The hyperbolic gamma function $G$ is a degeneration of the elliptic gamma function $\Gamma_e$,
\begin{equation}\label{degG}
\lim_{r\searrow 0}\Gamma_e\bigl(\exp(2\pi irz); \exp(-2\pi \omega_1r), \exp(-2\pi\omega_2r)\bigr)
\exp\left(\frac{\pi(z-i\omega)}{6ir\omega_1\omega_2}\right)=G(z-i\omega;\omega_1,\omega_2)
\end{equation}
for $\omega_1,\omega_2>0$, see \cite[Prop. III.12]{Ruijs1}.

In this section we call a sequence of points a downward (respectively upward) sequence of points if it is of the form
$a+i\mathbb{Z}_{\leq 0}\omega_1+i\mathbb{Z}_{\leq 0}\omega_2$ (respectively
$a + i\mathbb{Z}_{\geq 0}\omega_1+i\mathbb{Z}_{\geq 0}\omega_2$) for some $a\in\mathbb{C}$.
Observe that the hyperbolic gamma function $G(\omega_1,\omega_2;z)$, viewed as meromorphic function in $z\in\mathbb{C}$,
has poles at the downward sequence $-i\omega+i\mathbb{Z}_{\leq 0}\omega_1+i\mathbb{Z}_{\leq 0}\omega_2$
of points and has zeros at the upward sequence
$i\omega+i\mathbb{Z}_{\geq 0}\omega_1+i\mathbb{Z}_{\geq 0}\omega_2$ of points.
The pole of $G(z;\omega_1,\omega_2)$ at $z=-i\omega$ is simple and
\begin{equation}\label{resG}
\underset{z=-i\omega}{\hbox{Res}}\bigl(G(z;\omega_1,\omega_2)\bigr)=\frac{i}{2\pi}\sqrt{\omega_1\omega_2}.
\end{equation}
All contours in this section will be chosen as deformations of the real line $\mathbb{R}$
separating the upward pole sequences of the integrand from the downward ones.

We will also need to know the asymptotic behavior of $G(z)$ as $\Re(z) \to \pm \infty$
(uniformly for $\Im(z)$ in compacta of $\mathbb{R}$).
For our purposes it is sufficient to know that for any $a,b \in \mathbb{C}$ we have
\begin{equation}\label{eqhypasympG}
\lim_{\Re(z)\to \infty} \frac{G(z-a;\omega_1,\omega_2)}{G(z-b;\omega_1,\omega_2)}
\exp\left(\frac{\pi i z}{\omega_1\omega_2} (b-a)\right) =
\exp\left(\frac{\pi i}{2\omega_1\omega_2}(b^2-a^2)\right),
\end{equation}
where the corresponding $o(\Re(z))$-tail as $\Re(z)\rightarrow \infty$ can be estimated uniformly
for $\Im(z)$ in compacta of $\mathbb{R}$, and that for periods
satisfying $\omega_1\omega_2\in\mathbb{R}_{>0}$,
\begin{equation}\label{asymptotics}
|G(u+x;\omega_1,\omega_2)|\leq
M\exp\left(\pi\Im\Bigl(\frac{u}{\omega_1\omega_2}\Bigr)|x|\right),\qquad
\forall\, x\in\mathbb{R}
\end{equation}
for some constant $M>0$, provided that the line $u+\mathbb{R}$
does not hit a pole of $G$. See \cite[Appendix A]{Ruijs1} for details
and for more precise asymptotic estimates.

\subsection{Symmetries of the hyperbolic hypergeometric function}

The univariate hyperbolic beta integral \cite{Stokman} is
\begin{equation}\label{eqhypbetint}
\int_{\mathcal{C}} \frac{G(i\omega \pm 2z;\omega_1,\omega_2)}{\prod_{j=1}^6 G(u_j \pm z;\omega_1,\omega_2)} dz =
2 \sqrt{\omega_1\omega_2} \prod_{1\leq j<k\leq 6} G(i\omega-u_j-u_k;\omega_1,\omega_2)
\end{equation}
for generic $u_1,\ldots,u_6\in \mathbb{C}$ satisfying the additive balancing condition
$\sum_{j=1}^6 u_j=4i\omega$. Note that
this integral converges since the asymptotic behaviour of the integrand at $z=\pm \infty$ is
$\mathcal{O}(\exp(-4\pi |z| \omega/\omega_1\omega_2))$ in view of \eqref{eqhypasympG}.

We define now the integrand of the hyperbolic hypergeometric function $I_h(u;z)=I_h(u;z;\omega_1,\omega_2)$
as
\[
I_h(u;z;\omega_1,\omega_2) = \frac{G(i\omega \pm 2z;\omega_1,\omega_2)}  {\prod_{j=1}^8 G(u_j \pm z;\omega_1,\omega_2)}
\]
for parameters $u\in\mathbb{C}^8$.  The hyperbolic hypergeometric
function $S_h(u)=S_h(u;\omega_1,\omega_2)$ is defined by
\[
S_h(u;\omega_1,\omega_2) = \int_{\mathcal{C}} I_h(u;z;\omega_1,\omega_2) dz
\]
for generic parameters $u\in\mathcal{G}_{2i\omega}$ (see Section \ref{notations} for the definition
of $\mathcal{G}_{2i\omega}$).
The asymptotic behaviour of $I_h(u;z)$ at $z=\pm \infty$ is again
$\mathcal{O}(\exp(-4\pi |z| \omega/\omega_1\omega_2))$, so
the integral absolutely converges. It follows from \eqref{resG} and the analytic difference equations
for the hyperbolic gamma function that $S_h(u)$ has a unique meromorphic extension to $u\in\mathcal{G}_{2i\omega}$,
cf the analysis for the elliptic hypergeometric function $S_e(t)$.
We thus can and will view $S_h(u)$ as a meromorphic function in $u\in\mathcal{G}_{2i\omega}$.
Note furthermore that the real line can be chosen as integration contour in the definition of $S_h(u)$
if $u\in\mathcal{G}_{2i\omega}$ satisfies $\Im(u_j-i\omega)<0$ for all $j$.
The hyperbolic hypergeometric function $S_h(u)$ ($u\in\mathcal{G}_{2i\omega}$) coincides with
Spiridonov's \cite[\S 5]{Spir2} hyperbolic analogue $s(\cdot)$ of the Gauss hypergeometric function.

Using \eqref{degG} and the reflection equation of $G$, we can obtain the
hyperbolic hypergeometric function $S_h(vu;\omega_1,\omega_2)=S_h(i\omega-u_1,\ldots,i\omega-u_8;\omega_1,\omega_2)$
($u\in\mathcal{G}_{2i\omega}$) formally as the degeneration $r\downarrow 0$ of the elliptic hypergeometric
function $S_e(t;p,q)$ with $p=\exp(-2\pi\omega_1r)$, $q=\exp(-2\pi\omega_2r)$ and
$t=\psi_{2i\omega}(2\pi iru)\in \mathcal{H}_{\exp(-4\pi r\omega)}=\mathcal{H}_{pq}$.
This degeneration, which turns out to preserve the $W(E_7)$-symmetry (see below), can be proven rigorously,
see \cite{Rainslimits}. This entails in particular a derivation of the hyperbolic beta integral \eqref{eqhypbetint}
as rigorous degeneration of the elliptic beta integral \eqref{eqellbeta} (see \cite[\S 5.4]{Stokman} for the formal
analysis).

Next we give the explicit $W(E_7)$ symmetries of $S_h(u)$ in terms
of the $W(E_7)$ action on $u\in\mathcal{G}_{2i\omega}$ from Section \ref{notations}.
\begin{thm}\label{lemhypsymS}
The hyperbolic hypergeometric function $S_h(u)$ \textup{(}$u\in \mathcal{G}_{2i\omega}$\textup{)}
is invariant under permutations of $(u_1,\ldots,u_8)$ and it satisfies
\[
S_h(u;\omega_1,\omega_2) = S_h(wu;\omega_1,\omega_2) \prod_{1\leq j<k\leq 4} G(i\omega-u_j-u_k;\omega_1,\omega_2)
\prod_{5\leq j<k\leq 8} G(i\omega-u_j-u_k;\omega_1,\omega_2)
\]
as meromorphic functions in $u\in\mathcal{G}_{2i\omega}$.
\end{thm}
\begin{proof}
The proof is analogous to the proof in the elliptic case (Theorem \ref{lemellsym}).
For the $w$-symmetry we consider for suitable $u\in\mathcal{G}_{2i\omega}$ the double integral
\[
\int_{\mathbb{R}^2}
\frac{G(i\omega \pm 2z,i\omega \pm 2x)}  {\prod_{j=1}^4 G(u_j \pm z) G(i\omega+s\pm x \pm z) \prod_{k=5}^8 G(u_k-s\pm x)}
dzdx
\]
with $s=i\omega-\frac{1}{2}(u_1+u_2+u_3+u_4)=-i\omega+\frac{1}{2}(u_5+u_6+u_7+u_8)$.
We impose the conditions $\Im(s)<0$ and
\begin{equation}\label{c2}
\Im\bigl(u_j-i\omega)<0\quad (j=1,\ldots,4),\qquad \Im(u_k-i\omega)<\Im(s)\quad (k=5,\ldots,8)
\end{equation}
on $u\in\mathcal{G}_{2i\omega}$
to ensure that the upward and downward pole sequences of the integrand of the double integral
are separated by $\mathbb{R}$. Next we show that the parameter
restraints
\begin{equation}\label{c1}
-\Re\bigl(\frac{\omega}{\omega_1\omega_2}\bigr)<\Im\bigl(\frac{s}{\omega_1\omega_2}\bigr)<0
\end{equation}
on $u\in\mathcal{G}_{2i\omega}$ suffice to ensure
absolute convergence of the double integral. Using the reflection equation and
asymptotics \eqref{eqhypasympG} of $G$ we obtain the estimate
\[\frac{1}{|G(i\omega+s\pm x\pm z)|}\leq
M\exp\left(-2\pi\Im\Bigl(\frac{s+i\omega}{\omega_1\omega_2}\Bigr)
\bigl(|z+x|+|z-x|\bigr)\right),\qquad
\forall\, (x,z)\in\mathbb{R}^2
\]
for some constant $M>0$. It follows that the factor $G(i\omega+s\pm x\pm z)^{-1}$
of the integrand is absolutely and uniformly bounded if $\Im\bigl((i\omega+s)/\omega_1\omega_2)\geq
0$, i.e. if $\Im\bigl(s/\omega_1\omega_2\bigr)\geq -\Re\bigl(\omega/\omega_1\omega_2\bigr)$
(observe that $\Re\bigl(\frac{\omega}{\omega_1\omega_2}\bigr)>0$ due to the imposed conditions $\Re(\omega_j)>0$
on the periods $\omega_j$ ($j=1,2$)).
The asymptotic behaviour of the remaining factors of the integrand (which breaks up in factors only depending on $x$
or on $z$) can easily be determined by \eqref{eqhypasympG}, leading finally to the parameter restraints
\eqref{c1} for the absolute convergence of the double integral.

It is easy to verify that the
parameter subset of $\mathcal{G}_{2i\omega}$ defined by the additional restraints $\Im(s)<0$, \eqref{c2} and \eqref{c1}
is non-empty (by e.g. constructing parameters $u\in\mathcal{G}_{2i\omega}$
with small associated balancing parameter $s$).
Using Fubini's Theorem and the hyperbolic beta integral \eqref{eqhypbetint},
we now reduce the double integral to a single integral
by either evaluating the integral over $x$,
or by evaluating the integral over $z$. Using furthermore that
\[wu=(u_1+s,u_2+s,u_3+s,u_4+s,u_5-s,u_6-s,u_7-s,u_8-s)
\]
for $u\in\mathcal{G}_{2i\omega}$, it follows that the resulting identity is the desired $w$-symmetry of $S_h$
for the restricted parameter domain. Analytic continuation now completes the proof.
\end{proof}

The symmetry of $S_h(u)$ ($u\in\mathcal{G}_{2i\omega}$) with respect to
the action of the longest Weyl group element $v\in W(E_7)$ is as follows.
\begin{cor}
The hyperbolic hypergeometric function $S_h$ satisfies
\begin{equation}\label{eqhyprefl}
S_h(u;\omega_1,\omega_2) = S_h(vu;\omega_1,\omega_2) \prod_{1\leq j<k\leq 8} G(i\omega-u_j-u_k;\omega_1,\omega_2)
\end{equation}
as meromorphic functions in $u\in \mathcal{G}_{2i\omega}$.
\end{cor}
\begin{proof}
This follows from Theorem \ref{lemhypsymS} and \eqref{valternative}.
\end{proof}
\subsection{Contiguous relations}

Contiguous relations for the hyperbolic hypergeometric function $S_h$ can be derived in nearly
exactly the same manner as we did for the elliptic hypergeometric function $S_e$ (see Section \ref{contE}
and \cite[\S 6]{Spir2}). We therefore only indicate the main steps.
Using the $p=0$ case of \eqref{eqellfund} we have
\[
s(x \pm v)s(y \pm z) + s(x \pm y)s(z \pm v) + s(x \pm z)s(v \pm y) = 0,
\]
where $s(x\pm v)=s(x+v)s(x-v)$. In this subsection we write $\tau_{jk}=\tau_{jk}^{i\omega_1}$
($1\leq j\not=k\leq 8$), which acts on $u\in\mathcal{G}_{2i\omega}$ by subtracting $i\omega_1$ from $u_j$
and adding $i\omega_1$ to $u_k$.
We now obtain in analogy to the elliptic case the difference equation
\[
\frac{s((u_8 + i \omega \pm (u_7 -i \omega))/\omega_2)}  {s((u_6-i\omega \pm (u_7 -i\omega))/\omega_2)}
 S_h(\tau_{68} u) + (u_6 \leftrightarrow u_7) = S_h(u)
\]
as meromorphic functions in  $u\in \mathcal{G}_{2i\omega}$.
Using \eqref{eqhyprefl} we subsequently obtain
\[
\frac{s((u_7 -u_8 + 2i\omega)/\omega_2)}{s((u_7 - u_6)/\omega_2)}
\prod_{j=1}^5 s((u_j + u_6)/\omega_2)S_h(\tau_{86}u) + (u_6 \leftrightarrow u_7)
=\prod_{j=1}^5 s((u_j+u_8-2i\omega)/\omega_2)S_h(u)
\]
as meromorphic functions in $u\in\mathcal{G}_{2i\omega}$.
Combining these contiguous relations and simplifying we obtain
\begin{equation}\label{eqhypcont3}
A(u)S_h(\tau_{87}u) - (u_7\leftrightarrow u_8)=B(u)S_h(u),\qquad
u\in\mathcal{G}_{2i\omega},
\end{equation}
where
\begin{align*}
A(u) &= s((2i \omega -u_7+u_8)/\omega_2) \prod_{j=1}^6 s((u_j + u_7)/\omega_2), \\
B(u) &=  \frac{s((u_8 \pm u_7)/\omega_2) s((2i \omega +u_8-u_7)/\omega_2)
s((2i \omega -u_8+u_7)/\omega_2)}{s((2i\omega +u_8-u_6)/\omega_2)s((2i \omega +u_7-u_6)/\omega_2)}
\prod_{j=1}^5 s((-2i\omega +u_j+u_6)/\omega_2) \\
& \qquad - \frac{s((2i\omega -u_8+u_7)/\omega_2) s((u_7-u_6)/\omega_2)
s((-2i\omega +u_6+u_7)/\omega_2)}{s((2i\omega+u_8-u_6)/\omega_2)} \prod_{j=1}^5 s((u_j+u_8)/\omega_2) \\
& \qquad + \frac{s((2i\omega +u_8-u_7)/\omega_2)s((u_8-u_6)/\omega_2)
s((-2i\omega+u_6+u_8)/\omega_2)}{s((2i\omega+u_7-u_6)/\omega_2)} \prod_{j=1}^5 s((u_j+u_7)/\omega_2).
\end{align*}
This leads to the following theorem.
\begin{thm}\label{thmhypcontrel}
We have
\begin{equation}\label{eqhypcont}
A(u) (S_h(\tau_{87}u)-S_h(u)) - (u_7 \leftrightarrow u_8) = B_2(u)S_h(u)
\end{equation}
as meromorphic functions in $u\in\mathcal{G}_{2i\omega}$,
where $A(u)$ is as above and with $B_2(u)$ defined by
\begin{equation}\label{eqhypb2}
B_2(u)=\frac{s((u_7 \pm u_8)/\omega_2)s((u_7-u_8 \pm 2i\omega)/\omega_2)}{4}
\left( \sum_{j=7}^8 s(2(i\omega + u_j)/\omega_2)-\sum_{j=1}^6 s(2(i\omega-u_j)/\omega_2)\right).
\end{equation}
\end{thm}
\begin{proof}
It follows from \eqref{eqhypcont3} that \eqref{eqhypcont} holds with $B_2(u)= B(u)-A(u)-A(s_{78}u)$.
The alternative expression \eqref{eqhypb2} for $B_2$ was obtained by Mathematica. Observe though that part of
the zero locus of $B_2(u)$ ($u\in\mathcal{G}_{2i\omega}$) can be predicted in advance. Indeed, the left hand side of \eqref{eqhypcont} vanishes if
$u_7=u_8$ (both terms then cancel each other), and it vanishes if $u_7=u_8\pm i\omega$
(one term vanishes due to a $s$-factor, while the other term vanishes since either $S_h(\tau_{87}u)=S_h(u)$
or $S_h(\tau_{78}u)=S_h(u)$). The zero of $B_2(u)$ at $u_7=-u_8$ can be predicted from the fact that
all hyperbolic hypergeometric functions
$S_h$ in \eqref{eqhypcont} can be evaluated for $u_7=-u_8$ using the hyperbolic beta
integral \eqref{eqhypbetint}.
\end{proof}


\subsection{The degeneration to the hyperbolic Barnes integral}

In this subsection we degenerate the hyperbolic hypergeometric function $S_h(u)$ ($u\in\mathcal{G}_{2i\omega}$)
along the highest root $\beta_{1278}$ of $R(E_7)$ with respect to
the basis $\Delta_2$ of $R(E_7)$ (see \eqref{dyndiae72}
for the associated Dynkin diagram). The resulting degenerate integral $B_h(u)$ thus inherits symmetries with
respect to the standard maximal parabolic subgroup
\[
W_2(D_6):=W(E_7)_{\beta_{1278}}\subset W(E_7),
\]
which is isomorphic to the Weyl group of type $D_6$ and is generated by the simple reflections $s_\alpha$
($\alpha\in\Delta_2\setminus \{\alpha_{23}^-\}$). The corresponding
Dynkin diagram is
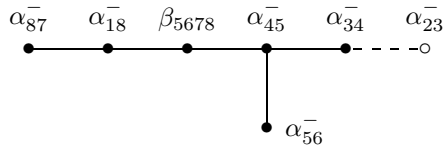
\begin{center}
\begin{picture}(150,46) 
\put(0,30){\bul}
\put(0,37){\makebox(0,0)[b]{$\alpha_{87}^-$}}
\put(30,30){\bul}
\put(30,37){\makebox(0,0)[b]{$\alpha_{18}^-$}}
\put(60,30){\bul}
\put(60,37){\makebox(0,0)[b]{$\beta_{5678}$}}
\put(90,30){\bul}
\put(90,37){\makebox(0,0)[b]{$\alpha_{45}^-$}}
\put(120,30){\bul}
\put(120,37){\makebox(0,0)[b]{$\alpha_{34}^-$}}
\put(150,30){\bulo}
\put(150,37){\makebox(0,0)[b]{$\alpha_{23}^-$}}
\put(90,0){\bul}
\put(97,0){\makebox(0,0)[l]{$\alpha_{56}^-$}}
\put(0,30){\line(1,0){120}}
\put(90,30){\line(0,-1){30}}
\multiput(123,30)(7,0){4}{\line(1,0){3}}
\end{picture}
\end{center}
Concretely, for generic parameters $u\in\mathcal{G}_{2i\omega}$ we define $B_h(u)=B_h(u;\omega_1,\omega_2)$ by
\[
B_h(u;\omega_1,\omega_2) = 2\int_{\mathcal{C}} \frac{\prod_{j=3}^6 G(z-u_j;\omega_1,\omega_2)}{\prod_{j=1,2,7,8} G(z+u_j;\omega_1,\omega_2)} dz.
\]
This integral converges absolutely since the asymptotic behaviour of the integrand at $z=\pm \infty$ is
$\exp(-4\pi \omega |z|/\omega_1\omega_2)$. We may take the real line as integration contour if $u\in\mathcal{G}_{2i\omega}$
satisfies $\Im(u_j-i\omega)<0$ for all $j$. Observe that
the integral $B_h(u)$ has a unique meromorphic extension to $u\in\mathcal{G}_{2i\omega}$. We call $B_h(u)$ the hyperbolic
Barnes integral since it is essentially Ruijsenaars' \cite{Ruijs2}
hyperbolic generalization of the Barnes integral representation of the Gauss hypergeometric
function, see Subsection \ref{reparh}.

\begin{remark}
The parameter space of the hyperbolic Barnes integral $B_h$ is in fact the quotient space $\mathcal{G}_{2i\omega}/\mathbb{C}\beta_{1278}$.
Indeed, for $\xi\in\mathbb{C}$ we have
\[B_h(u+\xi\beta_{1278})=B_h(u)
\]
as meromorphic functions in $u\in\mathcal{G}_{2i\omega}$, which follows by an easy application of \eqref{eqhypasympG} and Cauchy's Theorem.
\end{remark}

\begin{prop} \label{prophyplimB}
For $u\in \mathcal{G}_{2i\omega}$ satisfying $\Im(u_j-i\omega)<0$ \textup{(}$j=1,\ldots,8$\textup{)}
we have
\begin{align*}
\lim_{r\to \infty} &
S_h(u-r\beta_{1278})
\exp\left(\frac{2\pi r \omega}{\omega_1\omega_2}\right) \exp\left(\frac{\pi i }{2\omega_1 \omega_2}\Bigl(\sum_{j=1,2,7,8} u_j^2-\sum_{j=3}^6 u_j^2\Bigr)\right)
=B_h(u).
\end{align*}
\end{prop}
\begin{proof}
The conditions on the parameters $u\in\mathcal{G}_{2i\omega}$ allows us
to choose the real line as integration contour in the integral expression of $S_h(u-r\beta_{1278})$ ($r\in\mathbb{R}$) as well as
in the integral expression of $B_h(u)$.
Using that the integrand $I_h(u;z)$ of $S_h(u)$ is even in $x$, using the reflection equation for the hyperbolic gamma function,
and by a change of integration variable, we have
\begin{equation*}
\begin{split}
S_h(u-r\beta_{1278})e^{\frac{2\pi r \omega}{\omega_1\omega_2}}
&=e^{\frac{2\pi r \omega}{\omega_1\omega_2}}
\int_{-\infty}^\infty  \frac{G(i\omega \pm 2z)}{\prod_{j=1,2,7,8} G(u_j-\frac{r}{2} \pm z)\prod_{j=3}^6 G(u_j +\frac{r}{2} \pm z)} dz\\
&=2e^{\frac{2\pi r \omega}{\omega_1\omega_2}}
\int_0^\infty \frac{G(i\omega \pm 2z)}{\prod_{j=1,2,7,8} G(u_j-\frac{r}{2} \pm z)\prod_{j=3}^6 G(u_j +\frac{r}{2} \pm z)} dz\\
&=2 \int_{-\frac{r}{2}}^\infty k_1(2z+r) k_2(z+r) L(z) dz,
\end{split}
\end{equation*}
where
\begin{align*}
L(z) & = \frac{\prod_{j=3}^6 G(z-u_j)}{\prod_{j=1,2,7,8} G(z+u_j)}, \\
k_1(z) & = \frac{G(z+i\omega)}{G(z-i\omega)} e^{-\frac{2\pi \omega z }{\omega_1\omega_2}}
 = \bigl(1-e^{-2\pi z /\omega_1}\bigr)\bigl( 1-e^{-2\pi z/\omega_2}\bigr), \\
k_2(z) & = \frac{\prod_{j=1,2,7,8} G(z-u_j)}{\prod_{j=3}^6 G(z+u_j)} e^{\frac{4\pi \omega z}{\omega_1\omega_2}}.
\end{align*}
Here the second expression of $k_1$ follows from the analytic difference equations satisfied by $G$.
The pointwise limits of $k_1$ and $k_2$ are
\[
\lim_{z\to \infty} k_1(z) = 1 \qquad \lim_{z\to \infty} k_2(z)=e^{\frac{\pi i}{2\omega_1\omega_2}
(\sum_{j=3}^6 u_j^2 - \sum_{j=1,2,7,8} u_j^2)}.
\]
Moreover, observe that $k_1(z)$ is uniformly bounded for $z\in\mathbb{R}_{\geq 0}$ by $4$, and that $k_2(z)$, being a continuous function
on $\mathbb{R}_{\geq 0}$ with finite limit at infinity, is also uniformly bounded for $z\in\mathbb{R}_{\geq 0}$.

Denote by $\chi_{(-r/2,\infty)}(z)$ the indicator function of the interval $(-r/2,\infty)$.
By Lebesgue's theorem of dominated convergence we now conclude that
\begin{align*}
\lim_{r\to \infty} S_h(u-r\beta_{1278})e^{\frac{2\pi r \omega}{\omega_1\omega_2}} & =
2 \lim_{r\to \infty} \int_{-\frac{r}{2}}^\infty k_1(2z+r)k_2(z+r)L(z)dz  \\
 & = 2\int_{-\infty}^\infty \lim_{r\to \infty} \chi_{(-r/2,\infty)}(z) k_1(2z+s)k_2(z+s)L(z)dz  \\
& = 2e^{\frac{\pi i}{2\omega_1\omega_2}\bigl(\sum_{j=3}^6 u_j^2 - \sum_{j=1,2,7,8} u_j^2)} \int_{-\infty}^\infty L(z)dz \\
& =e^{\frac{\pi i}{2\omega_1\omega_2}(\sum_{j=3}^6 u_j^2 - \sum_{j=1,2,7,8} u_j^2)} B_h(u),
\end{align*}
as desired.
\end{proof}
In the following corollary we use Proposition \ref{prophyplimB} to degenerate the hyperbolic beta integral \eqref{eqhypbetint}.
The resulting integral evaluation formula is an hyperbolic analogue of the nonterminating Saalsch{\"u}tz
formula \cite[(2.10.12)]{GenR}, see Subsection \ref{D5}.

\begin{cor}\label{corhypbarneseval}
For generic $u\in\mathbb{C}^6$ satisfying $\sum_{j=1}^6u_j=4i\omega$ we have
\begin{equation}\label{eqhypbarneseval}
\int_{\mathcal{C}} \frac{G(z-u_4,z-u_5,z-u_6)}{G(z+u_1,z+u_2,z+u_3)} dz=
\sqrt{\omega_1\omega_2} \prod_{j=1}^3 \prod_{k=4}^6 G(i\omega-u_j-u_k).
\end{equation}
\end{cor}
\begin{proof}
Substitute the parameters $u^\prime=(u_1,u_2,u_4,u_5,u_6,0,u_3,0)$ in Proposition
\ref{prophyplimB} with $u_j\in\mathbb{C}$ satisfying
$\Im(u_j-i\omega)<0$ and $\sum_{j=1}^6u_j=4i\omega$. Then $B_h(u^\prime)$ is the left hand side of \eqref{eqhypbarneseval},
multiplied by $2$. On the other hand, by Proposition \ref{prophyplimB} and \eqref{eqhypbetint} we have
\begin{align*}
B_h(u^\prime) &=
\lim_{r\to \infty}
S_h(u^\prime-r\beta_{1278})
\exp\left(
\frac{2\pi r \omega}{\omega_1\omega_2} +
\frac{\pi i}{2\omega_1\omega_2}
(\sum_{j=1}^3 u_j^2 - \sum_{j=4}^6 u_j^2) \right) \\
&=  2\sqrt{\omega_1\omega_2}\prod_{j=1}^3 \prod_{k=4}^6 G(i\omega-u_j-u_k)\\
&\qquad\times\lim_{r\to \infty}\frac{\prod_{1\leq j<k\leq 3} G(i\omega-u_j-u_k+r)}{\prod_{4\leq j<k\leq 6}
G(u_j+u_k-i\omega +r)}\exp\left( \frac{2\pi r \omega}{\omega_1\omega_2} + \frac{\pi i}{2\omega_1\omega_2}
(\sum_{j=1}^3 u_j^2 - \sum_{j=4}^6 u_j^2) \right) \\
&= 2\sqrt{\omega_1\omega_2} \prod_{j=1}^3 \prod_{k=4}^6 G(i\omega-u_j-u_k),
\end{align*}
where the last equality follows from a straightforward but tedious computation using \eqref{eqhypasympG}.
The result for arbitrary generic parameters $u\in\mathbb{C}^6$
satisfying $\sum_{j=1}^6u_j=4i\omega$ now follows by analytic continuation.
\end{proof}

Next we determine the explicit $W_2(D_6)$-symmetries of $B_h(u)$.
\begin{prop}\label{lemhypsymB}
The hyperbolic Barnes integral $B_h(u)$ \textup{(}$u\in\mathcal{G}_{2i\omega}$\textup{)}
is invariant under permutations of $(u_1,u_2,u_7,u_8)$ and of
$(u_3,u_4,u_5,u_6)$ and it satisfies
\begin{equation}\label{eqhypsymBw}
B_h(u) = B_h(wu) \prod_{j=1,2}\prod_{k=3,4} G(i\omega - u_j-u_k) \prod_{j=5,6}\prod_{k=7,8} G(i\omega-u_j-u_k)
\end{equation}
as meromorphic functions in $u\in\mathcal{G}_{2i\omega}$
\end{prop}
\begin{proof}
The permutations symmetry is trivial.
The symmetry \eqref{eqhypsymBw} can be proven by degenerating the
corresponding symmetry of $S_h$, see Theorem \ref{lemhypsymS}.
We prove here the $w$-symmetry by considering the double integral
\[
\int_{\mathbb{R}^2} \frac{G(z-u_3, z-u_4, x-u_5+s, x-u_6+s, z-x-i\omega-s)}
{G(z+u_1,z+u_2,x+u_7-s, x+u_8-s, z-x+i\omega+s)} dz dx
\]
with $s=i\omega-\frac{1}{2}(u_1+u_2+u_3+u_4)=-i\omega+\frac{1}{2}(u_5+u_6+u_7+u_8)$, where
we impose on $u\in\mathcal{G}_{2i\omega}$ the additional conditions
\begin{equation}\label{c1B}
-\Re\bigl(\frac{\omega}{\omega_1\omega_2}\bigr)<\Im\bigl(\frac{s}{\omega_1\omega_2}\bigr)<\Re\bigl(\frac{\omega}{\omega_1\omega_2}\bigr)
\end{equation}
to ensure the absolute convergence of the double integral (this condition is milder than the corresponding condition \eqref{c1}
for $S_h$ due to the missing factors $G(i\omega\pm 2z,i\omega\pm 2x)$ in the numerator of the integrand),
and the conditions \eqref{c2}
to ensure that the upward and downward pole sequences are separated by $\mathbb{R}$.
Using Fubini's Theorem and the hyperbolic beta integral \eqref{eqhypbarneseval}, similarly as in the proof of
Theorem \ref{lemhypsymS}, yields \eqref{eqhypsymBw}.
\end{proof}

\subsection{The degeneration to the hyperbolic Euler integral}

In this subsection we degenerate the hyperbolic hypergeometric function $S_h(u)$ ($u\in\mathcal{G}_{2i\omega}$)
along the highest root $\alpha_{78}^-$ of $R(E_7)$ with respect to
the basis $\Delta_1$ of $R(E_7)$ (see \eqref{dyndiae71}
for the associated Dynkin diagram). The resulting degenerate integral $E_h(u)$ thus inherits symmetries with
respect to the standard maximal parabolic subgroup
\[
W_1(D_6):=W(E_7)_{\alpha_{78}^-}\subset W(E_7),
\]
which is isomorphic to the Weyl group of type $D_6$ and is generated by the simple reflections $s_\alpha$
($\alpha\in\Delta_1\setminus \{\alpha_{76}^-\}$). The corresponding
Dynkin diagram is
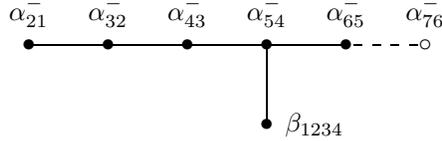
\begin{center}
\begin{picture}(150,46) 
\put(0,30){\bul}
\put(0,37){\makebox(0,0)[b]{$\alpha_{21}^-$}}
\put(30,30){\bul}
\put(30,37){\makebox(0,0)[b]{$\alpha_{32}^-$}}
\put(60,30){\bul}
\put(60,37){\makebox(0,0)[b]{$\alpha_{43}^-$}}
\put(90,30){\bul}
\put(90,37){\makebox(0,0)[b]{$\alpha_{54}^-$}}
\put(120,30){\bul}
\put(120,37){\makebox(0,0)[b]{$\alpha_{65}^-$}}
\put(150,30){\bulo}
\put(150,37){\makebox(0,0)[b]{$\alpha_{76}^-$}}
\put(90,0){\bul}
\put(97,0){\makebox(0,0)[l]{$\beta_{1234}$}}
\put(0,30){\line(1,0){120}}
\put(90,30){\line(0,-1){30}}
\multiput(123,30)(7,0){4}{\line(1,0){3}}
\end{picture}
\end{center}
By the conditions $\Re(\omega_j)>0$ on the periods $\omega_j$ ($j=1,2$) we have that  $\Re\bigl(\frac{\omega}{\omega_1\omega_2}\bigr)>0$.
For generic parameters $u=(u_1,\ldots,u_6)\in\mathbb{C}^6$ satisfying
\begin{equation}\label{eqhypcondeuler}
\Im\bigl(\frac{1}{\omega_1\omega_2} \sum_{j=1}^6 u_j\bigr)>\Re\bigl(\frac{2\omega}{\omega_1\omega_2}\bigr)
\end{equation}
we now define $E_h(u)=E_h(u;\omega_1,\omega_2)$ by
\[
E_h(u;\omega_1,\omega_2) = \int_{\mathcal{C}} \frac{G(i\omega \pm 2z;\omega_1,\omega_2)}{\prod_{j=1}^6 G(u_j \pm z;\omega_1,\omega_2) } dz.
\]
It follows from the asymptotics \eqref{eqhypasympG} of the hyperbolic gamma function that the condition \eqref{eqhypcondeuler} on the parameters
ensures the absolute convergence of $E_h(u)$. Furthermore, $E_h(u)$ admits a
unique meromorphic continuation to parameters $u\in\mathbb{C}^6$
satisfying \eqref{eqhypcondeuler} (in fact, it will be shown later that
$E_h(u)$ has a unique meromorphic continuation to $u\in\mathbb{C}^6$ by relating $E_h$ to the hyperbolic
Barnes integral $B_h$). Observe furthermore that $E_h(u)$ reduces to the hyperbolic beta integral \eqref{eqhypbetint}
when the parameters $u\in\mathbb{C}^6$ satisfy the balancing condition $\sum_{j=1}^6u_j=4i\omega$.
We call $E_h(u)$ the hyperbolic beta integral since its trigonometric analogue is a natural generalization of the Euler
integral representation of the Gauss hypergeometric function, see Subsection \ref{D5} and \cite[\S 6.3]{GenR}.
\begin{prop}\label{prophyplimE}
For $u\in\mathcal{G}_{2i\omega}$ satisfying $\Im(u_j-i\omega)<0$ \textup{(}$j=1,\ldots,8$\textup{)},
$\Im\bigl( (u_7+u_8)/\omega_1\omega_2\bigr)\geq 0$ and \eqref{eqhypcondeuler}, we have
\begin{align}\label{eqi2}
\lim_{r\to \infty}
S_h(u-r\alpha_{78}^-)
\exp\left(-\frac{\pi i}{\omega_1\omega_2} (u_7+u_8)(2r-u_7+u_8)\right) =
E_h(u_1,u_2,u_3,u_4,u_5,u_6).
\end{align}
\end{prop}

\begin{remark}
Proposition \ref{prophyplimE} is trivial when $u_7=-u_8$ due to the
reflection equation for $G$. The resulting limit is the
hyperbolic beta integral \eqref{eqhypbetint} (since the balancing
condition reduces to $\sum_{j=1}^6u_j=4i\omega$).
\end{remark}

\begin{proof}
The assumptions on the parameters ensure that the integration contours in $S_h$ and $E_h$ can be chosen as the real line.
We denote the integrand of the Euler integral by
\[
J(z) = \frac{G(iw \pm 2z)}{\prod_{j=1}^6 G(u_j \pm z)}
\]
and we set
\[
H(z) = \frac{G(z-u_7)}{G(z+u_8)} \exp\left(-\frac{\pi i z}{\omega_1\omega_2}(u_7+u_8) \right).
\]
This allows us to write write
\[I_h(u-r\alpha_{78}^-;z)\exp\left(-\frac{2\pi ir}{\omega_1\omega_2} (u_7+u_8) \right)=J(z)H(r+z)H(r-z),
\]
where (recall) $I_h(u;z)$ is the integrand of the hyperbolic hypergeometric function $S_h(u)$.
Observe that $H$ is a continuous function on $\mathbb{R}$ satisfying
\begin{align*}
\lim_{z\to \infty} H(z) &= \exp\left(\frac{\pi i}{2\omega_1\omega_2}(u_8^2-u_7^2) \right), \\
\lim_{z\to -\infty} H(z) \exp\left(\frac{2\pi i z}{\omega_1\omega_2} (u_7+u_8)\right)  &= \exp\left(\frac{\pi i}{2\omega_1\omega_2}(u_7^2-u_8^2)\right)
\end{align*}
by \eqref{eqhypasympG} and by the reflection equation for the hyperbolic gamma function. Moreover,
$H$ is uniformly bounded on $\mathbb{R}$ in view of the parameter condition $\Im\bigl(u_7+u_8/\omega_1\omega_2\bigr)\geq 0$
on the parameters, and we have
\[
\lim_{r\to \infty} H(r+z)H(r-z) = \exp\left(\frac{\pi i}{\omega_1\omega_2}(u_8^2-u_7^2) \right)
\]
for fixed $z\in\mathbb{R}$.

By Lebesgue's theorem of dominated convergence we conclude that
\begin{align*}
\lim_{r\rightarrow\infty}S_h\bigl(u-r\alpha_{78}^-\bigr) \exp\left(-\frac{2\pi i r}{\omega_1\omega_2}(u_7+u_8)\right) &=
\lim_{r\to \infty} \int_{\mathbb{R}} J(z)H(r+ z)H(r-z) dz \\
& = \int_{\mathbb{R}} J(z) \lim_{r\to \infty} H(r+z)H(r-z) dz \\ &
= E_h(u_1,\ldots,u_6) \exp\left(\frac{\pi i}{\omega_1\omega_2}(u_8^2-u_7^2)\right),
\end{align*}
as desired.
\end{proof}

As a corollary of Proposition \ref{prophyplimE} we obtain the hyperbolic beta integral of Askey-Wilson type,
initially independently proved in \cite{Ruijs4} and in \cite{Stokman}.
\begin{cor}
For generic $u=(u_1,u_2,u_3,u_4)\in\mathbb{C}^4$ satisfying
$\Im\bigl(\frac{1}{\omega_1\omega_2} \sum_{j=1}^4 u_j\bigr)>\Re\bigl(\frac{2\omega}{\omega_1\omega_2}\bigr)$
we have
\begin{equation}\label{eqhypawint}
\int_{\mathcal{C}} \frac{G(i\omega \pm 2z)}{\prod_{j=1}^4 G(u_j \pm z)}\,
dz = 2 \sqrt{\omega_1\omega_2}\,
G(u_1+u_2+u_3+u_4-3i\omega)\prod_{1\leq j<k\leq 4} G(i\omega-u_j-u_k).
\end{equation}
\end{cor}
\begin{proof}
Apply Proposition \ref{prophyplimE} under the additional condition $u_5=-u_6$ on the associated parameters $u\in\mathcal{G}_{2i\omega}$.
Using the reflection equation for the hyperbolic gamma function we see that the right hand side of \eqref{eqi2} becomes
the hyperbolic Askey-Wilson integral. On the other hand,
$S_h(u-r\alpha_{78}^-)$ can be evaluated by the hyperbolic beta integral \eqref{eqhypbetint}, resulting in
\begin{equation*}
\begin{split}
\int_{\mathcal{C}} \frac{G(i\omega \pm 2z)}{\prod_{j=1}^4 G(u_j \pm z)}\, dz&=
2\sqrt{\omega_1\omega_2}\,G(i\omega-u_7-u_8)\prod_{1\leq j<k\leq 4}G(i\omega-u_j-u_k)\\
&\times\lim_{r\rightarrow \infty}\exp\left(-\frac{\pi i}{\omega_1\omega_2}(u_7+u_8)(2r-u_7+u_8)\right)
\prod_{j=1}^4\frac{G(i\omega-u_j-u_7+r)}{G(-i\omega+u_j+u_8+r)}\\
&=2 \sqrt{\omega_1\omega_2}\,
G(u_1+u_2+u_3+u_4-3i\omega)\prod_{1\leq j<k\leq 4} G(i\omega-u_j-u_k),
\end{split}
\end{equation*}
where we used the balancing condition on $u$ and the asymptotics \eqref{eqhypasympG} of the hyperbolic gamma function
to obtain the last equality. The additional parameter restrictions which we have imposed in order to be able to apply Proposition \ref{prophyplimE}
can now be removed by analytic continuation.
\end{proof}

Since both the Euler and Barnes integrals are limits of the hyperbolic hypergeometric function we can connect
them according to the following theorem.
\begin{thm}\label{lemhypEtoB}
We have
\begin{equation}\label{eqhypEtoB}
B_h(u) =E_h(u_2-s,u_7-s,u_8-s,u_3+s,u_4+s,u_5+s)
\prod_{j=3}^5 G(i\omega-u_1-u_j) \prod_{j=2,7,8} G(i\omega-u_6-u_j).
\end{equation}
as meromorphic functions in $\{ u\in\mathcal{G}_{2i\omega}\, | \, \Im\bigl((u_1+u_6)/\omega_1\omega_2\bigr)<\Re\bigl(2\omega/\omega_1\omega_2\bigr)\}$,
where
\[
s =\frac{1}{2}(u_2 + u_6 + u_7 + u_8) - i\omega=i\omega-\frac{1}{2}(u_1+u_3+u_4+u_5).
\]
\end{thm}
\begin{proof}
This theorem can be proved by degenerating a suitable $E_7$-symmetry of $S_h$
using Proposition \ref{prophyplimB} and Proposition
\ref{prophyplimE}. We prove the theorem here directly by analyzing the double
integral
\[
\frac{1}{\sqrt{\omega_1\omega_2}} \int_{\mathbb{R}^2}
\frac{G(i\omega \pm 2z)\prod_{j=3}^5 G(x-u_j)}
{G(i\omega+s+x\pm z)G(x+u_1)\prod_{j=2,7,8} G(u_j-s \pm z)}\, dzdx
\]
for $\Re(\omega_1),\Re(\omega_2)>0$, $u\in \mathcal{G}_{2i\omega}$ and
$s=\frac{1}{2}(u_2+u_6+u_7+u_8)-i\omega$, where we impose the additional parameter
restraints $\omega_1\omega_2\in\mathbb{R}_{>0}$ and
\[|\Im(s)|<\Re(\omega),\qquad \Im(u_6+s)<0
\]
to ensure absolute convergence of the double integral (which
follows from a straightforward analysis of the integrand using \eqref{eqhypasympG} and
\eqref{asymptotics}, cf. the proof of Theorem \ref{lemhypsymS}), and
\[\Im(s)<0,\qquad \Im(i\omega-u_j)>0\quad (j=1,3,4,5),\qquad
\Im(i\omega-u_k+s)>0 \quad (k=2,7,8)
\]
to ensure pole sequence separation by the integration contours.
Note that these parameter restraints imply the parameter condition
$\Im(u_1+u_6)<2\Re(\omega)$ needed for the hyperbolic Euler
integral in the right hand side of \eqref{eqhypEtoB} to be
defined. Integrating the double integral first over $x$ and using the
integral evaluation formula \eqref{eqhypbarneseval} of Barnes
type, we obtain an expression of the double integral as a multiple of
$E_h(u_2-s,u_3+s,u_4+s,u_5+s,u_7-s, u_8-s)$. Integrating first over $z$
and using the hyperbolic Askey-Wilson integral \eqref{eqhypawint},
we obtain an expression of the double integral as a multiple of
$B_h(u)$. The resulting identity is \eqref{eqhypEtoB} for a
restricted parameter domain. Analytic continuation now completes
the proof.
\end{proof}

\begin{cor}
The hyperbolic Euler integral $E_h(u)$ has a unique meromorphic
continuation to $u\in\mathbb{C}^6$ \textup{(}which we also denote by $E_h(u)$\textup{)}.
\end{cor}

{}From the degeneration from $S_h$ to $E_h$ (see Proposition
\ref{prophyplimE}) it is natural to interpret the parameter domain
$\mathbb{C}^6$ as $\mathcal{G}_{2i\omega}/\mathbb{C}\alpha_{78}^-$
via the bijection
\begin{equation}\label{identE6}
\mathbb{C}^6\ni u\mapsto
\bigl(u_1,\ldots,u_6,2i\omega-\sum_{j=1}^6u_j,0\bigr)+\mathbb{C}\alpha_{78}^-.
\end{equation}
We use this identification to transfer the natural
$W_1(D_6)=W(E_7)_{\alpha_{78}^-}$-action on $\mathcal{G}_{2i\omega}/\mathbb{C}\alpha_{78}^-$
to the parameter space $\mathbb{C}^6$ of the hyperbolic Euler integral.
It is generated by permutations of $(u_1,\ldots,u_6)$ and by the
action of $w\in W_1(D_6)$, which is given explicitly by
\begin{equation}\label{w6}
w\bigl(u)=(u_1+s,u_2+s,u_3+s,u_4+s,u_5-s, u_6-s),\qquad
u\in\mathbb{C}^6,
\end{equation}
where $s=i\omega-\frac{1}{2}(u_1+u_2+u_3+u_4)$.
An interesting feature of $W_1(D_6)$-symmetries of the hyperbolic Euler integral
(to be derived in Corollary \ref{inteqEcor}), is the fact
that the nontrivial $w$-symmetry of $E_h$ generalizes to the following
explicit integral transformation for $E_h$.
\begin{prop}\label{inteqE}
For periods $\omega_1, \omega_2\in \mathbb{C}$ with $\Re(\omega_1), \Re(\omega_2)>0$
and $\omega_1\omega_2\not\in\mathbb{R}_{>0}$ and for
parameters $s\in\mathbb{C}$ and $u=(u_1,\ldots,u_6)\in\mathbb{C}^6$ satisfying
\begin{equation}\label{condinteq}
\Im\bigl(\frac{s}{\omega_1\omega_2}\bigr)>-\Re\bigl(\frac{\omega}{\omega_1\omega_2}\bigr),\qquad
\Im\bigl(\frac{u_1+u_2+u_3+u_4}{\omega_1\omega_2}\bigr),
\Im\bigl(\frac{u_5+u_6-2s}{\omega_1\omega_2}\bigr)>\Re\bigl(\frac{2\omega}{\omega_1\omega_2}\bigr)
\end{equation}
and
\begin{equation}\label{polecondinteq}
\Im(u_j-i\omega)<0 \quad (j=1,\ldots,4),\qquad
\Im(u_k-i\omega)<\Im(s)<0 \quad (k=5,6),
\end{equation}
we have
\begin{equation}\label{inttrans}
\begin{split}
\int_{\mathbb{R}}E_h(u_1,u_2,u_3,u_4,i\omega+s+x,i\omega+s-x)&\frac{G(i\omega\pm 2x)}{G(u_5-s\pm x, u_6-s\pm x)}\,dx\\
&=2\sqrt{\omega_1\omega_2}\,\frac{G(i\omega-u_5-u_6+2s)}{G(i\omega-u_5-u_6, i\omega+2s)}\,E_h(u).
\end{split}
\end{equation}
\end{prop}
\begin{proof}
Observe that the requirement $\omega_1\omega_2\not\in\mathbb{R}_{>0}$
ensures the existence of parameters $u\in\mathbb{C}^6$ and $s\in\mathbb{C}$
satisfying the restraints \eqref{condinteq} and
\eqref{polecondinteq}. Furthermore, \eqref{condinteq} ensures that
\[\Im\bigl(\frac{1}{\omega_1\omega_2}\bigl(\sum_{j=1}^6u_j\bigr)\bigr),
\Im\bigl(\frac{1}{\omega_1\omega_2}\bigl(\sum_{j=1}^4u_j+2i\omega+2s\bigr)\bigr)>
\Re\bigl(\frac{2\omega}{\omega_1\omega_2}\bigr),
\]
hence both hyperbolic Euler integrals in \eqref{inttrans} are
defined. We derive the integral transformation \eqref{inttrans}
by considering the double integral
\[\int_{\mathbb{R}^2}\frac{G(i\omega\pm 2z, i\omega\pm 2x)}
{G(i\omega+s\pm x\pm z)\prod_{j=1}^4G(u_j\pm z)\prod_{k=5}^6G(u_k-s\pm
x)}\,dzdx,
\]
which absolutely converges by \eqref{condinteq}.
Integrating the double integral first
over $x$ using the hyperbolic Askey-Wilson integral \eqref{eqhypawint}
yields the right hand side of \eqref{inttrans}. Integrating first over
$z$ results in the left hand side of \eqref{inttrans}.
\end{proof}
\begin{cor}\label{inteqEcor}
The hyperbolic Euler integral $E_h(u)$ \textup{(}$u\in\mathbb{C}^6$\textup{)}
is symmetric in $(u_1,\ldots,u_6)$ and it satisfies
\begin{equation}\label{eqhypsymwE}
E_h(u) = E_h(wu) G(i\omega - u_5-u_6) G(\sum_{j=1}^6 u_j -3i\omega) \prod_{1\leq j<k\leq 4}
G(i\omega-u_j-u_k)
\end{equation}
as meromorphic functions in $u\in\mathbb{C}^6$.
\end{cor}
\begin{proof}
The permutation symmetry is trivial. For \eqref{eqhypsymwE} we
apply Proposition \ref{inteqE} with
$s=i\omega-\frac{1}{2}(u_1+u_2+u_3+u_4)$.
The hyperbolic Euler integral in the left hand side of the
integral transformation \eqref{inttrans} can now be evaluated by the
hyperbolic beta integral \eqref{eqhypbetint}. The remaining
integral is an explicit multiple of $E_h(wu)$. The resulting
identity yields \eqref{eqhypsymwE} for a restricted parameter domain.
Analytic continuation completes the proof.
\end{proof}

\begin{remark}
The $w$-symmetry \eqref{eqhypsymwE} of $E_h$ can also be proved
by degenerating the $w$-symmetry of $S_h$, or by relating
\eqref{eqhypsymwE} to a $W_2(D_6)$-symmetry of
$B_h$ using Theorem \ref{lemhypEtoB}.
\end{remark}

The longest Weyl group element $v_1\in W_1(D_6)$ and the longest Weyl group element $v\in W(E_7)$
have the same action on $\mathcal{G}_{2i\omega}/\mathbb{C}\alpha_{78}^-$. Consequently, under the identification
\eqref{identE6}, $v_1$ acts on $\mathbb{C}^6$ by
\[v_1(u)=(i\omega-u_1,\ldots,i\omega-u_6),\qquad u\in\mathbb{C}^6.
\]

\begin{cor}\label{vsymE}
The symmetry of the hyperbolic Euler integral $E_h(u)$
with respect to the longest Weyl group element $v_1\in W_1(D_6)$
is
\[E_h(u)=E_h(v_1u)G(-3i\omega+\sum_{j=1}^6u_j)
\prod_{1\leq j<k\leq 6}G(i\omega-u_j-u_k)
\]
as meromorphic functions in $u\in\mathbb{C}^6$.
\end{cor}
\begin{proof}
For parameters $u\in\mathcal{G}_{2i\omega}$ such that both $u$ and
$vu$ satisfy the parameter restraints of Proposition
\ref{prophyplimE}, we degenerate the $v$-symmetry \eqref{eqhyprefl} of $S_h$ using \eqref{eqi2}.
Analytic continuation completes the proof.
\end{proof}
The contiguous relations for $S_h$ degenerate to the following contiguous relations for $E_h$.
\begin{lemma}\label{lemhypcontE}
We have
\begin{equation}\label{Econteq}
\begin{split}
\frac{\prod_{j=1}^4 s((u_j+u_5)/\omega_2)}{s((u_5-u_6+2i\omega)/\omega_2)}
(E_h(\tau^{i\omega_1}_{65} u) - E_h(u)) &- (u_5 \leftrightarrow u_6) \\
=&s((u_5 \pm u_6)/\omega_2)s((2i\omega - \sum_{j=1}^6u_j)/\omega_2) E_h(u)
\end{split}
\end{equation}
as meromorphic functions in $u\in\mathbb{C}^6$.
\end{lemma}
\begin{proof}
Use Proposition \ref{prophyplimE} to degenerate the contiguous relation
\eqref{eqhypcont} for the hyperbolic hypergeometric function $S_h$ to
$E_h$.
\end{proof}


\subsection{Ruijsenaars' $R$-function}\label{reparh}

Motivated by the theory of quantum integrable, relativistic particle systems
on the line, Ruijsenaars \cite{Ruijs2}, \cite{Ruijs3}, \cite{Ruijs4} introduced and
studied a generalized hypergeometric $R$-function $R$, which
is essentially the hyperbolic Barnes integral $B_h(u)$ with
respect to a suitable reparametrization (and re-interpretation)
of the parameters $u\in\mathcal{G}_{2i\omega}$.
The new parameters will be denoted by
$(\gamma,x,\lambda)\in \mathbb{C}^6$ with
$\gamma=(\gamma_1,\ldots,\gamma_4)^{T} \in \mathbb{C}^4$,
where $x$ (respectively $\lambda$) is viewed as the geometric (respectively spectral)
parameter, while the four parameters $\gamma_j$ are viewed as coupling
constants. As a consequence of the results derived in the previous
subsections, we will re-derive many of the properties
of the generalized hypergeometric $R$-function, and we
obtain a new integral representation of $R$ in terms of the
hyperbolic Euler integral $E_h$.

Set
\[
N(\gamma) = \prod_{j=1}^3 G(i\gamma_0+i\gamma_j+i\omega).
\]
Ruijsenaars' \cite{Ruijs2} generalized hypergeometric
function $R(\gamma;x,\lambda;\omega_1\omega_2) = R(\gamma,\lambda)$ is defined by
\begin{equation}\label{eqhypdefR}
R(\gamma;x,\lambda) = \frac{1}{2\sqrt{\omega_1\omega_2}} \frac{ N(\gamma) }
{G(i\gamma_0 \pm x,i\hat\gamma_0 \pm \lambda)}B_h(u)
\end{equation}
where $u \in \mathcal{G}_{2i\omega} / \mathbb{C}\beta_{1278}$ with
\begin{equation}\label{eqhypspecialt}
\begin{matrix}
&u_1=i\omega,\quad &u_2=i\omega + i\gamma_0 + i\gamma_1,\quad
&u_3=-i\gamma_0+x,\quad &u_4=-i\gamma_0-x,\\
&u_5=-i\hat\gamma_0+\lambda,\quad
&u_6=-i\hat\gamma_0-\lambda,\quad
&u_7=i\omega + i\gamma_0 + i\gamma_2,\quad
&u_8=i\omega + i\gamma_0 + i\gamma_3.
\end{matrix}
\end{equation}
Note that $R(\gamma;x,\lambda;\omega_1,\omega_2)$ is invariant under permuting the role of the two periods $\omega_1$ and $\omega_2$.
Observe furthermore  that the map  $(\gamma,x,\lambda) \to u+\mathbb{C}\beta_{1278}$, with $u$ given by
\eqref{eqhypspecialt}, defines a bijection
$\mathbb{C}^6 \stackrel{\sim}{\to}
\mathcal{G}_{2i\omega}/\mathbb{C}\beta_{1278}$.

We define the dual parameters $\hat \gamma$ by
\begin{equation}\label{eqhypdefdualpar}
\hat \gamma = \frac{1}{2} \left( \begin{array}{cccc} 1 & 1 & 1 & 1 \\ 1& 1 & -1 & -1 \\ 1 & -1 & 1 & -1 \\ 1 & -1 & -1 & 1
                                 \end{array} \right)  \gamma.
\end{equation}
We will need the following auxiliary function
\begin{align*}
c(\gamma;y) &= \frac{1}{G(2y+i\omega)} \prod_{j=0}^3 G(y-i\gamma_j).
\end{align*}

The following proposition was derived by different methods in \cite{Ruijs3}.
\begin{prop}\label{Rsymmetry}
$R$ is even in $x$ and $\lambda$ and self-dual, i.e.
\begin{align*}
R(\gamma;x,\lambda) = R(\gamma;-x,\lambda) = R(\gamma;x,-\lambda) =
R(\hat \gamma;\lambda,x).
\end{align*}
Furthermore, for an element $\sigma\in W(D_4)$, where $W(D_4)$ is the Weyl-group of type $D_4$ acting
on the parameters $\gamma$ by permutations and even numbers of sign flips, we have
\[
\frac{R(\gamma;x,\lambda)}{c(\gamma;x)c(\hat \gamma;\lambda)N(\gamma)} =
\frac{R(\sigma\gamma;x,\lambda)}{c(\sigma\gamma;x)c(\widehat{\sigma\gamma};\lambda)N(\sigma\gamma)}.
\]
\end{prop}
\begin{proof}
These symmetries are all direct consequences of the $W_2(D_6)$-symmetries
of the hyperbolic Barnes integral $B_h$ (see Proposition \ref{lemhypsymB}).
Concretely, note that the
$W_2(D_6)$-action on $\mathbb{C}^6\simeq \mathcal{G}_{2i\omega}/\mathbb{C}\beta_{1278}$ is
given by
\begin{equation*}
\begin{split}
s_{78}(\gamma,x,\lambda)&=(\gamma_0,\gamma_1,\gamma_3,\gamma_2,x,\lambda),\\
s_{18}(\gamma,x,\lambda)&=(-\gamma_3,\gamma_1,\gamma_2,-\gamma_0,x,\lambda),\\
w(\gamma,x,\lambda)&=(\gamma_1,\gamma_0,\gamma_2,\gamma_3,x,\lambda),\\
s_{45}(\gamma,x,\lambda)&=\bigl(\frac{1}{2}(\gamma_0+\hat{\gamma}_0)+
\frac{i}{2}(x+\lambda),\frac{1}{2}(\gamma_1+\hat{\gamma}_1)-\frac{i}{2}(x+\lambda),
\frac{1}{2}(\gamma_2+\hat{\gamma}_2)-\frac{i}{2}(x+\lambda),\\
&\quad \quad\frac{1}{2}(\gamma_3+\hat{\gamma}_3)-\frac{i}{2}(x+\lambda),
\frac{i}{2}(\hat{\gamma}_0-\gamma_0)+\frac{1}{2}(x-\lambda),
\frac{i}{2}(\hat{\gamma}_0-\gamma_0)+\frac{1}{2}(\lambda-x)\bigr),\\
s_{34}(\gamma,x,\lambda)&=(\gamma,-x,\lambda),\\
s_{56}(\gamma,x,\lambda)&=(\gamma,x,-\lambda).
\end{split}
\end{equation*}
The fact that $R(\gamma;x,\lambda)$ is even in $x$ and $\lambda$
follows now from the $s_{34}\in W_2(D_6)$ and $s_{56}\in W_2(D_6)$ symmetry of $B_h$,
respectively (see Proposition \ref{lemhypsymB}). Similarly, the duality is obtained
similarly from the action of $s_{35}s_{46}$ and using that
$\gamma_0+\gamma_i = \hat \gamma_0 + \hat \gamma_i$ ($i=1,2,3$), while
the $W(D_4)$-symmetry in $\gamma$ follows from considering the action of
$s_{27}\in W_2(D_6)$ (which interchanges $\gamma_1\leftrightarrow \gamma_2$),
$s_{78}$, $s_{18}$ and $w$.
\end{proof}
\begin{remark}
Corollary \ref{corhypbarneseval} implies the explicit evaluation
formula
\[R(\gamma;i\omega+i\gamma_3,\lambda;\omega_1,\omega_2)=
\prod_{j=1}^2\frac{G(i\omega+i\gamma_0+i\gamma_j)}
{G(i\omega+i\gamma_j+i\gamma_3)G(i\hat{\gamma}_j\pm\lambda)}.
\]
Using the $W(D_4)$-symmetry of $R$, this implies
\[R(\gamma;i\omega+i\gamma_0,\lambda;\omega_1,\omega_2)=1,
\]
in accordance with \cite[(3.26)]{Ruijs2}.
\end{remark}

Using Proposition \ref{lemhypsymB} and Theorem \ref{lemhypEtoB} we
can derive several different integral representations of the
$R$-function. First we derive the integral representation of $R$ which was previously derived in \cite{vdBult}
by relating $R$ to matrix coefficients of representations
of the modular double of the quantum group
$\mathcal{U}_q(\mathfrak{sl}_2(\mathbb{C}))$.
\begin{prop}
We have
\begin{align*}
R(\gamma;x,\lambda) = \frac{N(\gamma)}{2\sqrt{\omega_1\omega_2}}
\frac{G(x-i\gamma_0,x-i\gamma_1,\lambda-i\hat\gamma_0,\lambda-i\hat\gamma_1)}
{G(x+i\gamma_2,x+i\gamma_3,\lambda+i\hat\gamma_2,\lambda+i\hat\gamma_3)}
B_h(\upsilon),
\end{align*}
where
\begin{align*}
\upsilon_{1/2} & = x-\frac{\lambda}{2} + \frac{i\omega}{2} \pm \frac{i}{2} (\gamma_0-\gamma_1), &
\upsilon_{3/4} & = -x - \frac{\lambda}{2} + \frac{i\omega}{2} \pm \frac{i}{2} (\gamma_3-\gamma_2), \\
\upsilon_{5/6} &= \frac{\lambda}{2} + \frac{i\omega}{2} \pm \frac{i}{2} (-\gamma_0-\gamma_1), &
\upsilon_{7/8} & = \frac{\lambda}{2} + \frac{i\omega}{2} \pm \frac{i}{2} (\gamma_2+\gamma_3)
\end{align*}
and $\upsilon_{j/k}=\alpha\pm\beta$ means $\upsilon_j=\alpha+\beta$ and $\upsilon_k=\alpha-\beta$.
\end{prop}
\begin{proof}
Express $B_h(s_{36}ws_{35}s_{28}ws_{18}u)$ in terms of $B_h(u)$
using the $W_2(D_6)$-symmetries of the hyperbolic Barnes integral $B_h$ (see Proposition \ref{lemhypsymB})
and specialize $u$ as in \eqref{eqhypspecialt}. This gives the desired equality.
\end{proof}

Moreover we can express $R$ in terms of the hyperbolic Euler integral $E_h$,
which leads to a previously unknown integral representation.
\begin{thm}\label{Eulerreph}
We have
\begin{align*}
R(\gamma;x,\lambda) &= \frac{1}{2\sqrt{\omega_1\omega_2}}
\frac{\prod_{j=1}^3 G(i\gamma_0+i\gamma_j+i\omega,\lambda-i\hat \gamma_j)}
{G(\lambda+i\hat \gamma_0)}E_h(u) \\
&= \frac{1}{2\sqrt{\omega_1\omega_2}} \frac{ \prod_{j=1}^3 G(i\gamma_0+i\gamma_j+i\omega,\lambda - i\hat \gamma_j)}
{G(\lambda + i\hat \gamma_0)\prod_{j=0}^3 G(i\gamma_j \pm x)}E_h(\upsilon),
\end{align*}
where $u\in\mathbb{C}^6$ is given by
\begin{align*}
u_j = \frac{i\omega}{2}+i\gamma_{j-1} -\frac{i\hat\gamma_0}{2}  +\frac{\lambda}{2}, \qquad (j=1,\ldots,4),
\qquad \qquad u_{5/6}  = \frac{i\omega}{2} \pm x +\frac{i\hat \gamma_0}{2} -\frac{\lambda}{2},
\end{align*}
and $\upsilon\in\mathbb{C}^6$ is given by
\begin{align*}
 \upsilon_j = \frac{i\omega}{2}-i\gamma_{j-1} +\frac{i\hat\gamma_0}{2} & +\frac{\lambda}{2}, \qquad (j=1,\ldots,4),
 \qquad \qquad
 \upsilon_{5/6}  = \frac{i\omega}{2} \pm x -\frac{i\hat \gamma_0}{2} -\frac{\lambda}{2}.
\end{align*}
\end{thm}
\begin{proof}
To prove the first equation, express $R(\gamma;x,\lambda)$ in terms of
$R(-\gamma_3,\gamma_1,\gamma_2,-\gamma_0;x,\lambda)$ using the
$W(D_4)$-symmetry of $R$ (see Proposition \ref{Rsymmetry}).
Subsequently use the identity relating $B_h$ to $E_h$, see Theorem \ref{lemhypEtoB}.
To obtain the second equation, apply the symmetry of $E_h$ with
respect to the longest Weyl-group element $v_1\in W_1(D_6)$ (see Corollary \ref{vsymE})
in the first equation and use that $R$ is even in $\lambda$.
\end{proof}
The contiguous relation for $E_h$ (Lemma \ref{lemhypcontE}) now becomes the following result.
\begin{prop}[\cite{Ruijs2}]\label{AWdifflemma}
Ruijsenaars' $R$-function satisfies the Askey-Wilson second order difference equation
\begin{equation}\label{AW1}
A(\gamma;x;\omega_1,\omega_2) (R(\gamma;x+i\omega_1, \lambda)-R(\gamma;x,\lambda)) +
(x\leftrightarrow -x)=
B(\gamma;\lambda;\omega_1,\omega_2) R(\gamma;x,\lambda),
\end{equation}
where
\begin{align*}
A(\gamma;x;\omega_1,\omega_2) & = \frac{\prod_{j=0}^3s((i \omega + x + i\gamma_j)/\omega_2)}
{s(2x/\omega_2)s(2(x+i\omega)/\omega_2)}, \\
B(\gamma;\lambda;\omega_1,\omega_2) &= s((\lambda - i\omega-i\hat\gamma_0)/\omega_2)
s((\lambda+i\omega+i\hat\gamma_0)/\omega_2).
\end{align*}
\end{prop}
\begin{remark}
As is emphasized in \cite{Ruijs2}, $R$ satisfies four
Askey-Wilson second order difference equations; two equations
acting on the geometric variable $x$ (namely \eqref{AW1}, and \eqref{AW1} with
the role of $\omega_1$ and $\omega_2$ interchanged), as well as
two equations acting on the spectral parameter $\lambda$ by
exploring the duality of $R$ (see Proposition \ref{Rsymmetry}).
\end{remark}

For later purposes, it is convenient to rewrite \eqref{AW1} as the eigenvalue equation
\[\bigl(\mathcal{L}_{\gamma}^{\omega_1,\omega_2}R(\gamma;\,\cdot\,,\lambda;\omega_1,\omega_2)\bigr)(x)=
B(\gamma;\lambda;\omega_1,\omega_2)R(\gamma;x,\lambda;\omega_1,\omega_2)
\]
for the Askey-Wilson second order difference operator
\begin{equation}\label{AWoperator}
\bigl(\mathcal{L}_{\gamma}^{\omega,\omega_2}f\bigr)(x):=A(\gamma;x;\omega_1,\omega_2)\bigl(f(x+i\omega_1)-f(x)\bigr)+\bigl(x\leftrightarrow -x).
\end{equation}

\section{Trigonometric hypergeometric integrals}\label{trigsection}

\subsection{Basic hypergeometric series}
In this section we assume that the base $q$ satisfies $0<|q|<1$.
The trigonometric gamma function \cite{Ruijs1}
is essentially the $q$-gamma function $\Gamma_q(x)$, see \cite{GenR}.
For ease of presentation we express all the results in terms of the $q$-shifted factorial $\bigl(z;q\bigr)_{\infty}$,
which are related to $\Gamma_q(x)$ by
\[\Gamma_q(x)=\frac{\bigl(q;q\bigr)_{\infty}}{\bigl(q^x;q\bigr)_{\infty}}(1-q)^{1-x}
\]
(with a proper interpretation of the right hand side). The $q$-shifted factorial is the $p=0$ degeneration of
the elliptic gamma function,
\begin{equation}\label{elltrigred}
\Gamma_e(z;0,q)=\frac{1}{\bigl(z;q\bigr)_{\infty}}.
\end{equation}
while the role of the first order analytic difference
equation is taken over by
\[
\bigl(z;q\bigr)_\infty = (1-z)\bigl(qz;q\bigr)_\infty.
\]
However there is no reflection equation anymore; its role is taken over by the product formula
for Jacobi's (renormalized) theta function
\[\theta(z;q)=\bigl(z,q/z;q\bigr)_\infty.
\]
As a function of $z$ the $q$-shifted factorial $\bigl(z;q\bigr)_{\infty}$ is entire with zeros at $z=q^{-n}$ for
$n\in \mathbb{Z}_{\geq 0}$. In this section we call a sequence of the form $aq^{-n}$ ($n\in \mathbb{Z}_{\geq 0}$) an upward
sequence (since they diverge to infinity for large $n$) and a sequence of the form $aq^n$ ($n\in \mathbb{Z}_{\geq 0}$) a
downward sequence (as the elements converge to zero for large $n$).

We will use standard notations for basic hypergeometric series
from \cite{GenR}. In particular, the ${}_{r+1}\phi_r$ basic
hypergeometric series is
\[{}_{r+1}\phi_r\left(\begin{matrix} a_1,\ldots,a_{r+1}\\
b_1,\ldots
b_r\end{matrix};q,z\right)=\sum_{n=0}^{\infty}\frac{\bigl(a_1,\ldots,a_{r+1};q\bigr)_n}
{\bigl(q,b_1,\ldots,b_r;q\bigr)_n}z^n,\qquad |z|<1,
\]
where $\bigl(a;q\bigr)_n=\prod_{j=0}^{n-1}(1-aq^j)$ and with the usual
convention regarding products of such expressions. The very-well-poised ${}_{r+1}\phi_r$ basic hypergeometric
series is
\[{}_{r+1}W_r\bigl(a_1;a_4,a_5,\ldots,a_{r+1};q,z\bigr)=
{}_{r+1}\phi_r\left(\begin{matrix} a_1, qa_1^{\frac{1}{2}},
-qa_1^{\frac{1}{2}}, a_4,\ldots, a_{r+1}\\
a_1^{\frac{1}{2}}, -a_1^{\frac{1}{2}}, qa_1/a_4,\ldots,
qa_1/a_{r+1}\end{matrix}; q,z\right).
\]
Finally, the bilateral basic hypergeometric series ${}_r\psi_r$ is defined as
\begin{equation*}
\begin{split}
{}_r\psi_r\left(\begin{matrix} a_1,a_2,\ldots,a_r\\ b_1,b_2,\ldots,b_r\end{matrix};q,z\right)&=
\sum_{n=0}^{\infty}\frac{\bigl(a_1,a_2,\ldots,a_r;q\bigr)_n}{\bigl(b_1,b_2,\ldots,b_r;q\bigr)_n}z^n\\
&+\sum_{n=1}^{\infty}\frac{\bigl(q/b_1,q/b_2,\ldots,q/b_r;q\bigr)_n}
{\bigl(q/a_1,q/a_2,\ldots,q/a_r;q\bigr)_n}\left(\frac{b_1\cdots
b_r}{a_1\cdots a_rz}\right)^n,
\end{split}
\end{equation*}
provided that $|b_1\cdots b_r/a_1\cdots a_r|<|z|<1$ to ensure absolute and uniform convergence.

We end this introductory subsection by an elementary lemma which we will enable us to rewrite trigonometric
integrals with compact integration cycle in terms of trigonometric integrals with noncompact integration cycle.
Let $\mathbb{H}_+$ be the upper half plane in $\mathbb{C}$. In this section we choose
$\tau\in\mathbb{H}_+$ such that $q=e(\tau)$ once and for all,
where $e(x)$ is a shorthand notation for $\exp(2\pi ix)$. We
furthermore write $\Lambda=\mathbb{Z}+\mathbb{Z}\tau$.

\begin{lemma}\label{unfoldlemma}
Let $u,v\in\mathbb{C}$ such that $u\not\in v+\Lambda$. There
exists an $\eta=\eta(u,v)\in \mathbb{C}$, unique up to
$\Lambda$-translates, such that
\begin{equation}\label{unfoldformula}
\begin{split}
\frac{\theta\bigl(e(u+v-\eta-x), e(x-\eta);q\bigr)}{\theta\bigl(e(u-\eta), e(v-\eta);q\bigr)}&=
\frac{\bigl(e((v-u)/\tau)-1\bigr)\theta\bigl(e(x-u), e(v-x);q\bigr)}{\tau\bigl(q,q;q\bigr)_{\infty}
\theta\bigl(e(v-u);q\bigr)}\\
&\times\sum_{n=-\infty}^{\infty}\frac{1}{\bigl(1-e((v-x-n)/\tau)\bigr)\bigl(e((x+n-u)/\tau)-1\bigr)}.
\end{split}
\end{equation}
\end{lemma}
\begin{proof}
Set $\widetilde{q}=e(-1/\tau)$.
The bilateral sum
\begin{equation*}
\begin{split}
f(x)&=\sum_{n=-\infty}^{\infty}\frac{1}{\bigl(1-e((v-x-n)/\tau)\bigr)\bigl(e((x+n-u)/\tau)-1\bigr)}\\
&=\frac{1}{(1-e((v-x)/\tau))(e((x-u)/\tau)-1)}\,{}_2\psi_2
\left(\begin{matrix} e((v-x)/\tau), e((u-x)/\tau)\\ \widetilde{q}e((v-x)/\tau), \widetilde{q}e((u-x)/\tau)
\end{matrix}; \widetilde{q}, \widetilde{q}\right)
\end{split}
\end{equation*}
defines an elliptic function on $\mathbb{C}/\Lambda$,
with possible poles at most simple and located at $u+\Lambda$ and
at $v+\Lambda$. Hence there exists a $\eta\in\mathbb{C}$ (unique up to $\Lambda$-translates) and a constant
$C_\eta\in\mathbb{C}$ such that
\[f(x)=C_\eta\frac{\theta\bigl(e(u+v-\eta-x), e(x-\eta);q\bigr)}{\theta\bigl(e(x-u), e(v-x);q\bigr)}.
\]
We now compute the residue of $f$ at $u$ in two different ways:
\[\underset{x=u}{\hbox{Res}}(f)=\frac{\tau}{2\pi i}\frac{1}{\bigl(1-e((v-u)/\tau)\bigr)}
\]
from the bilateral series expression of $f$, and
\[\underset{x=u}{\hbox{Res}}(f)=-\frac{C_\eta}{2\pi i}\frac{\theta\bigl(e(u-\eta), e(v-\eta);q\bigr)}
{\bigl(q,q;q\bigr)_{\infty}\theta\bigl(e(v-u);q\bigr)}
\]
from the expression of $f$ as a quotient of theta-functions. Combining both identities yields an explicit expression
of the constant $C_\eta$ in terms of $\eta$, resulting in the formula
\[f(x)=\frac{\tau \bigl(q,q;q\bigr)_{\infty}\theta\bigl(e(v-u);q\bigr)}{\bigl(e((v-u)/\tau)-1\bigr)}
\frac{\theta\bigl(e(u+v-\eta-x), e(x-\eta);q\bigr)}{\theta\bigl(e(x-u), e(v-x), e(u-\eta), e(v-\eta);q\bigr)}
\]
for $f$. Rewriting this identity yields the desired result.
\end{proof}


\subsection{Trigonometric hypergeometric integrals with $E_6$ symmetries}\label{E6}

We consider trigonometric degenerations of $S_e(t)$ ($t\in\mathcal{H}_{pq}$) along root vectors $\alpha\in R(E_8)$ lying in the
$W(E_7)=W(E_8)_{\delta}$-orbit
\begin{equation}\label{orbit}
\mathcal{O}:=W(E_7)\bigl(\alpha_{18}^+\bigr)=\{\alpha_{jk}^+, \gamma_{jk} \,\, | \,\, 1\leq j<k\leq 8\},
\end{equation}
cf. Section \ref{notations}.
The degenerations relate to the explicit bijection
\begin{equation}\label{parametertbij}
\mathcal{G}_0\overset{\sim}{\longrightarrow}\mathcal{G}_{\log(pq)},\qquad
(u_1,\ldots,u_8)\mapsto (u_1,\ldots,u_8)+\log(pq)\alpha
\end{equation}
on the parameter spaces (in logarithmic form) of the associated integrals.
We obtain two different trigonometric degenerations, depending on whether we degenerate
along an orbit vector of the form $\alpha=\alpha_{jk}^+$, or of the form $\gamma_{jk}$.

Specifically, we consider the trigonometric degenerations $S_t(t)$ respectively $U_t(t)$ ($t\in\mathcal{H}_1$) of $S_e(t)$
($t\in\mathcal{H}_{pq}$) along the orbit vector $\alpha_{18}^+$ and $\gamma_{18}$ respectively.
The orbit vector $\alpha_{18}^+$ (respectively $\gamma_{18}$) is the additional simple root turning
the basis $\Delta_1$ (respectively $\Delta_2$) of $R(E_7)$ into the basis $\overline{\Delta}_1$ (respectively $\overline{\Delta}_2$) of $R(E_8)$,
see Section \ref{notations}. The induced symmetry group of $S_t(t)$ ($t\in\mathcal{H}_1$)
is the isotropy subgroup $W(E_7)_{\alpha_{18}^+}$ of $W(E_7)$, while the induced symmetry group of $U_t(t)$
($t\in\mathcal{H}_1$) is $W(E_7)_{\gamma_{18}}$.
It follows from the analysis in Section \ref{notations} that
$W(E_7)_{\alpha_{18}^+}=W(E_7)_{\gamma_{18}}$ is a maximal, standard parabolic subgroup
of $W(E_7)$ with respect to both bases $\Delta_1$ and $\Delta_2$, isomorphic to the Weyl group $W(E_6)$ of type $E_6$,
with corresponding simple roots $\Delta_1^\prime=\Delta_1\setminus \{\alpha_{21}^-\}$ and $\Delta_2^\prime=\Delta_2\setminus \{\alpha_{87}^-\}$, and
with corresponding Dynkin diagrams
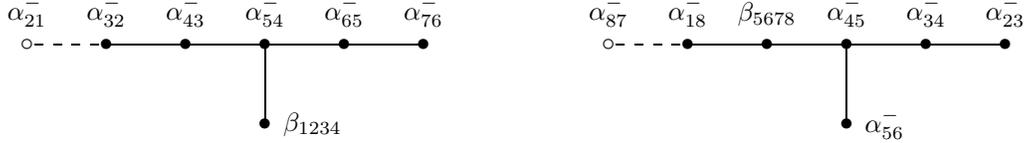
\begin{center}
\begin{picture}(370,46)

\put(0,0){\begin{picture}(150,46)  
\put(0,30){\bulo}
\put(0,37){\makebox(0,0)[b]{$\alpha_{21}^-$}}
\put(30,30){\bul}
\put(30,37){\makebox(0,0)[b]{$\alpha_{32}^-$}}
\put(60,30){\bul}
\put(60,37){\makebox(0,0)[b]{$\alpha_{43}^-$}}
\put(90,30){\bul}
\put(90,37){\makebox(0,0)[b]{$\alpha_{54}^-$}}
\put(120,30){\bul}
\put(120,37){\makebox(0,0)[b]{$\alpha_{65}^-$}}
\put(150,30){\bul}
\put(150,37){\makebox(0,0)[b]{$\alpha_{76}^-$}}
\put(90,0){\bul}
\put(97,0){\makebox(0,0)[l]{$\beta_{1234}$}}
\put(30,30){\line(1,0){120}}
\put(90,30){\line(0,-1){30}}
\multiput(3,30)(7,0){4}{\line(1,0){3}}
\end{picture}}

\put(220,0){\begin{picture}(150,46)  
\put(0,30){\bulo}
\put(0,37){\makebox(0,0)[b]{$\alpha_{87}^-$}}
\put(30,30){\bul}
\put(30,37){\makebox(0,0)[b]{$\alpha_{18}^-$}}
\put(60,30){\bul}
\put(60,37){\makebox(0,0)[b]{$\beta_{5678}$}}
\put(90,30){\bul}
\put(90,37){\makebox(0,0)[b]{$\alpha_{45}^-$}}
\put(120,30){\bul}
\put(120,37){\makebox(0,0)[b]{$\alpha_{34}^-$}}
\put(150,30){\bul}
\put(150,37){\makebox(0,0)[b]{$\alpha_{23}^-$}}
\put(90,0){\bul}
\put(97,0){\makebox(0,0)[l]{$\alpha_{56}^-$}}
\put(30,30){\line(1,0){120}}
\put(90,30){\line(0,-1){30}}
\multiput(3,30)(7,0){4}{\line(1,0){3}}
\end{picture}}

\end{picture}
\end{center}
Observe that $\alpha_{18}^-$ and $\alpha_{76}^-$ are the highest roots
of the standard parabolic root system $R(E_6)$ of type $E_6$ in $R(E_7)$ corresponding to the bases $\Delta_1^\prime$ and $\Delta_2^\prime$
respectively. {}From now on we write
\[
W(E_6):=W(E_7)_{\alpha_{18}^+}=W(E_7)_{\gamma_{18}}.
\]

We first introduce the trigonometric hypergeometric integrals $S_t(t)$ and $U_t(t)$
($t\in\mathcal{H}_1$) explicitly. Their integrands are defined by
\begin{equation*}
\begin{split}
I_t(t;z) &= \frac{(z^{\pm 2},t_1^{-1} z^{\pm 1},t_8^{-1} z^ {\pm 1};q)_\infty}
{\prod_{j=2}^7 (t_j z^{\pm 1};q)_\infty},\\
J_t^\mu(t;z) &= 2\frac{\theta(t_1t_8/\mu z,
z/\mu;q)}{\theta(t_1/\mu, t_8/\mu;q)}\left(1-\frac{z^2}{q}\right)
\prod_{i=2}^7\frac{\bigl(z/t_i;q\bigr)_{\infty}}{\bigl(t_iz;q\bigr)_{\infty}}\prod_{j=1,8}\frac{1}{\bigl(t_jz/q,
t_j/z;q\bigr)_{\infty}},
\end{split}
\end{equation*}
where $t=(t_1,\ldots,t_8)\in \bigl(\mathbb{C}^\times\bigr)^8$.
For generic $t=(t_1,\ldots,t_8)\in\mathbb{C}^8$ satisfying $\prod_{j=1}^8t_j=1$
and generic $\mu\in\mathbb{C}$
we now define the resulting trigonometric hypergeometric integrals
as
\[
S_t(t) = \int_{\mathcal{C}} I_t(t;z) \frac{dz}{2\pi i z},\qquad
U_t^\mu(t)=\int_{\mathcal{C}^\prime} J_t^\mu(t;q)\frac{dz}{2\pi i z}
\]
where $\mathcal{C}$ (respectively $\mathcal{C}^\prime$)
is a deformation of the positively oriented unit circle $\mathbb{T}$ including the pole sequences
$t_jq^{\mathbb{Z}_{\geq 0}}$ ($j=2,\ldots,7$) of $I_t(t;z)$ and excluding their reciprocals
(respectively including the pole sequences $t_jq^{\mathbb{Z}_{\geq 0}}$ ($j=1,8$) of $J_t^\mu(t;z)$ and excluding
the pole sequences $t_j^{-1}q^{\mathbb{Z}_{\leq 1}}$ ($j=1,8$) and $t_i^{-1}q^{\mathbb{Z}_{\leq 0}}$ ($i=2,\ldots,7$)).
As in the elliptic and hyperbolic cases, one observes that
$S_t(t)$ (respectively $U_t^\mu(t)$) admits a unique meromorphic
extension to the parameter domain $\{t\in\mathbb{C}^8 \, | \, \prod_{j=1}^8t_j=1\}$
(respectively $\{(\mu,t)\in \mathbb{C}^\times\times\mathbb{C}^8\, | \, \prod_{j=1}^8t_j=1\}$).
We call $S_t(t)$ the trigonometric hypergeometric function.
\begin{lemma}\label{muindep}
The integral $U_t^\mu(t)$ is independent of
$\mu\in\mathbb{C}^\times$.
\end{lemma}
\begin{proof}
There are several different, elementary arguments
to prove the lemma, we give here the argument based on Liouville's Theorem.
Note that $U_t^{q\mu}(t)=U_t^\mu(t)$, and that the possible poles of
$\mu\mapsto U_t^\mu(t)$ are at $t_jq^{\mathbb{Z}}$ ($j=1,8$). Without loss of generality we
assume the generic conditions on the parameters $t\in\mathbb{C}^8$ ($\prod_{j=1}^8t_j=1$)
such that $U_t^\mu(t)$ admits the integral representation as above, and such that $t_1\not\in
t_8q^{\mathbb{Z}}$. The latter condition ensures that the possible poles $t_1q^{\mathbb{Z}}, t_8q^{\mathbb{Z}}$
of $\mu\mapsto U_t^\mu(t)$ are at most simple. But the residue of $U_t^\mu(t)$ at $\mu=t_j$ ($j=1,8$) is zero,
since it is an integral over a deformation $\mathcal{C}^\prime$ of $\mathbb{T}$ whose integrand
is analytic  within the integration contour $\mathcal{C}^\prime$ and vanishes at the origin. Hence $\mathbb{C}^\times \ni\mu\mapsto U_t^\mu(t)$ is
bounded and analytic, hence constant by Liouville's Theorem.
\end{proof}
In view of Lemma \ref{muindep}, we omit the $\mu$-dependence in the notation for $U_t^\mu(t)$.
Since $I_t(-t;z)=I_t(t;-z)$ and
$J_t^\mu(-t;z)=J_t^{-\mu}(t;-z)$, we may and will view $S_t$ and $U_t$
as meromorphic function on $\mathcal{H}_1$.

By choosing a special value of $\mu$, we are able to derive another, "unfolded" integral representation of $U_t(t)$ as follows.
Let $\mathbb{H}_+$ by the upper half plane in $\mathbb{C}$. Choose $\tau\in\mathbb{H}_+$ such that $q=e(\tau)$,
where $e(x)$ is a shorthand notation for $\exp(2\pi ix)$. Recall the surjective map
$\psi_0: \mathcal{G}_0\rightarrow \mathcal{H}_1$ from Section \ref{notations}.
\begin{cor}\label{unfold}
For generic parameters $u\in\mathcal{G}_0$ we have
\begin{equation*}
\begin{split}
U_t(\psi_0(2\pi iu))&=\frac{2}{\tau\bigl(q,q;q\bigr)_{\infty}}\frac{\bigl(e((u_8-u_1)/\tau)-1\bigr)}{e(u_1)
\theta(e(u_8-u_1);q)}\\
&\times\int_{\mathcal{L}}\left\{\left(1-\frac{e(2x)}{q}\right)\prod_{j=2}^7\frac{\bigl(e(x-u_j);q\bigr)_{\infty}}
{\bigl(e(x+u_j);q\bigr)_{\infty}}\frac{\bigl(qe(x-u_1),
qe(x-u_8);q\bigr)_{\infty}}{\bigl(q^{-1}e(x+u_1),
q^{-1}e(x+u_8);q\bigr)_{\infty}}\right.\\
&\qquad\qquad\qquad\quad\qquad\qquad\qquad\qquad\times\left.\frac{e(x)}{\bigl(1-e((u_8-x)/\tau)\bigr)
\bigl(1-e((x-u_1)/\tau)\bigr)}\right\}dx
\end{split}
\end{equation*}
where the integration contour $\mathcal{L}$ is some translate $\xi+\mathbb{R}$
\textup{(}$\xi\in i\mathbb{R}$\textup{)} of the real line with a finite number of indentations, such that
$\mathcal{C}$ separates the pole sequences $-u_1+\mathbb{Z}+\mathbb{Z}_{\leq 1}\tau$, $-u_8+\mathbb{Z}+\mathbb{Z}_{\leq 1}\tau$ and
$-u_j+\mathbb{Z}+\mathbb{Z}_{\leq 0}\tau$ \textup{(}$j=2,\ldots,7$\textup{)} of the integrand from the pole sequences
$u_1+\mathbb{Z}_{\geq 0}\tau$ and $u_8+\mathbb{Z}_{\geq 0}\tau$.
\end{cor}
\begin{remark}
Note that always $\xi\not=0$ in Corollary \ref{unfold}.
Due to the balancing condition $\sum_{j=1}^8u_j=0$, there are no parameter choices for which
$\mathcal{L}=\mathbb{R}$ can be taken as integration contour.
This is a reflection of the fact that there are no parameters $t\in\mathcal{H}_1$ such that the unit circle
$\mathbb{T}$ can be chosen as integration cycle in the original integral representation
$U_t(t)=\int_{\mathcal{C}^\prime}J_t^\mu(t;z)\frac{dz}{2\pi iz}$ of $U_t(t)$.
\end{remark}

\begin{proof}
In the integral expression
\[U_t(\psi_c(2\pi iu))=\int_{\mathcal{C}^\prime}J_t^\mu(\psi_0(2\pi iu);z)\frac{dz}{2\pi iz},
\]
we change the integration variable to $z=e(x)$, take $\mu=e(\eta(u_1,u_8))$, and we use Lemma
\ref{unfoldlemma} to rewrite the quotient of theta-functions in the integrand as a bilateral sum.
Changing the integration over the intented line segment with the bilateral sum using Fubini's Theorem,
we can rewrite the resulting expression as a single integral over a
noncompact integration cycle $\mathcal{L}$. This leads directly to the desired result.
\end{proof}
In the following lemma we show that $U_t(t)$ can be expressed as a
sum of two nonterminating very-well-poised ${}_{10}\phi_9$ series.

\begin{lemma}\label{series}
As meromorphic functions in $t\in \mathcal{H}_1$, we have
\begin{equation*}
U_t(t)=\frac{2}{\bigl(q,t_1^2, t_1t_8/q, t_8/t_1;q\bigr)_{\infty}}\prod_{j=2}^7\frac{\bigl(t_1/t_j;q\bigr)_{\infty}}
{\bigl(t_1t_j;q\bigr)_{\infty}}\,{}_{10}W_9\bigl(\frac{t_1^2}{q};t_1t_2, t_1t_3,\ldots, t_1t_7,
\frac{t_1t_8}{q};q,q\bigr)+
\bigl(t_1\leftrightarrow t_8\bigr).
\end{equation*}
\end{lemma}
\begin{proof}
For generic $t\in \mathcal{H}_1$ we shrink the contour $\mathcal{C}^\prime$ in the integral representation
of $U_t^\mu(t)=\int_{\mathcal{C}^\prime}J_t^\mu(t;z)\frac{dz}{2\pi iz}$
to the origin while picking up the residues at the pole sequences $t_1q^{\mathbb{Z}_{\geq 0}}$ and $t_8q^{\mathbb{Z}_{\geq 0}}$
of the integrand $J_t^\mu(t;z)$. The resulting sum of residues can be directly rewritten as a sum of two very-well-poised ${}_{10}\phi_9$
series, leading to the desired identity (cf. the general residue techniques in \cite[\S 4.10]{GenR}).
\end{proof}
\begin{remark}\label{Phi}
Lemma \ref{series} yields that $U_t(t)$ is, up to an explicit rescaling factor,
an integral form of the particular sum $\Phi$ of two very-well-poised ${}_{10}\phi_9$ series
as e.g. studied in \cite{GuptaMasson} and \cite{LiJ} (see \cite[(1.8)]{GuptaMasson}, \cite[(9c)]{LiJ}).
Note furthermore that the explicit $\mu$-dependent quotient of theta-functions in the integrand of
$U_t^\mu(t)$ has the effect that it balances the very-well-poised ${}_{10}\phi_9$ series
when picking up the residues of $J_t^\mu(t;z)$ at the two pole sequences $t_1q^{\mathbb{Z}_{\geq 0}}$ and $t_8q^{\mathbb{Z}_{\geq 0}}$.
\end{remark}

In the following proposition we show that $S_t$ (respectively $U_t$) is the degeneration of $S_e$
along the root vector $\alpha_{18}^+$ (respectively $\gamma_{18}$).

\begin{prop}\label{elldegtrig}
Let $t=(t_1,\ldots,t_8)\in\mathbb{C}^8$ be generic parameters satisfying the balancing condition
$\prod_{j=1}^8t_j=1$. Then
\begin{equation}\label{degenerationformulas}
\begin{split}
S_t(t)&=\lim_{p\rightarrow 0}S_e(pqt_1,t_2,\ldots,t_7, pqt_8),\\
U_t(t)&=\lim_{p\rightarrow 0}\theta\bigl(t_1t_8/pq;q\bigr)S_e\bigl((pq)^{-\frac{1}{2}}t_1, (pq)^{\frac{1}{2}}t_2,\ldots, (pq)^{\frac{1}{2}}t_7,
(pq)^{-\frac{1}{2}}t_8).
\end{split}
\end{equation}
\end{prop}
\begin{proof}
For the degeneration to $S_t(t)$ we use that
\[I_e(pqt_1,t_2,\ldots,t_7,pqt_8;z)=\frac{\prod_{j=2}^7\Gamma_e(t_jz^{\pm 1};p,q)}{\Gamma_e(z^{\pm 2}, t_1^{-1}z^{\pm 1}, t_8^{-1}z^{\pm 1};p,q)}
\]
in view of the reflection equation for $\Gamma_e$, which (pointwise) tends to $I_t(t;z)$ as $p\rightarrow 0$ in view of \eqref{elltrigred}.
A standard application of Lebesgue's dominated convergence theorem leads to the limit of the associated integrals.

The degeneration to $U_t(t)$ is more involved, since one needs to use a nontrivial symmetry argument to cancel some unwanted sequences of
poles of $I_e(t;z)$. To ease the notations we set
\[t_p=\bigl((pq)^{-\frac{1}{2}}t_1,(pq)^{\frac{1}{2}}t_2,\ldots,(pq)^{\frac{1}{2}}t_7, (pq)^{-\frac{1}{2}}t_8\bigr)
\]
 and we denote
\[
Q(z)=\frac{\theta\bigl((pq)^{-\frac{1}{2}}t_1t_8/\mu z, (pq)^{-\frac{1}{2}}t_1z, (pq)^{-\frac{1}{2}}t_8z, (pq)^{\frac{1}{2}}z/\mu;q\bigr)}
{\theta\bigl(z^2;q\bigr)}.
\]
By \eqref{eqellfund}, we have the identity
\[Q(z)+Q(z^{-1})=\theta\bigl(t_1t_8/pq, t_1/\mu, t_8\mu;q\bigr).
\]
Since the integrand $I_e(\tau_p;z)$ is invariant under $z\mapsto z^{-1}$, we can consequently write
\[\theta\bigl(t_1t_8/pq;q\bigr)S_e(t_p;z)=\frac{2}{\theta\bigl(t_1/\mu, t_8/\mu;q\bigr)}\int_{\mathcal{C}}Q(z)I_e(t_p;z)\frac{dz}{2\pi iz},
\]
with $\mathcal{C}$ a deformation of the positively oriented unit circle $\mathbb{T}$ seperating the downward pole sequences of the integrand
from the upward pole sequences. Taking $(pq)^{\frac{1}{2}}z$ as a new integration variable and using the functional equation
and reflection equation of $\Gamma_e$, we obtain the integral representation
\begin{equation}\label{almost}
\begin{split}
&\theta\bigl(t_1t_8/pq;q\bigr)S_e(t_p;z)=\\
&\qquad =2\int_{\mathcal{C}}\frac{\theta\bigl(t_1t_8/\mu z, z/\mu;q\bigr)}
{\theta\bigl(t_1/\mu, t_8/\mu;q\bigr)}\theta\bigl(z^2/q;p\bigr)\prod_{j=2}^7\frac{\Gamma_e(t_jz;p,q)}{\Gamma_e(z/t_j;p,q)}
\prod_{j=1,8}\Gamma_e\bigl(t_jz/q, t_j/z;p,q\bigr)\frac{dz}{2\pi iz},
\end{split}
\end{equation}
where $\mathcal{C}$ is a deformation of the positively oriented unit circle $\mathbb{T}$ which
includes the pole sequences
$t_1p^{\mathbb{Z}_{\geq 0}}q^{\mathbb{Z}_{\geq 0}}, t_8p^{\mathbb{Z}_{\geq 0}}q^{\mathbb{Z}_{\geq 0}}$ and
$t_jp^{\mathbb{Z}_{\geq 1}}q^{\mathbb{Z}_{\geq 1}}$ ($j=2,\ldots,7$), and which excludes the pole sequences
$t_1^{-1}p^{\mathbb{Z}_{\leq 0}}q^{\mathbb{Z}_{\leq 1}}$, $t_8^{-1}p^{\mathbb{Z}_{\leq 0}}q^{\mathbb{Z}_{\leq 1}}$ and
$t_j^{-1}p^{\mathbb{Z}_{\leq 0}}q^{\mathbb{Z}_{\leq 0}}$ ($j=2,\ldots,7$).
We can now take the limit $p\rightarrow 0$ in \eqref{almost} with $p$-independent, fixed integration contour $\mathcal{C}$,
leading to the desired limit relation
\[\lim_{p\rightarrow 0}\theta\bigl(t_1t_8/pq;q\bigr)S_e(t_p)=U_t^\mu(t).
\]
\end{proof}

\begin{remark}
Observe that Lemma \ref{series} and the proof of Proposition \ref{elldegtrig} entail independent proofs of
Lemma \ref{muindep}.
\end{remark}

By specializing the parameters $t\in\mathcal{H}_1$ in
Proposition \ref{elldegtrig} further, we arrive at trigonometric integrals which can be evaluated
by \eqref{eqellbeta}. The resulting trigonometric degenerations
lead immediately to the trigonometric Nassrallah-Rahman integral evaluation formula \cite[(6.4.1)]{GenR}
and Gasper's integral evaluation formula \cite[(4.11.4)]{GenR}:
\begin{cor}\label{evaltrig}
For generic $t=(t_1,\ldots,t_6)\in\mathbb{C}^6$ satisfying the balancing condition
$\prod_{j=1}^6t_j=1$ we have
\begin{equation*}
\begin{split}
&\int_{\mathcal{C}}\frac{\bigl(z^{\pm 2}, t_1^{-1}z^{\pm 1};q\bigr)_{\infty}}
{\prod_{j=2}^6\bigl(t_jz^{\pm 1};q\bigr)_{\infty}}\frac{dz}{2\pi iz}=\frac{2\prod_{j=2}^6\bigl(1/t_1t_j;q\bigr)_{\infty}}
{\bigl(q;q\bigr)_{\infty}\prod_{2\leq j<k\leq 6}\bigl(t_jt_k;q\bigr)_{\infty}},\\
&\int_{\mathcal{C}^\prime}\frac{\theta\bigl(t_5t_6/\mu z, z/\mu;q\bigr)}{\theta\bigl(t_5/\mu, t_6/\mu;q\bigr)}
\left(1-\frac{z^2}{q}\right)\prod_{j=1}^4\frac{\bigl(z/t_j;q\bigr)_{\infty}}{\bigl(t_jz;q\bigr)_{\infty}}
\prod_{k=5,6}\frac{1}{\bigl(t_kz/q, t_k/z;q\bigr)_{\infty}}\frac{dz}{2\pi iz}=\\
&\qquad\qquad\qquad\qquad\qquad\qquad\qquad=\frac{\prod_{1\leq j<k\leq 4}\bigl(1/t_jt_k;q\bigr)_{\infty}}{\bigl(q,t_5t_6/q;q\bigr)_{\infty}
\prod_{j=1}^4\prod_{k=5}^6\bigl(t_jt_k;q\bigr)_{\infty}},
\end{split}
\end{equation*}
where $\mathcal{C}$ is the deformation of $\mathbb{T}$ seperating the pole sequences
$t_jq^{\mathbb{Z}_{\geq 0}}$ \textup{(}$j=2,\ldots,6$\textup{)} of the integrand from their reciprocals,
and where $\mathcal{C}^\prime$ is the deformation of $\mathbb{T}$ seperating the pole
sequences $t_5q^{\mathbb{Z}_{\geq 0}}, t_6q^{\mathbb{Z}_{\geq 0}}$ of the integrand
from the pole sequences $t_j^{-1}q^{\mathbb{Z}_{\leq 0}}$ \textup{(}$j=1,\ldots,4$\textup{)},
$t_5^{-1}q^{\mathbb{Z}_{\leq 1}}$ and $t_6^{-1}q^{\mathbb{Z}_{\leq 1}}$.
\end{cor}
\begin{proof}
For the first integral evaluation, take $t\in\mathcal{H}_1$ and $t_7=t_8^{-1}$ in the degeneration from
$S_e$ to $S_t$, and use the elliptic Nassrallah-Rahman integral evaluation formula \eqref{eqellbeta}.

For the second integral evaluation, take $t\in\mathcal{H}_1$ with $t_6=t_7^{-1}$ in the degeneration from
$S_e$ to $U_t^\mu$ and again use \eqref{eqellbeta} to evaluate the elliptic integral. It leads to the second integral evaluation formula
with generic parameters $(t_1,t_2,t_3,t_4,t_5,t_8)\in\mathbb{C}^6$ satisfying $t_1\cdots t_5t_8=1$.
\end{proof}

The second integral in Corollary \ref{evaltrig} can be unfolded using Corollary \ref{unfold}.
We obtain for generic parameters $u\in\mathbb{C}^6$ satisfying $\sum_{j=1}^6u_j=0$,
\begin{equation*}
\begin{split}
\int_{\mathcal{L}}&\left\{\left(1-\frac{e(2x)}{q}\right)\prod_{j=1}^4\frac{\bigl(e(x-u_j);q\bigr)_{\infty}}
{\bigl(e(x+u_j);q\bigr)_{\infty}}\frac{\bigl(qe(x-u_5),
qe(x-u_6);q\bigr)_{\infty}}{\bigl(q^{-1}e(x+u_5),
q^{-1}e(x+u_6);q\bigr)_{\infty}}\right.\\
&\qquad\qquad\qquad\qquad\qquad\qquad\qquad\left.\times\frac{e(x)}{\bigl(1-e((u_6-x)/\tau)\bigr)
\bigl(1-e((x-u_5)/\tau)\bigr)}\right\}dx\\
&\qquad\qquad\qquad\qquad=\frac{\tau
t_5\theta\bigl(t_6/t_5;q\bigr)}{\bigl(e((u_6-u_5)/\tau)-1\bigr)}\frac{\bigl(q;q\bigr)_{\infty}\prod_{1\leq
j<k\leq 4}\bigl(1/t_jt_k;q\bigr)_{\infty}}{\bigl(t_5t_6/q\bigr)_{\infty}\prod_{j=1}^4
\prod_{k=5}^6\bigl(t_jt_k;q\bigr)_{\infty}},
\end{split}
\end{equation*}
where $\tau\in\mathbb{H}_+$ such that $q=e(\tau)$, where $t_j=e(u_j)$ ($j=1,\ldots,6$) and
where the integration contour $\mathcal{L}$ is some translate $\xi+\mathbb{R}$
\textup{(}$\xi\in i\mathbb{R}$\textup{)} of the real line with a finite number of indentations, such that
$\mathcal{C}$ separates the pole sequences $-u_5+\mathbb{Z}+\mathbb{Z}_{\leq 1}\tau$,
$-u_6+\mathbb{Z}+\mathbb{Z}_{\leq 1}\tau$ and
$-u_j+\mathbb{Z}+\mathbb{Z}_{\leq 0}\tau$ \textup{(}$j=1,\ldots,4$\textup{)} of the integrand from the pole sequences
$u_5+\mathbb{Z}_{\geq 0}\tau$ and $u_6+\mathbb{Z}_{\geq 0}\tau$.
This is Agarwal's identity \cite[(4.7.5)]{GenR}.

Furthermore, using Lemma \ref{series} the second integral in Corollary \ref{evaltrig} can be written as
a sum of two very-well-poised ${}_8\phi_7$-series. We obtain for
generic $t\in\mathbb{C}^6$ satisfying $\prod_{j=1}^6t_j=1$,
\begin{equation*}
\begin{split}
\frac{1}{\bigl(q,t_5^2, t_5t_6/q, t_6/t_5;q\bigr)_{\infty}}
\prod_{j=1}^4\frac{\bigl(t_5/t_j;q\bigr)_{\infty}}{\bigl(t_5t_j;q\bigr)_{\infty}}
\,&{}_8W_7\bigl(\frac{t_5^2}{q}; t_5t_1, t_5t_2, t_5t_3, t_5t_4,
\frac{t_5t_6}{q};q,q\bigr)+ (t_5\leftrightarrow t_6)\\
&=\frac{\prod_{1\leq j<k\leq 4}\bigl(1/t_jt_k;q\bigr)_{\infty}}
{\bigl(q,
t_5t_6/q;q\bigr)_{\infty}\prod_{j=1}^4\prod_{k=5}^6\bigl(t_jt_k;q\bigr)_{\infty}}
\end{split}
\end{equation*}
which is Bailey's summation formula \cite[(2.11.7)]{GenR} of the
sum of two very-well-poised ${}_8\phi_7$ series.

We can now compute the (nontrivial) $W(E_6)$-symmetries of the trigonometric hypergeometric integrals
$S_t$ and $U_t$ by taking limits of the corresponding symmetries on the elliptic level using
Proposition \ref{elldegtrig}. We prefer to give a derivation based on the trigonometric evaluation formulas
(see Corollary \ref{evaltrig}), in analogy to our approach in the elliptic and hyperbolic cases.
\begin{prop}\label{relation}
The trigonometric integrals $S_t(t)$ and $U_t(t)$ \textup{(}$t\in\mathcal{H}_1$\textup{)} are invariant under permutations of $(t_1,t_8)$
and of $(t_2,\ldots,t_7)$. Furthermore,
\begin{equation}\label{eqtrigsym}
\begin{split}
S_t(t)&= S_t(wt) \frac{\bigl(1/t_1t_2,1/t_1t_3,1/t_1t_4,1/t_8t_5,1/t_8t_6,1/t_8t_7;q\bigr)_\infty}
{\bigl(t_2t_3,t_2t_4,t_3t_4,t_5t_6,t_5t_7,t_6t_7;q\bigr)_\infty},\\
U_t(t)&= U_t(wt) \frac{\bigl(1/t_2t_3, 1/t_2t_4, 1/t_3t_4, 1/t_5t_6, 1/t_5t_7, 1/t_6t_7;q\bigr)_{\infty}}
{\bigl(t_1t_2, t_1t_3, t_1t_4, t_5t_8, t_6t_8, t_7t_8;q\bigr)_{\infty}}
\end{split}
\end{equation}
as meromorphic functions in $t\in\mathcal{H}_1$.
\end{prop}
\begin{proof}
In order to derive the $w$-symmetry of $S_t(t)$ we consider the double integral
\begin{align*}
\int_{\mathbb{T}^2} \frac{\bigl(z^{\pm 2},x^{\pm 2},t_1^{-1}z^{\pm 1}, st_8^{-1}x^{\pm 1};q\bigr)_{\infty}}
{\bigl(t_2z^{\pm 1},t_3z^{\pm 1},t_4z^{\pm 1}, sz^{\pm 1} x^{\pm 1},s^{-1}t_5x^{\pm 1},s^{-1}t_6x^{\pm 1},s^{-1}t_7x^{\pm 1};q\bigr)_\infty}
\frac{dz}{2\pi i z} \frac{dx}{2\pi i x}
\end{align*}
for parameters $(t_1,\ldots,t_8)\in \mathbb{C}^8$ satisfying $\prod_{j=1}^8t_j=1$,
where $s^2t_1t_2t_3t_4=1=s^{-2}t_5t_6t_7t_8$ and where we assume
the additional parameter restraints
\[|t_2|, |t_3|, |t_4|, |s|, |t_5/s|, |t_6/s|, |t_7/s|<1
\]
to ensure that the integration contour $\mathbb{T}$ separates the
downward sequences of poles of from the upward sequences.
The desired transformation then follows by either integrating the double integral
first to $x$, or first to $z$, using in each case the trigonometric Nassrallah-Rahman integral evaluation formula
(see Corollary \ref{evaltrig}).

The proof of the $w$-symmetry of $U_t(t)$ follows the same line of arguments. For $\epsilon>0$ we denote $\epsilon\mathbb{T}$
for the positively oriented circle in the complex plane with radius $\epsilon$ and centered at the origin.
The $w$-symmetry
\[U_t^{\mu/s}(t)=U_t^{\mu}(wt)\frac{\bigl(1/t_2t_3, 1/t_2t_4, 1/t_3t_4, 1/t_5t_6, 1/t_5t_7, 1/t_6t_7;q\bigr)_{\infty}}
{\bigl(t_1t_2, t_1t_3, t_1t_4, t_5t_8, t_6t_8, t_7t_8;q\bigr)_{\infty}}
\]
for $t\in\mathcal{H}_1$, where $s^2t_1t_2t_3t_4=1=s^{-2}t_5t_6t_7t_8$, by considering for
 $(t_1,\ldots,t_8)\in \mathbb{C}^8$ satisfying $\prod_{j=1}^8t_j=1$ the double integral
\begin{equation*}
\begin{split}
&\int_{\bigl(|qs|\mathbb{T}\bigr)^2}\left\{\frac{\theta\bigl(st_1t_8/\mu z, t_8z/\mu x, x/\mu;q\bigr)}
{\theta\bigl(st_8/\mu, st_1/\mu, t_8/s\mu;q\bigr)}\left(1-\frac{z^2}{q}\right)\left(1-\frac{x^2}{q}\right)\right.\\
&\;\; \times\left.\frac{\bigl(z/t_2, z/t_3, z/t_4, xz/s, sx/t_5, sx/t_6, sx/t_7;q\bigr)_{\infty}}
{\bigl(t_1z/q, t_1/z, t_2z, t_3z, t_4z, sxz/q, sz/x, sx/z, t_5x/s, t_6x/s, t_7x/s, t_8x/qs, t_8/sx;q\bigr)_{\infty}} \! \right\} \!
\frac{dz}{2\pi i z}\frac{dx}{2\pi ix}
\end{split}
\end{equation*}
with $s^2t_1t_2t_3t_4=1=s^{-2}t_5t_6t_7t_8$, where we assume the additional parameter restraints
\[0<|s|\ll |q^{\frac{1}{2}}|,\qquad |t_1|, |t_2^{-1}|, |t_3^{-1}|, |t_4^{-1}| <|qs|,\qquad
|t_5|, |t_6|, |t_7|<|q^{-1}|,\qquad |t_8|<|qs^2|
\]
to ensure a proper separation by the integration contours of the upward sequences of poles from the downward sequences.
Using the second trigonometric integral
evaluation formula of Corollary \ref{evaltrig} then yields the desired result for the restricted parameter
domain. Analyic continuation completes the proof.
\end{proof}

\begin{remark}
Rewriting $U_t(t)$ as a sum of two very-well-poised ${}_{10}\phi_9$ series (see Lemma \ref{series}
and Remark \ref{Phi}),
the $w$-symmetry of $U_t(t)$ becomes Bailey's four-term transformation formula \cite[(2.12.9)]{GenR},
see also \cite{GuptaMasson}.
The identification of the symmetry group of $U_t$ with the Weyl
group of type $E_6$ has been derived by different methods in \cite{LiJ}.
\end{remark}

Finally we relate the two trigonometric integrals $S_t$ and $U_t$. We can obtain the following
proposition as a degeneration of a particular $W(E_7)$-symmetry of $S_e$, but we prefer here
to give a direct proof using double integrals.
\begin{prop}\label{SU}
As meromorphic function in $t\in\mathcal{H}_1$, we have
\begin{align*}
S_t(t)
\frac{\prod_{2\leq j <k \leq 5} (t_jt_k;q)_\infty (t_6t_7;q)_\infty }
{(1/qt_1t_8,1/t_1t_6,1/t_1t_7,1/t_8t_6,1/t_8t_7;q)_\infty}
    =
U_t(t_6/s,st_2,st_3,st_4,st_5,t_1/s,t_8/s,t_7/s),
\end{align*}
where $t_2t_3t_4t_5s^2=1=t_1t_6t_7t_8/s^2$
\end{prop}
\begin{proof}
For $(t_1,\ldots,t_8)\in\mathbb{C}^8$ satisfying $\prod_{j=1}^8t_j=1$ we consider the double integral
\begin{align*}
\int_{z\in \eta\mathbb{T}}\int_{x\in\mathbb{T}}&
\frac{\theta(\mu z,t_6t_7\mu/s^2z)}{\theta(t_6\mu/s,t_7\mu/s)}
\left(1-\frac{z^2}{q}\right) \prod_{j=2}^5\frac{1}{(t_j x^{\pm 1};q)_\infty} \\ & \times
\frac{(x^{\pm 2},zx^{\pm 1}/s,sz/t_1,sz/t_8;q)_\infty}{(szx^{\pm 1},t_7z/sq,t_7/sz,t_6z/sq,t_6/sz,t_1z/s,t_8z/s;q)_\infty}
\frac{dx}{2\pi i x} \frac{dz}{2 \pi i z},
\end{align*}
with $s^2t_2t_3t_4t_5=1=s^{-2}t_1t_6t_7t_8$ and with $0<\eta<\min\bigl(|s^{-1}|,|q^{\frac{1}{2}}|\bigr)$, where we
assume the additional parameter restraints
\[|t_2|, |t_3|, |t_4|, |t_5|<1,\qquad |t_6|, |t_7|<\eta |s|,\qquad
|t_1|, |t_8|<\eta^{-1}|s|
\]
to ensure a proper separation by the integration contours of the upward sequences of poles from the downward sequences.
Using Corollary \ref{evaltrig}, we can first integrate over $x$ using the trigonometric Nassrallah-Rahman integral evaluation formula,
or first integrate over $z$ using the second integral evaluation formula of Corollary \ref{evaltrig}.
The resulting identity gives the desired result for restricted parameter values. Analytic continuation completes the proof.
\end{proof}

\begin{remark}\label{biorthtrig}
{\bf (i)} Combining Proposition \ref{SU} with Lemma \ref{series} we obtain an expression of $S_t(t)$ as a sum of
two very-well-poised ${}_{10}\phi_9$ series, which is originially due to Rahman \cite[(6.4.8)]{GenR}.

{\bf (ii)} For e.g. $t_1t_6=q^m$ ($m\in\mathbb{Z}_{\geq 0}$), it follows from {\bf (i)}  (see also \cite{Ra} and \cite[(6.4.10)]{GenR})
that the $S_t(t;p,q)$ essentially coincides with
the biorthogonal rational function of Rahman \cite{Ra}, which is explicitly given as a terminating very-well-poised ${}_{10}\phi_9$ series.
\end{remark}



\subsection{Contiguous relations}

The fundamental equation on this level equals
\begin{equation}\label{fundreltrig}
\frac{1}{y}(1-vx^{\pm 1})(1-yz^{\pm 1}) + \frac{1}{z}(1-vy^{\pm 1})(1-zx^{\pm 1}) + \frac{1}{x}(1-vz^{\pm 1})(1-xy^{\pm 1}) =0
\end{equation}
where $(1-ax^{\pm 1})=(1-ax)(1-ax^{-1})$. The fundamental relation \eqref{fundreltrig} is the $p=0$ reduction of
\eqref{eqellfund}.
In this section $\tau_{ij}=\tau_{ij}^{-\log(q)}$ acts as in the elliptic case by multiplying $t_i$ by $q$ and
dividing $t_j$ by $q$. Formula \eqref{fundreltrig} leads as in the elliptic case to the difference equation
\begin{equation}\label{eqtrigdiff1}
\frac{(1-t_5t_6^{\pm 1}/q)}{(1-t_4t_6^{\pm 1})}  S_t(\tau_{45}t)
+ \frac{(1-t_5t_4^{\pm 1}/q)}{(1-t_6t_4^{\pm 1})}  S_t(\tau_{65}t) = S_t(t),\qquad t\in\mathcal{H}_1.
\end{equation}
To obtain a second difference equation between trigonometric hypergeometric functions where two times
the same parameter is multiplied by $q$, we can mimick the approach in the elliptic case with the role of the
longest Weyl group element taken over by the element $u=ws_{35}s_{46}w\in W_(E_6)$.
Alternatively, one can rewrite the difference equation \eqref{eqellcont2} for $S_e$
in the form
\begin{equation*}
\begin{split}
&\frac{\theta\bigl(t_3/qt_4,1/t_1t_5, 1/t_8t_5, t_2t_5/q, t_5t_6/q, t_5t_7/q;p\bigr)}
{\theta\bigl(t_3/t_5;p\bigr)}S_e(\tau_{45}\widetilde{t})+ (t_3\leftrightarrow t_5)=\\
&\qquad\qquad\qquad\qquad\qquad\qquad\qquad\qquad=\theta\bigl(1/qt_1t_4, 1/qt_8t_4, t_2t_4, t_4t_6, t_4t_7;p\bigr)S_e(\widetilde{t})
\end{split}
\end{equation*}
where $t\in \mathcal{H}_1$ and $\widetilde{t}=(pqt_1,t_2,\ldots, t_7, pqt_8)$,
and degenerate it using Proposition \ref{elldegtrig}. We arrive at
\begin{equation}
\frac{(1-t_3/qt_4)}{(1-t_3/t_5)}\prod_{j=1,8}\frac{(1-1/t_5t_j)}{(1-1/qt_4t_j)}
\prod_{j=2,6,7} \frac{(1-t_5t_j/q)}{(1-t_4t_j)}
S_t(\tau_{45}t) + (t_3 \leftrightarrow t_5) = S_t(t),\qquad t\in\mathcal{H}_1
\end{equation}
Together these equations imply the following result.
\begin{prop}
We have
\begin{equation}\label{eqtrigcont}
A(t) S_t(\tau_{45}t) + (t_4 \leftrightarrow t_5) = B(t)S_t(t)
\end{equation}
as meromorphic functions in $t\in\mathcal{H}_1$,
where
\begin{align*}
A(t) & = - \frac{(1-\frac{1}{t_1t_5})(1-\frac{1}{t_8t_5}) \prod_{j=2,3,6,7} (1-\frac{t_5t_j}{q})}
{t_4 (1-\frac{t_4t_5}{q})(1-\frac{t_5}{qt_4})(1-\frac{t_4}{t_5})} \\
B(t) & = -\frac{(1-\frac{1}{qt_1t_6})(1-\frac{1}{qt_8t_6})(1-t_3t_6)(1-t_7t_6)(1-t_2t_6)}{t_6(1-\frac{t_4}{qt_6})(1-\frac{t_5}{qt_6})} \\
& \qquad + \frac{(1-\frac{t_6}{t_4})(1-t_6t_4)\prod_{j=1,8}(1-\frac{1}{t_jt_5})
\prod_{j=2,3,7}(1-\frac{t_jt_5}{q})}{t_6(1-\frac{t_5}{qt_6})(1-\frac{t_4t_5}{q})(1-\frac{t_5}{qt_4})(1-\frac{t_4}{t_5})}  \\
& \qquad + \frac{(1-\frac{t_6}{t_5})(1-t_6t_5)\prod_{j=1,8}(1-\frac{1}{t_jt_4})
\prod_{j=2,3,7}(1-\frac{t_jt_4}{q})}{t_6(1-\frac{t_4}{qt_6})(1-\frac{t_4t_5}{q})(1-\frac{t_4}{qt_5})(1-\frac{t_5}{t_4})}.
\end{align*}
\end{prop}
Despite the apparent asymmetric expression $B$ still satisfies $B(s_{67}t)=B(t)$.

The contiguous relation for the elliptic hypergeometric function $S_e$ with step-size $p$
can also be degenerated to the trigonometric level. A direct derivation is as follows.
By \eqref{eqellfund} we have
\[
\frac{\theta(t_8^{-1}t_7^{\pm 1};q)}{\theta(t_6t_7^{\pm 1};q)} I_t(t_1,t_2,\ldots,t_5,q t_8,t_7,t_6/q;z)
+ \frac{\theta(t_8^{-1}t_6^{\pm 1};q)}{\theta(t_7t_6^{\pm 1};q)} I_t(t_1,t_2,\ldots,t_6,q t_8,t_7/q;z)
=I_t(t;z).
\]
Integrating this equation we obtain
\begin{equation}\label{eqtrigs3}
\frac{\theta(t_8^{-1}t_7^{\pm 1};q)}{\theta(t_6t_7^{\pm 1};q)} S_t(t_1,t_2,\ldots,t_5,q t_8,t_7,t_6/q)
+ (t_6\leftrightarrow t_7) = S_t(t)
\end{equation}
as meromorphic functions in $t\in\mathcal{H}_1$, a three term transformation for $S_t$.
The three term transformation \cite[(6.5)]{GuptaMasson} is equivalent to the sum of two equations of this type
(in which the parameters are chosen such that two terms coincide and two other terms
cancel each other).

\begin{remark}
In \cite{LiJ} it is shown that there are essentially five different types of
three term transformations for $\Phi$ (see Remark \ref{Phi}), or equivalently of the integrals $U_t$ and
$S_t$. The different types arise from a careful analysis of the
three term transformations in terms of the $W(E_7)$-action on
$\mathcal{H}_1$. It is likely that all five different
types of three term transformations for $\Phi$ can be re-obtained
by degenerating contiguous relations for $S_e$ with step-size $p$
(similarly as the derivation of \eqref{eqtrigs3}): concretely, the five prototypes
are in one-to-one correspondence to the orbits of
\[\{(\alpha,\beta,\gamma)\in\mathcal{O}^3 \, | \, \alpha,\beta,\gamma \hbox{ are pair-wise different}\}
\]
under the diagonal action of $W(E_7)$, where $\mathcal{O}$ is the $W(E_7)$-orbit \eqref{orbit}.
\end{remark}


\subsection{Degenerations with $D_5$ symmetries}\label{D5}

In this section we consider degenerations of $S_t$ and $U_t$ with
symmetries with respect to the Weyl group of type $D_5$. Compared to the
analysis on the hyperbolic level, we introduce a trigonometric analog of the
Euler and Barnes' type integrals, as well as a third, new type of
integral arising as degeneration of $U_t$. We first introduce the degenerate integrals explicitly.

For generic $t=(t_1,\ldots,t_6)\in\bigl(\mathbb{C}^\times\bigr)^6$ we define the
trigonometric Euler integral as
\begin{equation}
E_t(t)=\int_{\mathcal{C}}\frac{\bigl(z^{\pm 2}, t_1^{-1}z^{\pm
1};q\bigr)_{\infty}}{\prod_{j=2}^6\bigl(t_jz^{\pm
1};q\bigr)_{\infty}}\frac{dz}{2\pi iz},
\end{equation}
where $\mathcal{C}$ is a deformation of the positively oriented
unit circle $\mathbb{T}$ seperating the decreasing pole sequences $t_jq^{\mathbb{Z}_{\geq 0}}$
($j=2,\ldots,6$) of the integrand from their reciprocals. We have
$E_t(-t)=E_t(t)$, and $E_t$ has a unique meromorphic extension to
$\bigl(\mathbb{C}^\times\bigr)^6$. The resulting meromorphic function on
$\bigl(\mathbb{C}^\times\bigr)^6/C_2$ is denoted also by $E_t$.

For generic $\mu\in\mathbb{C}^\times$ and generic $t=(t_1,\ldots,t_8)\in\mathbb{C}^8$ satisfying the
balancing condition $\prod_{j=1}^8t_j=1$ we define the
trigonometric Barnes integral as
\begin{equation}
B_t(t) = 2\int_{\mathcal{C}} \frac{\theta(t_2t_7/\mu z, z/\mu;q)}{\theta(t_2/\mu,t_7/\mu;q)}
\frac{(z/t_1,z/t_8;q)_{\infty}}{\prod_{j=3}^6 (t_jz;q)_\infty (t_2/z,t_7/z;q)_\infty}
\frac{dz}{2\pi i z},
\end{equation}
where $\mathcal{C}$ is a deformation of $\mathbb{T}$ seperating the decreasing pole
sequences $t_2q^{\mathbb{Z}_{\geq 0}}$ and $t_7q^{\mathbb{Z}_{\geq 0}}$ of the integrand
from the increasing pole sequences $t_j^{-1}q^{\mathbb{Z}_{\leq 0}}$ ($j=3,\ldots,6$).
Analogously to the analysis of the integral $U_t(t)$, we have that
the trigonometric Barnes integral $B_t(t)$ uniquely extends to a meromorphic function in
$\{(\mu,t)\in \mathbb{C}^\times\times \mathbb{C}^8 \, | \, \prod_{j=1}^8t_j=1\}$ which is
independent of $\mu$ (cf. Lemma \ref{muindep}). Furthermore, by a change of integration variable
we have $B_t(-t)=B_t(t)$, hence $B_t$ may (and will) be interpreted as meromorphic
function on $\mathcal{H}_1$.

Finally, for generic $t=(t_1,\ldots,t_6)\in\bigl(\mathbb{C}^\times\bigr)^6$ we
consider
\begin{equation}
V_t(t)=2\int_{\mathcal{C}}\frac{\theta\bigl(qt_2t_3t_4t_5t_6z;q\bigr)}
{\theta\bigl(qt_1t_2t_3t_4t_5t_6;q\bigr)}\left(1-\frac{z^2}{q}\right)
\prod_{j=2}^6\frac{\bigl(z/t_j;q\bigr)_{\infty}}{\bigl(t_jz;q\bigr)_{\infty}}
\frac{1}{\bigl(t_1z/q, t_1/z;q\bigr)_{\infty}}\frac{dz}{2\pi iz},
\end{equation}
where $\mathcal{C}$ is a deformation of $\mathbb{T}$ seperating
the decreasing pole sequence $t_1q^{\mathbb{Z}_{\geq 0}}$ of the
integrand from the remaining (increasing) pole sequences. As before, $V_t$
unique extends to a meromorphic function on $\bigl(\mathbb{C}^\times\bigr)^6/C_2$.

Similarly as for $U_t(t)$, the trigonometric Barnes integral $B_t(t)$ can be unfolded.
Recall that $q=e(\tau)$ with $\tau\in\mathbb{H}_+$, where $e(x)$ is a shorthand notation
for $\exp(2\pi ix)$.
\begin{lemma}\label{unfoldB}
For generic parameters $u\in\mathcal{G}_0$ we have
\begin{equation*}
\begin{split}
B_t(\psi_0(2\pi iu))&=\frac{2}{\tau\bigl(q,q;q\bigr)_{\infty}}\frac{\bigl(e((u_7-u_2)/\tau)-1\bigr)}{e(u_2)
\theta(e(u_7-u_2);q)}\\
&\times\int_{\mathcal{L}}\left\{\frac{\bigl(e(x-u_1), qe(x-u_2), qe(x-u_7),
e(x-u_8);q\bigr)_{\infty}}
{\bigl(e(x+u_3), e(x+u_4), e(x+u_5),
e(x+u_6);q\bigr)_{\infty}}\right.\\
&\left.\qquad\qquad\qquad\qquad\qquad\qquad\qquad\times
\frac{e(x)}{\bigl(1-e((u_7-x)/\tau)\bigr)\bigl(1-e((x-u_2)/\tau)\bigr)}\right\}dx,
\end{split}
\end{equation*}
where the integration contour $\mathcal{L}$ is some translate $\xi+\mathbb{R}$
\textup{(}$\xi\in i\mathbb{R}$\textup{)} of the real line with a finite number of indentations, such that
$\mathcal{C}$ separates the pole sequences $-u_j+\mathbb{Z}+\mathbb{Z}_{\leq 0}\tau$
\textup{(}$j=3,\ldots,6$\textup{)} of the integrand from the pole sequences
$u_2+\mathbb{Z}_{\geq 0}\tau$ and $u_7+\mathbb{Z}_{\geq 0}\tau$.
\end{lemma}
\begin{proof}
The proof is similar to the proof of Corollary \ref{unfold}.
\end{proof}
For $B_t(t)$ and $V_t(t)$ we have the following series expansions
in balanced ${}_4\phi_3$'s (respectively in a very-well-poised
${}_8\phi_7$).
\begin{lemma}\label{series2}
{\bf (a)} We have
\[B_t(t)=\frac{2\bigl(t_2/t_1, t_2/t_8;q\bigr)_{\infty}}
{\bigl(q,t_7/t_2, t_2t_3, t_2t_4, t_2t_5,
t_2t_6;q\bigr)_{\infty}}\,
{}_4\phi_3\left(\begin{matrix}
t_2t_3, t_2t_4, t_2t_5, t_2t_6\\
qt_2/t_7, t_2/t_1, t_2/t_8\end{matrix};q,q\right)+ (t_2\leftrightarrow
t_7)
\]
as meromorphic functions in $t\in\mathcal{H}_1$.\\
{\bf (b)} We have
\[V_t(t)=\frac{2}{\bigl(q,t_1^2;q\bigr)_{\infty}}\prod_{j=2}^6\frac{\bigl(t_1/t_j;q\bigr)_{\infty}}
{\bigl(t_1t_j;q\bigr)_{\infty}}\,
{}_8W_7\left(\frac{t_1^2}{q};t_1t_2, t_1t_3,\ldots,
t_1t_6;q,\frac{1}{t_1t_2t_3t_4t_5t_6}\right)
\]
as meromorphic functions in $t\in\bigl(\mathbb{C}^\times\bigr)^6/C_2$: $|t_1t_2t_3t_4t_5t_6|>1$.
\end{lemma}
\begin{proof}
This follows by a straightforward residue computation as in the
proof of Lemma \ref{series} (cf. also \cite[\S 4.10]{GenR}). For {\bf (a)} one picks up the
residues at the increasing pole sequences $t_2q^{\mathbb{Z}_{\geq 0}}$ and $t_7q^{\mathbb{Z}_{\geq
0}}$ of the integrand of $B_t(t)$; for {\bf (b)} one picks up the residues at the single increasing pole
sequence $t_1q^{\mathbb{Z}_{\geq 0}}$ of the integrand of
$V_t(t)$.
\end{proof}

\begin{prop}\label{trigellprop}
For generic $t\in\mathcal{H}_1$ we have
\begin{equation*}
\begin{split}
\lim_{u\rightarrow 0}
&S_t(t_1,\ldots,t_6,t_7u,t_8/u)=E_t(t_1,\ldots,t_6),\\
\lim_{u\rightarrow 0}&\bigl(t_2t_7/u;q\bigr)_{\infty}
S_t(t_1u^{-\frac{1}{2}}, t_2u^{-\frac{1}{2}},
t_3u^{\frac{1}{2}}, t_4u^{\frac{1}{2}}, t_5u^{\frac{1}{2}},
t_6u^{\frac{1}{2}}, t_7u^{-\frac{1}{2}},
t_8u^{-\frac{1}{2}})=B_t(t),\\
\lim_{u\rightarrow
0}&\bigl(t_1t_7/u;q\bigr)_{\infty}U_t(t_1,\ldots,t_6,t_7/u,
t_8u)=V_t(t_1,\ldots,t_6),\\
\lim_{u\rightarrow 0}& U_t(t_2u^{\frac{1}{2}}, t_1u^{\frac{1}{2}},
t_3u^{-\frac{1}{2}}, t_4u^{-\frac{1}{2}}, t_5u^{-\frac{1}{2}},
t_6u^{-\frac{1}{2}}, t_8u^{\frac{1}{2}},
t_7u^{\frac{1}{2}})=B_t(t).
\end{split}
\end{equation*}
\end{prop}
\begin{proof}
The first limit is direct. For the second limit, we
follow the same approach as in the proof of Proposition
\ref{elldegtrig}.
Define $Q(z)$ as
\[
Q(z) =  \frac{\theta(t_2zs^{-\frac12},t_7z s^{-\frac12},\mu zs^{\frac12},t_2t_7 \mu s^{-\frac12}z^{-1};q)}{\theta(z^2;q)}.
\]
Using \eqref{eqellfund} we obtain the equation
\[
Q(z)+ Q(z^{-1}) = \theta(t_2t_7/s,t_2\mu,t_7\mu;q),
\]
and hence, as in the proof of Proposition \ref{elldegtrig},
\[
(t_2t_7/u;q)_\infty S_t(t_u)
 =2 \frac{(t_2t_7/u;q)_\infty}{\theta(t_2t_7/u,t_2\mu,t_7\mu;q)}
 \int_{\mathcal{C}} I_t(t_u;z) Q(z) \frac{dz}{2\pi i z}
\]
for an appropriate contour $\mathcal{C}$,
where we use the abbreviated notation
\[t_u=(t_1 u^{-\frac12},
t_2u^{-\frac12},t_3u^{\frac12}, t_4 u^{\frac12},
t_5 u^{\frac12},t_6u^{\frac12},t_7u^{-\frac12},t_us^{-\frac12}).
\]
Taking $u^{-\frac{1}{2}}z$ as new integration variable we obtain
\begin{align*}
(ut_2t_7;q)_\infty S_t(t_u)
&= 2
\int_{\mathcal{C}} \frac{ \theta(\mu z ,t_2t_7\mu/z;q) }{ \theta(t_2\mu,t_7\mu;q)} \frac{(z/t_1,z/t_8;q)_\infty}
{\prod_{j=3}^6 (t_j z;q)_\infty (t_2/z,t_7/z ;q)_\infty }
 \\ & \qquad \qquad \qquad \times
(1-\frac{u}{z^2}) \frac{(ut_1^{-1}/z,ut_8^{-1}/z,qut_2^{-1}/z,qut_7^{-1}/z ;q)_\infty}
{(qu/t_2t_7;q)_\infty \prod_{j=3}^6 (t_ju/z;q)_\infty}
\frac{dz}{2\pi i z},
\end{align*}
where, for $u$ small enough, we take $\mathcal{C}$ to be a
$u$-independent deformation of $\mathbb{T}$ seperating the
decreasing pole sequences $t_2q^{\mathbb{Z}_{\geq 0}}, t_7q^{\mathbb{Z}_{\geq 0}}$
and $t_juq^{\mathbb{Z}_{\geq 0}}$ ($j=3,\ldots,6$) of the
integrand from the decreasing pole sequences
$t_j^{-1}q^{\mathbb{Z}_{\leq 0}}$ ($j=3,\ldots,6$). The limit $u\rightarrow 0$
can be taken in the resulting integral, leading to the desired
result.

To prove the third limit, we set $\mu=q/t_1t_7t_8$ in the integral
expression of
$U_t(t)=\int_{\mathcal{C}^\prime}J_t^\mu(t;z)\frac{dz}{2\pi iz}$ to remove the
contribution $\bigl(t_7z;q\bigr)_{\infty}$ in the denominator of
the integrand:
\[U_t(t)=2\int_{\mathcal{C}^\prime}\frac{\theta(t_1t_7t_8/z;q)}{\theta(t_1t_7,t_7t_8;q)}
\left(1-\frac{z^2}{q}\right)\prod_{j=2}^6\frac{\bigl(z/t_j;q\bigr)_{\infty}}{\bigl(t_jz;q\bigr)_{\infty}}
\frac{\bigl(z/t_7, q/t_7z;q\bigr)_{\infty}}{\bigl(t_1z/q, t_1/z,
t_8z/q, t_8/z;q\bigr)_{\infty}}\frac{dz}{2\pi iz}.
\]
In the resulting integral the desired limit can be taken
directly, leading to the desired result.

For the fourth limit, one easily verifies that
\[B_t(t)=\lim_{u\rightarrow 0} U_t^{\mu u^{\frac{1}{2}}}(t_2u^{\frac{1}{2}}, t_1u^{\frac{1}{2}},
t_3u^{-\frac{1}{2}}, t_4u^{-\frac{1}{2}}, t_5u^{-\frac{1}{2}},
t_6u^{-\frac{1}{2}}, t_8u^{\frac{1}{2}},
t_7u^{\frac{1}{2}})
\]
for generic $t\in\mathcal{H}_1$ after changing integration variable $z$ to $zu^{\frac{1}{2}}$
on the right hand side.
\end{proof}

Proposition \ref{trigellprop} and Corollary
\ref{evaltrig} immediately lead to the following three trigonometric integral
evaluations (of which the first is the well known Askey-Wilson
integral evaluation \cite[(6.1.4)]{GenR}).
\begin{cor}\label{trigevalcor}
For generic parameters $t=(t_1,t_2,t_3,t_4)\in\mathbb{C}^4$ we
have
\begin{equation*}
\begin{split}
\int_{\mathcal{C}}\frac{\bigl(z^{\pm
2};q\bigr)_{\infty}}{\prod_{j=1}^4\bigl(t_jz^{\pm
1};q\bigr)_{\infty}}\frac{dz}{2\pi
iz}&=\frac{2\bigl(t_1t_2t_3t_4;q\bigr)_{\infty}}{\bigl(q;q\bigr)_{\infty}\prod_{1\leq
j<k\leq 4}\bigl(t_jt_k;q\bigr)_{\infty}},\\
\int_{\mathcal{C}^\prime}\frac{\theta\bigl(qt_2t_3t_4z;q\bigr)}{\theta\bigl(qt_1t_2t_3t_4;q\bigr)}
\left(1-\frac{z^2}{q}\right)&\frac{\bigl(z/t_2, z/t_3,
z/t_4;q\bigr)_{\infty}}{\bigl(t_1/z, t_1z/q, t_2z, t_3z, t_4z;q\bigr)_{\infty}}\frac{dz}{2\pi
iz}\\
&=\frac{\bigl(qt_1t_2t_3t_4, 1/t_2t_3, 1/t_2t_4,
1/t_3t_4;q\bigr)_{\infty}}{\bigl(q, t_1t_2, t_1t_3, t_1t_4;q\bigr)_{\infty}}
\end{split}
\end{equation*}
with $\mathcal{C}$ \textup{(}respectively $\mathcal{C}^\prime$\textup{)} a deformation of $\mathbb{T}$
separating the sequences $t_jq^{\mathbb{Z}_{\geq 0}}$ \textup{(}$j=1,\ldots,4$\textup{)} from their reciprocals
\textup{(}respectively separating $t_1q^{\mathbb{Z}_{\geq 0}}$
from $t_1^{-1}q^{\mathbb{Z}_{\leq 1}}, t_2^{-1}q^{\mathbb{Z}_{\leq 0}}, t_3^{-1}q^{\mathbb{Z}_{\leq 0}}$
and $t_4^{-1}q^{\mathbb{Z}_{\leq 0}}$\textup{)}.

For generic $\mu\in\mathbb{C}^\times$ and $t\in\mathbb{C}^6$ satisfying $\prod_{j=1}^6t_j=1$
we have
\[\int_{\mathcal{C}}\frac{\theta\bigl(t_1t_5/\mu z,
z/\mu;q\bigr)}{\theta\bigl(t_1/\mu, t_5/\mu;q\bigr)}
\frac{\bigl(z/t_6;q\bigr)_{\infty}}{\bigl(t_1/z, t_2z, t_3z, t_4z,
t_5/z;q\bigr)_{\infty}}\frac{dz}{2\pi iz}=
\frac{1}{\bigl(q;q\bigr)_{\infty}}\prod_{j=2}^4\frac{\bigl(1/t_jt_6;q\bigr)_{\infty}}
{\bigl(t_1t_j, t_jt_5;q\bigr)_{\infty}},
\]
with $\mathcal{C}$ a deformation of $\mathbb{T}$ separating the
pole sequences $t_1q^{\mathbb{Z}_{\geq 0}}, t_5q^{\mathbb{Z}_{\geq 0}}$
from $t_2^{-1}q^{\mathbb{Z}_{\leq 0}}, t_3^{-1}q^{\mathbb{Z}_{\leq 0}}$ and $t_4^{-1}q^{\mathbb{Z}_{\leq
0}}$.
\end{cor}
\begin{proof}
Specializing the degeneration from $S_t$ to $E_t$ in Proposition
\ref{trigellprop} to generic parameters
$t\in\mathcal{H}_1$ under the additional condition $t_1t_2=1$ and using the trigonometric
Nassrallah-Rahman integral evaluation (Corollary \ref{evaltrig})
leads to the Askey-Wilson integral evaluation with corresponding
parameters $(t_3,t_4,t_5,t_6)$.

Similarly, specializing the degeneration
from $U_t$ to $V_t$ (respectively $S_t$ to $B_t$) to generic parameters
$t\in\mathcal{H}_1$ under the additional condition $t_2t_3=1$ (respectively
$t_1t_3=1$) and using the
Nassrallah-Rahman integral evaluation we obtain the second
(respectively third) integral evaluation with parameters
$(t_1,t_4,t_5,t_6)$ (respectively $(t_2,t_4,t_5,
t_6,t_7,t_8)$).
\end{proof}
Various well-known identities are direct consequences of Corollary
\ref{trigevalcor}. Firstly, analogous to the unfolding of the integrals
$U_t$ and $V_t$ (see Corollary \ref{unfold} and
Lemma \ref{unfoldB}), the left hand side of the third integral evaluation
can be unfolded. We obtain for
generic $u\in\mathbb{C}^6$ with $\sum_{j=1}^6u_j=0$,
\begin{equation*}
\begin{split}
\int_{\mathcal{L}}\frac{\bigl(qe(x-u_1), qe(x-u_5),
e(x-u_6);q\bigr)_{\infty}}{\bigl(e(x+u_2), e(x+u_3),
e(x+u_4);q\bigr)_{\infty}}&\frac{e(x)}{\bigl(1-e((u_5-x)/\tau)\bigr)\bigl(1-e((x-u_1)/\tau)\bigr)}\,dx\\
=\frac{\tau t_1\theta\bigl(t_5/t_1;q\bigr)}{\bigl(e((u_5-u_1)/\tau)-1\bigr)}
&\frac{\bigl(q,1/t_2t_6, 1/t_3t_6,
1/t_4t_6;q\bigr)_{\infty}}
{\bigl(t_1t_2, t_1t_3, t_1t_4, t_2t_5, t_3t_5,
t_4t_5;q\bigr)_{\infty}}
\end{split}
\end{equation*}
where $\tau\in\mathbb{H}_+$ such that $q=e(\tau)$, where $t_j=e(u_j)$ ($j=1,\ldots,6$) and
where the integration contour $\mathcal{L}$ is some translate $\xi+\mathbb{R}$
\textup{(}$\xi\in i\mathbb{R}$\textup{)} of the real line with a finite number of indentations such that
$\mathcal{C}$ separates the pole sequences $-u_j+\mathbb{Z}+\mathbb{Z}_{\leq 0}\tau$
\textup{(}$j=2,3,4$\textup{)} of the integrand from the pole sequences
$u_1+\mathbb{Z}_{\geq 0}\tau$ and $u_5+\mathbb{Z}_{\geq 0}\tau$. This integral identity
is Agarwal's \cite[(4.4.6)]{GenR} trigonometric analogue of
Barnes' second lemma.

The left hand side of the second integral evaluation
in Corollary \ref{trigevalcor} can be rewritten as a unilateral
sum by picking up the
residues at $t_1q^{\mathbb{Z}_{\geq 0}}$, cf. Lemma \ref{series}.
The resulting identity is
\begin{equation*}
{}_6\phi_5\left(\begin{matrix}
q^{-1}t_1^2, q^{\frac{1}{2}}t_1, -q^{\frac{1}{2}}t_1, t_1t_2,
t_1t_3, t_1t_4\\
q^{-\frac{1}{2}}t_1, -q^{-\frac{1}{2}}t_1, t_1/t_2, t_1/t_3,
t_1/t_4\end{matrix}; q,\frac{1}{t_1t_2t_3t_4}\right)=
\frac{\bigl(t_1^2, 1/t_2t_3, 1/t_2t_4, 1/t_3t_4;q\bigr)_{\infty}}
{\bigl(1/t_1t_2t_3t_4, t_1/t_2, t_1/t_3,
t_1/t_4;q\bigr)_{\infty}},
\end{equation*}
for generic $t\in\mathbb{C}^4$ satisfying $|t_1t_2t_3t_4|>1$,
which is the ${}_6\phi_5$ summation formula
\cite[(2.7.1)]{GenR}.

For generic $t\in\mathbb{C}^6$ satisfying $\prod_{j=1}^6t_j=1$
the left hand side of the third integral evaluation in Corollary \ref{trigevalcor}
can be written as a sum of two unilateral series by picking up the poles
of the integrand at the decreasing sequences $t_1q^{\mathbb{Z}_{\geq 0}}$ and
$t_5q^{\mathbb{Z}_{\geq 0}}$ of poles of the integrand. The resulting
identity is
\begin{equation*}
\frac{\bigl(t_1/t_6;q\bigr)_{\infty}}{\bigl(t_5/t_1, t_1t_2,
t_1t_3, t_1t_4;q\bigr)_{\infty}}\,{}_3\phi_2\left(\begin{matrix}
t_1t_2, t_1t_3, t_1t_4\\ qt_1/t_5, t_1/t_6\end{matrix};q,q\right)
+\bigl(t_1\leftrightarrow t_5\bigr)=
\prod_{j=2}^4
\frac{\bigl(1/t_jt_6;q\bigr)_{\infty}}
{\bigl(t_1t_j, t_jt_5;q\bigr)_{\infty}}
\end{equation*}
for generic $t\in\mathbb{C}^6$ satisfying $\prod_{j=1}^6t_j=1$,
which is the nonterminating version \cite[(2.10.12)]{GenR}
of Saalsch{\"u}tz formula.

We now return to the three trigonometric hypergeometric integrals
$E_t$, $B_t$ and $V_t$.
Recall that the symmetry group of $S_t$ and $U_t$ is the subgroup
$W(E_6)=W(E_7)_{\alpha_{18}^+}=W(E_7)_{\gamma_{18}}$,
which is a maximal standard parabolic subgroup of $W(E_7)$
with respect to both bases $\Delta_1$ and $\Delta_2$ of $R(E_7)$ (see Section
\ref{notations}), with corresponding sub-bases
$\Delta_1^\prime=\Delta_1\setminus \{\alpha_{12}^-\}$ and
$\Delta_2^\prime=\Delta_2\setminus \{\alpha_{87}^-\}$
respectively.
The four limits of Proposition \ref{trigellprop} now imply
that the trigonometric integrals $E_t$, $B_t$ and $V_t$ have
symmetry groups $W(E_6)_{\alpha_{78}^-}$ or
$W(E_6)_{\beta_{1278}}$.
The stabilizer subgroup $W(E_6)_{\alpha_{78}^-}$ is a standard maximal parabolic subgroup of $W(E_6)$
with respect to both bases $\Delta_1^\prime$ or $\Delta_2^\prime$,
with corresponding sub-basis
\[\Delta(D_5)=\Delta_1^\prime\setminus
\{\alpha_{76}^-\}=\Delta_2^\prime\setminus \{\alpha_{18}^-\}
\]
and with corresponding Dynkin sub-diagrams
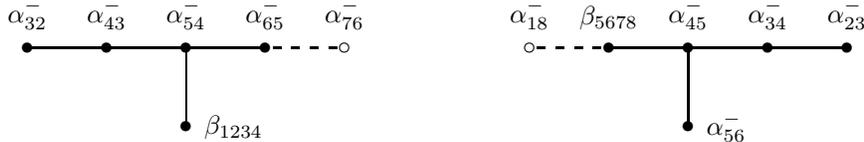
\begin{center}

\begin{picture}(310,46)

\put(0,0){\begin{picture}(120,46) 
\put(0,30){\bul}
\put(0,37){\makebox(0,0)[b]{$\alpha_{32}^-$}}
\put(30,30){\bul}
\put(30,37){\makebox(0,0)[b]{$\alpha_{43}^-$}}
\put(60,30){\bul}
\put(60,37){\makebox(0,0)[b]{$\alpha_{54}^-$}}
\put(90,30){\bul}
\put(90,37){\makebox(0,0)[b]{$\alpha_{65}^-$}}
\put(120,30){\bulo}
\put(120,37){\makebox(0,0)[b]{$\alpha_{76}^-$}}
\put(60,0){\bul}
\put(67,0){\makebox(0,0)[l]{$\beta_{1234}$}}
\put(0,30){\line(1,0){90}}
\put(60,30){\line(0,-1){30}}
\multiput(93,30)(7,0){4}{\line(1,0){3}}
\end{picture}}

\put(190,0){\begin{picture}(120,46) 
\put(0,30){\bulo}
\put(0,37){\makebox(0,0)[b]{$\alpha_{18}^-$}}
\put(30,30){\bul}
\put(30,37){\makebox(0,0)[b]{$\beta_{5678}$}}
\put(60,30){\bul}
\put(60,37){\makebox(0,0)[b]{$\alpha_{45}^-$}}
\put(90,30){\bul}
\put(90,37){\makebox(0,0)[b]{$\alpha_{34}^-$}}
\put(120,30){\bul}
\put(120,37){\makebox(0,0)[b]{$\alpha_{23}^-$}}
\put(60,0){\bul}
\put(67,0){\makebox(0,0)[l]{$\alpha_{56}^-$}}
\put(30,30){\line(1,0){90}}
\put(60,30){\line(0,-1){30}}
\multiput(3,30)(7,0){4}{\line(1,0){3}}
\end{picture}}

\end{picture}
\end{center}
respectively. Similarly, $W(E_6)_{\beta_{1278}}$ is a standard maximal
parabolic subgroup of $W(E_6)$ with respect to the basis
$\Delta_2^\prime$,
with corresponding sub-basis
\[\Delta^\prime(D_5)=\Delta_2^\prime\setminus \{\alpha_{23}^-\}
\]
and with corresponding Dynkin sub-diagram
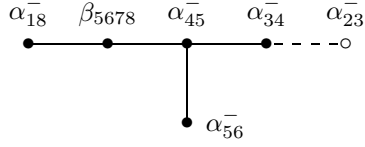
\begin{center}

\begin{picture}(120,46) 
\put(0,30){\bul}
\put(0,37){\makebox(0,0)[b]{$\alpha_{18}^-$}}
\put(30,30){\bul}
\put(30,37){\makebox(0,0)[b]{$\beta_{5678}$}}
\put(60,30){\bul}
\put(60,37){\makebox(0,0)[b]{$\alpha_{45}^-$}}
\put(90,30){\bul}
\put(90,37){\makebox(0,0)[b]{$\alpha_{34}^-$}}
\put(120,30){\bulo}
\put(120,37){\makebox(0,0)[b]{$\alpha_{23}^-$}}
\put(60,0){\bul}
\put(67,0){\makebox(0,0)[l]{$\alpha_{56}^-$}}
\put(0,30){\line(1,0){90}}
\put(60,30){\line(0,-1){30}}
\multiput(93,30)(7,0){4}{\line(1,0){3}}
\end{picture}
\end{center}
We write
\[W(D_5)=W(E_6)_{\alpha_{78}^-},\qquad W^\prime(D_5)=W(E_6)_{\beta_{1278}}\]
for the corresponding isotropy group, which are both isomorphic to
the Weyl group of type $D_5$.

The isotropy group $W(D_5)$ acts on
$\bigl(\mathbb{C}^\times\bigr)^6/C_2$: the simple reflections corresponding to
roots of the form $\alpha_{ij}^-\in\Delta(D_5)$ act by permuting the $i$th and $j$th coordinate,
while $w$ acts by
\[w(\pm t)=\pm\bigl(st_1,st_2,st_3,st_4,t_5/s,t_6/s\bigr),\qquad s^2=1/t_1t_2t_3t_4.
\]
With this action, the degenerations
to $E_t$ and $V_t$ in Proposition \ref{trigellprop} are
$W(D_5)$-equivariant in an obvious sense.

We can now directly compute the $W(D_5)$-symmetries of the trigonometric integrals
$E_t$ and $V_t$, as well as $W^\prime(D_5)$-symmetries of $B_t$, by taking limits of the corresponding
symmetries for $S_t$ and $U_t$ using Proposition
\ref{trigellprop}. This yields the following result.

\begin{prop}\label{EVBsymmetry}
{\bf a)} The trigonometric hypergeometric integrals $E_t(t)$ and $V_t(t)$
\textup{(}$t\in \bigl(\mathbb{C}^\times\bigr)^6\!/C_2$\textup{)} are invariant under
permutations of $(t_2,\ldots,t_6)$. Furthermore,
\begin{align*}
E_t(t)&=E_t(wt)\frac{\left(1/t_1t_2,1/t_1t_3,1/t_1t_4,t_1t_2t_3t_4t_5t_6;q\right)_{\infty}}
{\left(t_2t_3,t_2t_4,t_3t_4,t_5t_6;q\right)_\infty}, \\
V_t(t)&=V_t(wt)\frac{\left(1/t_2t_3, 1/t_2t_4, 1/t_3t_4,
1/t_5t_6;q\right)_{\infty}}
{\left(t_1t_2, t_1t_3, t_1t_4,
1/t_1t_2t_3t_4t_5t_6;q\right)_{\infty}}
\end{align*}
as meromorphic functions in $t\in \bigl(\mathbb{C}^\times\bigr)^6/C_2$.\\
{\bf b)} The trigonometric Barnes integral $B_t(t)$
\textup{(}$t\in\mathcal{H}_1$\textup{)} is invariant under
permutations of the pairs $(t_1,t_8)$, $(t_2,t_7)$ and of $(t_3,t_4,t_5,t_6)$.
Furthermore,
\[
B_t(t) = B_t(wt) \frac{\bigl(1/t_1t_3,1/t_1t_4,1/t_5t_8,1/t_6t_8;q\bigr)_\infty}
{\bigl(t_2t_4,t_2t_3,t_5t_7,t_6t_7;q\bigr)_\infty}
\]
as meromorphic functions in $t\in\mathcal{H}_1$.
\end{prop}
\begin{remark}
The $w$-symmetry of $V_t$, rewritten in series form using Lemma
\ref{series2}, gives the transformation formula
\cite[(2.10.1)]{GenR} for very-well-poised ${}_8\phi_7$ basic
hypergeometric series.
\end{remark}
Similarly as in the hyperbolic theory, the $w$-symmetry of $E_t$
generalizes to the following integral transformation formula for
the trigonometric Euler integral $E_t$.
\begin{prop}\label{inteqEt}
For $t\in\bigl(\mathbb{C}^\times\bigr)^6$ and $s\in\mathbb{C}^\times$
satisfying
\[|t_2|, |t_3|, |t_4|, |s|, |t_5/s|, |t_6/s|<1
\]
we have
\[\int_{\mathbb{T}}E_t(t_1,t_2,t_3,t_4,sx,sx^{-1})\frac{\bigl(x^{\pm
2};q\bigr)_{\infty}}{\bigl(t_5x^{\pm 1}/s, t_6x^{\pm
1}/s;q\bigr)_{\infty}}\,\frac{dx}{2\pi ix}=
\frac{2\bigl(t_5t_6;q\bigr)_{\infty}}{\bigl(q,s^2,t_5t_6/s^2;q\bigr)_{\infty}}\,
E_t(t_1,\ldots,t_6).
\]
\end{prop}
\begin{proof}
The proof is similar to the hyperbolic case (see Proposition
\ref{inteqE}), now using the double integral
\[\int_{\mathbb{T}^2}\frac{\bigl(z^{\pm 2}, z^{\pm 1}/t_1,
x^{\pm 2};q\bigr)_{\infty}}
{\bigl(sz^{\pm 1}x^{\pm
1};q\bigr)_{\infty}\prod_{j=2}^4\bigl(t_jz^{\pm
1};q\bigr)_{\infty}\prod_{k=5}^6\bigl(t_kx^{\pm
1}/s;q\bigr)_{\infty}}\,\frac{dz}{2\pi iz}\,\frac{dx}{2\pi ix}.
\]
\end{proof}
Specializing $s^2=1/t_1t_2t_3t_4$ in Proposition \ref{inteqEt} and
using the trigonometric Nassrallah-Rahman integral (see Corollary
\ref{trigevalcor}), we re-obtain the $w$-symmetry of $E_t$
(see Proposition \ref{EVBsymmetry}).

The three trigonometric integrals $E_t$, $B_t$ and $V_t$ are
interconnected as follows.

\begin{prop}\label{interconnectionprop}
We have
\begin{equation*}
\begin{split}
B_t(t)&=\frac{\bigl(1/t_1t_6, 1/t_8t_6;q\bigr)_{\infty}}
{\bigl(t_2t_3, t_2t_4, t_2t_5, t_6t_7;q\bigr)_{\infty}}\,
V_t\bigl(t_7/s, t_3s, t_4s, t_5s, t_1/s, t_8/s\bigr),\\
&=\frac{\bigl(1/t_1t_3, 1/t_1t_4, 1/t_1t_5,
1/t_6t_8;q\bigr)_{\infty}}{\bigl(t_2t_6,
t_7t_6;q\bigr)_{\infty}}\,E_t\bigl(t_8/v, t_7/v, t_3v, t_4v, t_5v,
t_2/v\bigr),
\end{split}
\end{equation*}
as meromorphic functions in $t\in\mathcal{H}_1$, where
$s^2=t_1t_6t_7t_8=1/t_2t_3t_4t_5$ and $v^2=t_2t_6t_7t_8=1/t_1t_3t_4t_5$.
\end{prop}
\begin{proof}
This follows by combining Proposition \ref{SU} and Proposition \ref{trigellprop}.
Concretely, to relate $B_t$ and $V_t$ one computes for
generic $t\in\mathcal{H}_1$ and with $s^2=1/t_2t_3t_4t_5$,
\begin{equation*}
\begin{split}
B_t(t)&=\lim_{u\rightarrow 0}\bigl(t_2t_7/u;q\bigr)_{\infty}\,S_t\bigl(t_1u^{-\frac{1}{2}},
t_2u^{-\frac{1}{2}}, t_3u^{\frac{1}{2}}, t_4u^{\frac{1}{2}},
t_5u^{\frac{1}{2}}, t_6u^{\frac{1}{2}}, t_7u^{-\frac{1}{2}},
t_8u^{-\frac{1}{2}}\bigr)\\
&=\frac{\bigl(1/t_1t_6, 1/t_8t_6;q\bigr)_{\infty}}
{\bigl(t_2t_3, t_2t_4, t_2t_5, t_6t_7;q\bigr)_{\infty}}
\lim_{u\rightarrow 0}\bigl(t_2t_7/u;q\bigr)_{\infty}\,
U_t\bigl(t_7/s, t_3s, t_4s, t_5s, t_1/s, t_8/s, t_2s/u,
t_6u/s\bigr)\\
&=\frac{\bigl(1/t_1t_6, 1/t_8t_6;q\bigr)_{\infty}}
{\bigl(t_2t_3, t_2t_4, t_2t_5, t_6t_7;q\bigr)_{\infty}}\,
V_t\bigl(t_7/s, t_3s, t_4s, t_5s, t_1/s, t_8/s\bigr),
\end{split}
\end{equation*}
where the first and third equality follows from Proposition \ref{trigellprop}
and the second equality follows from Proposition \ref{SU}.
To relate $B_t$ and $E_t$, we first note that Proposition \ref{SU}
is equivalent to the identity
\begin{equation}\label{inverted}
U_t(t)=\frac{\bigl(1/t_6t_7;q\bigr)_{\infty}\prod_{2\leq j<k\leq
5}\bigl(1/t_jt_k;q\bigr)_{\infty}}{\bigl(t_1t_8/q, t_1t_6, t_8t_6, t_1t_7,
t_8t_7;q\bigr)_{\infty}}\,
S_t\bigl(t_6/s, t_2s, t_3s, t_4s, t_5s, t_1/s, t_8/s, t_7/s\bigr),
\end{equation}
where $t\in\mathcal{H}_1$ and $s^2=1/t_2t_3t_4t_5$. For generic $t\in\mathcal{H}_1$
and with $v^2=1/t_1t_3t_4t_5$ we then compute
\begin{equation*}
\begin{split}
B_t(t)&=\lim_{u\rightarrow 0}U_t\bigl(t_2u^{\frac{1}{2}}, t_1u^{\frac{1}{2}},
t_3u^{-\frac{1}{2}}, t_4u^{-\frac{1}{2}}, t_5u^{-\frac{1}{2}},
t_6u^{-\frac{1}{2}}, t_8u^{\frac{1}{2}},
t_7u^{\frac{1}{2}}\bigr)\\
&=\frac{\bigl(1/t_1t_3, 1/t_1t_4, 1/t_1t_5,
1/t_6t_8;q\bigr)_{\infty}}
{\bigl(t_2t_6, t_7t_6;q\bigr)_{\infty}}\lim_{u\rightarrow 0}
S_t\bigl(t_8/v, t_7/v, t_3v, t_4v, t_5v, t_2/v, t_1vu,
t_6/vu\bigr)\\
&=\frac{\bigl(1/t_1t_3, 1/t_1t_4, 1/t_1t_5,
1/t_6t_8;q\bigr)_{\infty}}{\bigl(t_2t_6,
t_7t_6;q\bigr)_{\infty}}\,E_t\bigl(t_8/v, t_7/v, t_3v, t_4v, t_5v,
t_2/v\bigr),
\end{split}
\end{equation*}
where the first and third equality follows from Proposition \ref{trigellprop}
and the second equality follows from \eqref{inverted}.
\end{proof}
\begin{remark}
{\bf a)} Combining the interconnection between $E_t$ and $V_t$
from Proposition \ref{interconnectionprop} with the expression of
$V_t$ as very-well-poised ${}_8\phi_7$ series from Lemma \ref{series2}
yields the Nassrallah-Rahman integral representation \cite[(6.3.7)]{GenR}.

{\bf b)} Similarly, combining the interconnection between $B_t$ and $V_t$
from Proposition \ref{interconnectionprop} with their series
expressions from Lemma \ref{series2} yields the expression \cite[(2.10.10)]{GenR}
of a very-well-poised ${}_8\phi_7$ series as a sum of two balanced ${}_4\phi_3$
series.
\end{remark}

Degenerating the contiguous relations of $S_t$ using Proposition \ref{trigellprop}
leads directly to contiguous relations for $E_t$, $B_t$ and $V_t$. For instance,
we obtain

\begin{prop}\label{lemtrigdiffE}
We have
\[
A(t) E_t(t_1,t_2,t_3,qt_4,t_5/q,t_6) + (t_4 \leftrightarrow t_5) = B(t) E_t(t)
\]
as meromorphic functions in $t\in\bigl(\mathbb{C}^\times\bigr)^6/C_2$, where
\begin{align*}
A(t) & = -\frac{(1-\frac{1}{t_1t_5})(1-\frac{t_2t_5}{q})(1-\frac{t_3t_5}{q})(1-\frac{t_6t_5}{q})}
{t_4(1-\frac{t_4t_5}{q})(1-\frac{t_5}{qt_4})(1-\frac{t_4}{t_5})}, \\
B(t) &= \frac{q}{t_1t_4t_5} - t_2t_3t_6 +\frac{qt_4}{t_5}A(t)+\frac{qt_5}{t_4}A(s_{45}t).
\end{align*}
\end{prop}
\begin{proof}
Substitute $t=(t_1,t_2,t_3,t_4,t_5,t_7 u,t_6,t_8/u)$ with $\prod_{j=1}^8t_j=1$
in \eqref{eqtrigcont} and take the limit $u\to 0$.
\end{proof}
For later purposes, we also formulate the corresponding result for $B_t(t)$.
We substitute $t=\bigl(t_1u^{-\frac{1}{2}}, t_7u^{-\frac{1}{2}}, t_3u^{\frac{1}{2}}, t_4u^{\frac{1}{2}}, t_5u^{\frac{1}{2}},
t_2u^{-\frac{1}{2}}, t_6u^{\frac{1}{2}}, t_8u^{-\frac{1}{2}}\bigr)$ for generic $t\in\mathbb{C}^8$
satisfying $\prod_{j=1}^8t_j=1$ in \eqref{eqtrigcont}, multiply the resulting equation by $u^{\frac{1}{2}}\bigl(t_2t_7/u;q\bigr)_{\infty}$,
and take the limit $u\to 0$. We arrive at
\begin{equation}\label{BAW}
\alpha(t)B_t(t_1,t_2,t_3,qt_4,t_5/q,t_6,t_7,t_8)+ (t_4\leftrightarrow t_5)=\beta(t)B_t(t),\qquad t\in\mathcal{H}_1,
\end{equation}
where
\begin{equation*}
\begin{split}
\alpha(t)&=-\frac{\bigl(1-\frac{1}{t_1t_5}\bigr)\bigl(1-\frac{t_2t_5}{q}\bigr)\bigl(1-\frac{t_5t_7}{q}\bigr)\bigl(1-\frac{1}{t_5t_8}\bigr)}
{t_4\bigl(1-\frac{t_4}{t_5}\bigr)\bigl(1-\frac{t_5}{qt_4}\bigr)},\\
\beta(t)&=t_7\bigl(1-t_2t_3\bigr)\bigl(1-t_2t_6\bigr)
+\frac{\bigl(1-t_2t_4\bigr)}{\bigl(1-\frac{t_2t_5}{q}\bigr)}\alpha(t)+
\frac{\bigl(1-\frac{t_2t_4}{q}\bigr)}{\bigl(1-t_2t_5\bigr)}\alpha(s_{45}t).
\end{split}
\end{equation*}

The degeneration of \eqref{eqtrigs3} yields Bailey's \cite[(2.11.1)]{GenR}
three term transformation formula for very-well-poised ${}_8\phi_7$'s:
\begin{prop}
We have
\[
\frac{\theta(t_1^{-1}t_2^{\pm 1};q)}{\theta(t_3t_2^{\pm1};q)} E_t(t_3/q,qt_1,t_2,t_4,t_5,t_6)
+ (t_2 \leftrightarrow t_3) =E_t(t).
\]
\end{prop}
\begin{proof}
Consider \eqref{eqtrigs3} with $t_1$ and $t_8$, $t_2$ and $t_7$, and $t_3$ and $t_6$ interchanged.
Subsequently substitute the parameters $(t_1,t_2,t_3,t_4,t_5,t_6,t_7 u,t_8/u)$ with
$\prod_{j=1}^8t_j=1$ and take the limit $u\to 0$.
\end{proof}


\subsection{The Askey-Wilson function}

In this subsection we relate the trigonometric hypergeometric integrals with $D_5$ symmetry to
the nonpolynomial eigenfunction of the Askey-Wilson second order difference operator, known
as the Askey-Wilson function. The Askey-Wilson function is the trigonometric analog of Ruijsenaars'
$R$-function, and is closely related to harmonic analysis on the quantum $\hbox{SU}(1,1)$ group.

As for the $R$-function, we introduce the Askey-Wilson function in terms of the trigonometric
Barnes integral $B_t$.
Besides the usual Askey-Wilson parameters we also use logarithmic variables
in order to make the connection to the $R$-function more transparent. We write the base $q\in\mathbb{C}^\times$
with $|q|<1$ as $q=e(\omega_1/\omega_2)$ with $\tau=\omega_1/\omega_2\in\mathbb{H}_+$
and $e(x)=\exp(2\pi ix)$ as before.
{}From the previous subsection it follows that the
parameter space of $B_t(t)$ is $\mathcal{H}_1/\mathbb{C}^\times e(\beta_{1278})$.
In logarithmic coordinates, this relates to $\mathcal{G}_0/\mathbb{C}\beta_{1278}$.
We identify $\mathcal{G}_0/\mathbb{C}\beta_{1278}$ with $\mathbb{C}^6$
by assigning to the six-tuple $(\gamma,x,\lambda)=(\gamma_0,\gamma_1,\gamma_2,\gamma_3,\lambda,x)$
the class in $\mathcal{G}_0/\mathbb{C}\beta_{1278}$ represented by
$u=(u_1,\ldots,u_8)\in\mathcal{G}_0$ with
\begin{equation}\label{reparametrizelog}
\begin{split}
u_1&=\bigl(-\gamma_0-\gamma_1-2\omega\bigr)/\omega_2,\quad\,\,\,\,
u_2=0,\\
u_3&=\bigl(\hat{\gamma}_0+\omega-i\lambda\bigr)/\omega_2,\qquad\quad
u_4=\bigl(\gamma_0+\omega-ix\bigr)/\omega_2,\\
u_5&=\bigl(\gamma_0+\omega+ix\bigr)/\omega_2,\qquad\quad
u_6=\bigl(\hat{\gamma}_0+\omega+i\lambda\bigr)/\omega_2,\\
u_7&=\bigl(-\gamma_0-\gamma_3\bigr)/\omega_2,\qquad\qquad
u_8=\bigl(-\gamma_0-\gamma_2-2\omega\bigr)/\omega_2,
\end{split}
\end{equation}
where $\omega=\frac{1}{2}(\omega_1+\omega_2)$ as before.
We define the corresponding six-tuple of Askey-Wilson parameters
$(a,b,c,d,\mu,z)$ by
\begin{equation}\label{AWpar}
\begin{split}
(a,b,c,d)&=\bigl(e((\gamma_0+\omega)/\omega_2),e((\gamma_1+\omega)/\omega_2),
e((\gamma_2+\omega)/\omega_2),
e(\gamma_3+\omega)/\omega_2)\bigr),\\
(\mu,z)&=\bigl(e(-i\lambda/\omega_2), e(-ix/\omega_2)\bigr).
\end{split}
\end{equation}
The four-tuple $(a,b,c,d)$ represents the four parameter freedom in the
Askey-Wilson theory, while $z$ (respectively $\mu$) plays the role of
geometric (respectively spectral) parameter. Furthermore, we define
the dual Askey-Wilson parameters by
\[\bigl(\widetilde{a},\widetilde{b},\widetilde{c},\widetilde{d})=
\bigl(e((\hat{\gamma}_0+\omega)/\omega_2),
e((\hat{\gamma}_1+\omega)/\omega_2),
e((\hat{\gamma}_2+\omega)/\omega_2),
e((\hat{\gamma}_3+\omega)/\omega_2)\bigr),
\]
with $\hat{\gamma}$ the dual parameters defined by
\eqref{eqhypdefdualpar}. We furthermore associate to the
logarithmic parameters $u\in\mathcal{G}_0$ (see
\eqref{reparametrizelog}) the parameters $t=\psi_0(2\pi
iu)\in\mathcal{H}_1$, so that
\begin{equation}\label{reptrig}
t=\psi_0(2\pi iu)=
\bigl(1/ab,1,\widetilde{a}\mu, az, a/z, \widetilde{a}/\mu, q/ad, 1/ac\bigr)
\end{equation}

We define the Askey-Wilson function $\phi(\gamma;x,\lambda)=
\phi\bigl(\gamma;x,\lambda;\omega_1,\omega_2)$ by
\begin{equation}\label{phidef}
\phi(\gamma;x,\lambda)=\frac{\bigl(q,t_2t_3,t_2t_4,t_2t_5,t_2t_6;q\bigr)_{\infty}}{2}B_t(t)
\end{equation}
with $t=\psi_0(2\pi iu)\in\mathcal{H}_1$ and $u$ given by \eqref{reparametrizelog}.
Note the similarity to the definition of Ruijsenaars' $R$-function, see
\eqref{eqhypdefR}.

{}From the series expansion of $B_t$ as sum of two balanced ${}_4\phi_3$, we have
in terms of Askey-Wilson parameters \eqref{AWpar},
\begin{equation}\label{AWseries}
\begin{split}
\phi(\gamma;x,\lambda)=&\frac{\bigl(ab,ac;q\bigr)_{\infty}}{\bigl(q/ad;q\bigr)_{\infty}}\,
{}_4\phi_3
\left(\begin{matrix} az^{\pm 1},\widetilde{a}\mu^{\pm 1}\\
ab,ac,ad\end{matrix};q,q\right)\\
+&\frac{\bigl(qb/d, qc/d,\widetilde{a}\mu^{\pm 1}, az^{\pm 1};q\bigr)_{\infty}}
{\bigl(ad/q, q\widetilde{a}\mu^{\pm 1}/ad, qz^{\pm 1}/d;q\bigr)_{\infty}}\,{}_4\phi_3
\left(\begin{matrix} qz^{\pm 1}/d, q\widetilde{a}\mu^{\pm 1}/ad\\
q^2/ad, qb/d, qc/d\end{matrix};q,q\right),
\end{split}
\end{equation}
which shows that $\phi(\gamma;x,\lambda)$ is, up to a $(z,\gamma)$-independent rescaling factor,
the Askey-Wilson function as defined in e.g. \cite{KenS}.

We now re-derive several fundamental
properties of the Askey-Wilson function
using the results of the previous subsection.
Comparing the symmetries of the Askey-Wilson function $\phi(\gamma;x,\lambda)$
to the symmetries of the $R$-function (Proposition \ref{Rsymmetry}),
we have a broken symmetry in the parameters $\gamma$ is broken (from the Weyl group of type
$D_4$ to the Weyl group of type $D_3$).
The most important symmetry (self-duality) is also valid for the
Askey-Wilson function and has played a fundamental role in the study of the
associated generalized Fourier transform (see \cite{KenS}).
Self-duality of the Askey-Wilson function has a natural
interpretation in terms of Cherednik's theory on double affine
Hecke algebras, see \cite{Stok}. Concretely, the symmetries of the
Askey-Wilson function are as follows.
\begin{prop}\label{Awsymmetry}
The Askey-Wilson function $\phi(\gamma;x,\lambda)$ is even in $x$
and $\lambda$ and is self-dual,
\[\phi(\gamma;x,\lambda)=\phi(\gamma;-x,\lambda)=\phi(\gamma;x,-\lambda)=
\phi(\hat{\gamma};\lambda,x).
\]
Furthermore, $\phi(\gamma;x,\lambda)$ has a $W(D_3)$-symmetry in
the parameters $\gamma$, given by
\begin{equation*}
\begin{split}
\phi(\gamma_1,\gamma_0,\gamma_2,\gamma_3;x,\lambda)&=\frac{\bigl(e((-\hat{\gamma}_3+\omega\pm
i\lambda)/\omega_2);q\bigr)_{\infty}}
{\bigl(e((\hat{\gamma}_2+\omega\pm
i\lambda)/\omega_2);q\bigr)_{\infty}}\,\phi(\gamma;x,\lambda),\\
\phi(\gamma_0,\gamma_2,\gamma_1,\gamma_3;x,\lambda)&=\phi(\gamma;x,\lambda),\\
\phi(\gamma_0,\gamma_1,-\gamma_3,-\gamma_2;x,\lambda)&=
\frac{\bigl(e((-\gamma_3+\omega\pm ix)/\omega_2);q\bigr)_{\infty}}
{\bigl(e((\gamma_2+\omega\pm
ix)/\omega_2);q\bigr)_{\infty}}\,\phi(\gamma;x,\lambda).
\end{split}
\end{equation*}
\end{prop}
\begin{proof}
Similarly as for the $R$-function (see Proposition
\ref{Rsymmetry}), the symmetries of the Askey-Wilson function correspond to
the $W(D_5)$-symmetries of $B_t$. Alternatively, all symmetries follow
trivially from the series expansion \eqref{AWseries} of the Askey-Wilson function,
besides its symmetry with respect to $\gamma_0\leftrightarrow \gamma_1$
and $\gamma_2\leftrightarrow -\gamma_3$. These two symmetry
relations are equivalent under duality, since
\[(\gamma_1,\gamma_0,\gamma_2,\gamma_3)^{\hat{}}=(\hat{\gamma}_0,\hat{\gamma}_1,
-\hat{\gamma}_3, -\hat{\gamma}_2),
\]
so we only discuss the symmetry with respect to $\gamma_0\leftrightarrow
\gamma_1$.  By Proposition \ref{EVBsymmetry} we have, with
parameters $t$ given by \eqref{reptrig} and
\eqref{reparametrizelog},
\begin{equation*}
\begin{split}
\phi(\gamma_1,\gamma_0,\gamma_2,\gamma_3;x,\lambda)&=
\frac{\bigl(q,e((\gamma_1+\omega\pm ix)/\omega_2),
e((\hat{\gamma}_0+\omega\pm
i\lambda)/\omega_2);q\bigr)_{\infty}}{2} B_t(ws_{35}t)\\
&=\frac{\bigl(e((-\hat{\gamma}_3+\omega\pm
i\lambda)/\omega_2);q\bigr)_{\infty}}
{\bigl(e((\hat{\gamma}_2+\omega\pm
i\lambda)/\omega_2);q\bigr)_{\infty}}\,\phi(\gamma;x,\lambda),
\end{split}
\end{equation*}
as desired.
\end{proof}
Next we show that the Askey-Wilson function satisfies the same Askey-Wilson second order
difference equation (with step-size $i\omega_1$) as Ruijsenaars'
$R$-function, a result which has previously been derived
from detailed studies of the associated Askey-Wilson polynomials in \cite{IR}, cf.
also \cite{KenS}.
\begin{lemma}\label{AWdifflemmatrig}
The Askey-Wilson function $\phi(\gamma;x,\lambda)$ satisfies the
second order difference equation
\begin{align*}
A(\gamma;x;\omega_1,\omega_2) (\phi(\gamma;x+i\omega_1,\lambda;\omega_1,\omega_2)  -
\phi(\gamma;x,\lambda;\omega_1,\omega_2)) &  +
(x\leftrightarrow -x)\\ & =B(\gamma;\lambda;\omega_1,\omega_2) \phi(\gamma;x,\lambda;\omega_1,\omega_2),
\end{align*}
where $A$ and $B$ are given by
\begin{align*}
A(\gamma;x;\omega_1,\omega_2) & = \frac{\prod_{j=0}^3 \sinh(\pi(i \omega + x + i\gamma_j)/\omega_2)}
{\sinh(2\pi x/\omega_2)\sinh(2\pi(i\omega+x)/\omega_2)}, \\
B(\gamma;\lambda;\omega_1,\omega_2) &=\sinh(\pi(\lambda - i\omega-i\hat\gamma_0)/\omega_2)\sinh(\pi(\lambda+i\omega+i\hat\gamma_0)/\omega_2).
\end{align*}
\end{lemma}
\begin{proof}
Specialize the parameters according to \eqref{reptrig} and \eqref{reparametrizelog}
in \eqref{BAW}. Subsequently express $B_t(t)$, $B_t(\tau_{45}t)$ and
$B_t(\tau_{54}t)$ in terms of $\phi(\gamma;x,\lambda)$,
$\phi(\gamma;x+i\omega_1,\lambda)$ and $\phi(\gamma;x-i\omega_1,\lambda)$ respectively. The
resulting equation is the desired difference
equation.
\end{proof}
\begin{remark}
Denoting $\Phi(z;\mu)=\Phi(a,b,c,d;z,\mu)$ for the Askey-Wilson
function in the usual Askey-Wilson parameters, Lemma
\ref{AWdifflemmatrig} becomes the Askey-Wilson second order difference equation
\[
A(z)(\Phi(qz,\mu)-\Phi(z,\mu)) + A(z^{-1})(\Phi(z/q,\mu)-\Phi(z,\mu)) =
\bigl(\widetilde{a} (\gamma+\gamma^{-1})-1-\widetilde{a}^2\bigr)\Phi(z,\mu),
\]
where
\[
A(z)= \frac{(1-az)(1-bz)(1-cz)(1-dz)}{(1-qz^2)(1-z^2)}
\]
and $\widetilde{a}=e((\hat{\gamma}_0+\omega)/\omega_2)$.
\end{remark}

We have now seen that the $R$-function $R(\gamma;x,\lambda;\omega_1,\omega_2)$ as well as the Askey-Wilson function
$\phi(\gamma;x,\lambda;\omega_1,\omega_2)$ are solutions to the eigenvalue problem
\begin{equation}\label{eigenvalueproblem}
\mathcal{L}_\gamma^{\omega_1,\omega_2}f=B(\gamma;\lambda;\omega_1,\omega_2)f
\end{equation}
for the Askey-Wilson second order difference operator $\mathcal{L}_{\gamma}^{\omega_1,\omega_2}$
\eqref{AWoperator} with step-size $i\omega_1$. These two solutions have essentially different behaviour in the $i\omega_2$-step
direction: the Askey-Wilson function $\phi(\gamma;x,\lambda)$ is
$i\omega_2$-periodic, while the $R$-function $R(\gamma;x,\lambda)$ is $\omega_1\leftrightarrow\omega_2$
invariant (hence is also an eigenfunction of the Askey-Wilson second
order difference operator $\mathcal{L}_{\gamma}^{\omega_2,\omega_1}$ with step-size $i\omega_2$, with eigenvalue
$B(\gamma;\lambda;\omega_2,\omega_1)$). On the other hand,
note that $\widetilde{\tau}=-\omega_2/\omega_1\in\mathbb{H}_+$ and that
\[A(\gamma;x;\omega_2,\omega_1)=A(-\gamma;-x;-\omega_2,\omega_1),\qquad
B(\gamma;\lambda;\omega_2,\omega_1)=B(-\gamma;\lambda;-\omega_2,\omega_1)
\]
with $-\gamma=(-\gamma_0,-\gamma_1,-\gamma_2,-\gamma_3)$, so that the Askey-Wilson function
$\phi(-\gamma;x,\lambda;-\omega_2,\omega_1)$ (with associated modular inverted base $\widetilde{q}=e(-\omega_2/\omega_1)$)
does satisfy the Askey-Wilson second order difference equation
\begin{equation}\label{trigAWinversion}
\bigl(\mathcal{L}_{\gamma}^{\omega_2,\omega_1}\phi(-\gamma;\,\cdot\,,\lambda;-\omega_2,\omega_1)\bigr)(x)=
B(\gamma;\lambda;\omega_2,\omega_1)\phi(-\gamma;x,\lambda;-\omega_2,\omega_1),
\end{equation}
cf. \cite[\S 6.6]{Ruijsproc}. In the next section we match the hyperbolic theory to the trigonometric
theory, which in particular entails an explicit expression of the
$R$-function in terms of products of Askey-Wilson functions in base $q$ and base $\widetilde{q}$.

Note furthermore that Proposition \ref{Awsymmetry} hints at the fact that
the solution space to the Askey-Wilson eigenvalue problem \eqref{eigenvalueproblem} admits a
natural twisted $W(D_4)$-action on the parameters $\gamma$. In fact, the solution space to \eqref{eigenvalueproblem}
is invariant under permutations of $(\gamma_0,\gamma_1,\gamma_2,\gamma_3)$. Furthermore, a straightforward
computation shows that
\[g(\gamma;\,\cdot\,)^{-1}\circ\mathcal{L}_\gamma^{\omega_1,\omega_2}\circ g(\gamma;\,\cdot\,)=
\mathcal{L}_{(\gamma_0,\gamma_1,-\gamma_3,-\gamma_2)}^{\omega_1,\omega_2}+
B(\gamma;\lambda;\omega_1,\omega_2)-B(\gamma_0,\gamma_1,-\gamma_3,-\gamma_2;\lambda;\omega_1,\omega_2)
\]
for the gauge factor
\[g(\gamma;x)=\frac{\bigl(e((\gamma_2+\omega\pm ix)/\omega_2);q\bigr)_{\infty}}
{\bigl(e((-\gamma_3+\omega\pm ix)/\omega_2);q\bigr)_{\infty}},
\]
which implies that for a given solution $F_\lambda(\gamma_0,\gamma_1,-\gamma_3,-\gamma_2;\,\cdot\,)$ to
the eigenvalue probem
\[\mathcal{L}_{(\gamma_0,\gamma_1,-\gamma_3,-\gamma_2)}^{\omega_1,\omega_2}f=B(\gamma_0,\gamma_1,-\gamma_3,-\gamma_2;\lambda;\omega_1,\omega_2)f,
\]
we obtain a solution
\[\widetilde{F}_\lambda(\gamma;x):=g(\gamma;x)F_\lambda(\gamma_0,\gamma_1,-\gamma_3,-\gamma_2;x)
\]
to the eigenvalue problem \eqref{eigenvalueproblem}. A similar observation forms the starting point of Ruijsenaars' \cite{Ruijs3}
analysis of the $W(D_4)$-symmetries of the $R$-function (see also Section \ref{reparh}).


\begin{remark}\label{vb}
A convenient way to formalize the $W(D_4)$-symmetries of
the eigenvalue problem \eqref{eigenvalueproblem} (in the present trigonometric setting)
is by interpreting $i\omega_2$-periodic solutions to \eqref{eigenvalueproblem}, depending meromorphically on $(\gamma,x,\lambda)$,
as defining a sub-vectorbundle $\widetilde{\Gamma}^0(\omega_1,\omega_2)$
of the meromorphic vectorbundle $\Gamma^0(\omega_1,\omega_2)$ over
\[X=\bigl(\mathbb{C}/\mathbb{Z}\omega_2\bigr)^4\times \mathbb{C}/(\mathbb{Z}i\omega_1+\mathbb{Z}i\omega_2)\times
\mathbb{C}/\mathbb{Z}i\omega_2
\]
consisting of meromorphic functions in $(\gamma,x,\lambda)\in\bigl( \mathbb{C}/\mathbb{Z}\omega_2\bigr)^4\times \mathbb{C}/\mathbb{Z}i\omega_2\times
\mathbb{C}/\mathbb{Z}i\omega_2$.
The above analysis can now equivalently be reformulated as the following property of
$\widetilde{\Gamma}^0(\omega_1,\omega_2)$:
the sub-vectorbundle $\widetilde{\Gamma}^0(\omega_1,\omega_2)$ is $W(D_4)$-invariant with respect to the twisted
$W(D_4)$-action
\begin{equation}\label{twisted}
(\sigma\cdot
f)(\gamma;x,\lambda):=V_\sigma(\gamma;x,\lambda)^{-1}f(\sigma^{-1}\gamma;x,\lambda),\qquad
\sigma\in W(D_4)
\end{equation}
on $\Gamma^0(\omega_1,\omega_2)$,
where $V_\sigma(\gamma;x,\lambda)=h(\sigma^{-1}\gamma;x,\lambda)/h(\gamma;x,\lambda)$
($\sigma\in W(D_4)$) is the $1$-coboundary with $h(\gamma;x,\lambda)=h(\gamma;x,\lambda;\omega_1,\omega_2)$
e.g. given by
\begin{equation}\label{h}
h(\gamma;x,\lambda;\omega_1,\omega_2)=
\frac{\theta\bigl(e((\gamma_3-\hat{\gamma}_0+ix)/\omega_2);q\bigr)}
{\theta\bigl(e((\omega-\gamma_3-ix)/\omega_2);q\bigr)\prod_{j=0}^3\bigl(e((\omega-\gamma_j+ix)/\omega_2);q\bigr)_{\infty}}
\in\Gamma_0(\omega_1,\omega_2)^\times,
\end{equation}
and where $W(D_4)$ acts on the $\gamma$ parameters by permutations
and even sign changes. By a straightforward analysis using Casorati-determinants
and the asymptotically free solutions to the eigenvalue problem
\eqref{eigenvalueproblem}, one can furthermore show that $\widetilde{\Gamma}^0(\omega_1,\omega_2)$
is a (trivial) meromorphic vectorbundle over $X$ of rank two
(compare with the general theory on difference equations in \cite{Put}).
\end{remark}

We end this subsection by expressing the Askey-Wilson function $\phi(\gamma;x,\lambda)$
in terms of the trigonometric integrals $E_t$ and $V_t$ using
Proposition \ref{interconnectionprop}. Note its close resemblance
with the hyperbolic case, cf. Theorem \ref{Eulerreph}.
\begin{lemma}\label{phiEV}
{\bf a)} We have
\begin{equation*}
\begin{split}
\phi(\gamma;x,\lambda)&=\frac{\bigl(q;q\bigr)_{\infty}}{2}
\frac{\bigl(e((\hat{\gamma}_0+\omega-i\lambda)/\omega_2),
e((\hat{\gamma}_1+\omega+i\lambda)/\omega_2),
e((\hat{\gamma}_2+\omega+i\lambda)/\omega_2);q\bigr)_{\infty}}
{\bigl(e((-\hat{\gamma}_3+\omega-i\lambda)/\omega_2);q\bigr)_{\infty}}\\
&\times\frac{\prod_{j=0}^2\bigl(e((\gamma_j+\omega\pm
ix)/\omega_2);q\bigr)_{\infty}}{\bigl(e((-\gamma_3+\omega\pm
ix)/\omega_2);q\bigr)_{\infty}}E_t(t)
\end{split}
\end{equation*}
with
\begin{equation*}
\begin{split}
t_1&=e((-\frac{3\omega}{2}+\gamma_3-\frac{\hat{\gamma}_0}{2}+\frac{i\lambda}{2})/\omega_2),\qquad
t_2=e((\frac{\omega}{2}+\gamma_2-\frac{\hat{\gamma}_0}{2}+\frac{i\lambda}{2})/\omega_2),\\
t_3&=e((\frac{\omega}{2}+\gamma_1-\frac{\hat{\gamma}_0}{2}+\frac{i\lambda}{2})/\omega_2),\qquad\qquad
t_4=e((\frac{\omega}{2}+\gamma_0-\frac{\hat{\gamma}_0}{2}+\frac{i\lambda}{2})/\omega_2),\\
t_5&=e((\frac{\omega}{2}+\frac{\hat{\gamma}_0}{2}+ix-\frac{i\lambda}{2})/\omega_2),\qquad\qquad
t_6=e((\frac{\omega}{2}+\frac{\hat{\gamma}_0}{2}-ix-\frac{i\lambda}{2})/\omega_2).
\end{split}
\end{equation*}
{\bf b)} We have
\[
\phi(\gamma;x,\lambda)=\frac{\bigl(q;q\bigr)_{\infty}}{2}\frac{\bigl(e((\hat{\gamma}_0+\omega+i\lambda)/\omega_2),
e((\hat{\gamma}_1+\omega-i\lambda)/\omega_2), e((\hat{\gamma}_2+\omega-i\lambda)/\omega_2);q\bigr)_{\infty}}
{\bigl(e((-\hat{\gamma}_3+\omega+i\lambda)/\omega_2);q\bigr)_{\infty}}\,V_t(t)
\]
with
\begin{equation*}
\begin{split}
t_1&=e((\frac{3\omega}{2}-\gamma_3+\frac{\hat{\gamma}_0}{2}-\frac{i\lambda}{2})/\omega_2),\qquad
t_2=e((-\frac{\omega}{2}-\gamma_2+\frac{\hat{\gamma}_0}{2}-\frac{i\lambda}{2})/\omega_2),\\
t_3&=e((-\frac{\omega}{2}-\gamma_1+\frac{\hat{\gamma}_0}{2}-\frac{i\lambda}{2})/\omega_2),\qquad
t_4=e((-\frac{\omega}{2}-\gamma_0+\frac{\hat{\gamma}_0}{2}-\frac{i\lambda}{2})/\omega_2),\\
t_5&=e((-\frac{\omega}{2}-\frac{\hat{\gamma}_0}{2}+ix+\frac{i\lambda}{2})/\omega_2),\qquad
t_6=e((-\frac{\omega}{2}-\frac{\hat{\gamma}_0}{2}-ix+\frac{i\lambda}{2})/\omega_2).
\end{split}
\end{equation*}
\end{lemma}
\begin{proof}
{\bf a)}
We use Proposition \ref{Awsymmetry} to rewrite $\phi(\gamma;x,\lambda)$
in terms of $\phi(\gamma_0,\gamma_1,-\gamma_3,-\gamma_2;x,-\lambda)$.
Subsequently we use the defining expression of $\phi(\gamma_0,\gamma_1,-\gamma_3,-\gamma_2;x,-\lambda)$
to obtain
\[\phi(\gamma;x,\lambda)=\frac{\bigl(q,e((\hat{\gamma}_1+\omega\pm
i\lambda)/\omega_2), e((\gamma_0+\omega\pm ix)/\omega_2),
e((\gamma_2+\omega\pm ix)/\omega_2);q\bigr)_{\infty}}
{2\bigl(e((-\gamma_3+\omega\pm
ix)/\omega_2);q\bigr)_{\infty}}\,B_t(\xi)
\]
with
\begin{equation*}
\begin{split}
\xi=\bigl(e((-\gamma_0-\gamma_1-2\omega)/\omega_2),1, &e((\hat{\gamma}_1+\omega+i\lambda)/\omega_2),
e((\gamma_0+\omega-ix)/\omega_2),e((\gamma_0+\omega+ix)/\omega_2),\\
&e((\hat{\gamma}_1+\omega-i\lambda)/\omega_2),
e((\gamma_2-\gamma_0)/\omega_2),
e((\gamma_3-\gamma_0-2\omega)/\omega_2)\bigr).
\end{split}
\end{equation*}
With this specific ordered set $\xi$ of parameters we apply
Proposition \ref{interconnectionprop} to rewrite $B_t(\xi)$ in
terms of $E_t$, which results in the desired identity.\\
{\bf b)} This follows from applying Proposition \ref{interconnectionprop}
directly to the definition \eqref{phidef} of $\phi(\gamma;x,\lambda)$.
\end{proof}

Using the expression of the Askey-Wilson function in terms of $V_t$ and
using Lemma \ref{series2},
we thus obtain an expression of the Askey-Wilson function as
very-well-poised ${}_8\phi_7$ series.


\section{Hyperbolic versus trigonometric theory}


\subsection{Hyperbolic versus trigonometric gamma functions}

We fix throughout this section periods $\omega_1,\omega_2\in\mathbb{C}$
with $\Re(\omega_1)>0$, $\Re(\omega_2)>0$ and
$\tau=\omega_1/\omega_2\in\mathbb{H}_+$. We set
\[q=q_{\omega_1,\omega_2}=e(\omega_1/\omega_2),\qquad
\widetilde{q}=\widetilde{q}_{\omega_1,\omega_2}=e(-\omega_2/\omega_1)
\]
where $e(x)=\exp(2\pi ix)$ as before, so that $|q|, |\widetilde{q}|<1$.

Shintani's \cite{Shintani} product expansion is
\begin{equation}\label{Sh}
G(\omega_1,\omega_2;x)=e\left(-\frac{1}{48}\bigl(\frac{\omega_1}{\omega_2}+\frac{\omega_2}{\omega_1}\bigr)\right)
e\left(-\frac{x^2}{4\omega_1\omega_2}\right)\frac{\bigl(e((ix+\omega)/\omega_2);q\bigr)_{\infty}}
{\bigl(e((ix-\omega)/\omega_1);\widetilde{q}\bigr)_{\infty}},
\end{equation}
where $\omega=\frac{1}{2}(\omega_1+\omega_2)$ as before. For a
proof of \eqref{Sh}, see \cite[Prop. A.1]{Stokman}. In other
words, the product expansion \eqref{Sh} expresses the hyperbolic
gamma function as a quotient of two trigonometric gamma functions
(one in base $q$, the other in the modular inverted base
$\widetilde{q}$). In this section
we explicitly write the base-dependence; e.g. we write
$S_t(t;q)$ ($t\in\mathcal{H}_1$) to denote the trigonometric
hypergeometric function $S_t(t)$ in base $q$.


\subsection{Hyperbolic versus trigonometric hypergeometric integrals}

We explore \eqref{Sh} to relate the hyperbolic integrals
to their trigonometric analogs. We start with the hyperbolic
hypergeometric function $S_h(u)$ ($u\in\mathcal{G}_{2i\omega}$).
For $u\in\mathcal{G}_{2i\omega}$ we write
\begin{equation}\label{parameterrel}
t_j=e\bigl((iu_j+\omega)/\omega_2\bigr),\qquad
\widetilde{t}_j=e\bigl((iu_j-\omega)/\omega_1\bigr),\qquad
j=1,\ldots,8.
\end{equation}
Observe that $\prod_{j=1}^8t_j=q^2$ and
$\prod_{j=1}^8\widetilde{t}_j=\widetilde{q}^6$.

\begin{thm}\label{hyptrigS}
As meromorphic functions of $u\in\mathcal{G}_{2i\omega}$ we have
\begin{equation}\label{hyptrigSfor}
\begin{split}
S_h(u)&=\omega_2
e\bigl((2\omega^2+\sum_{j=1}^5u_j^2-u_6^2+u_7^2-u_8^2)/2\omega_1\omega_2\bigr) \frac{\bigl(\widetilde{q},\widetilde{q};\widetilde{q}\bigr)_{\infty}}{2} \\
& \qquad \times
\frac{\prod_{j=1}^5\theta\bigl(\widetilde{t}_j\widetilde{t}_7/\widetilde{q};\widetilde{q}\bigr)}
{\theta\bigl(\widetilde{t}_7/\widetilde{t}_6;\widetilde{q}\bigr)}
U_t\bigl(\widetilde{q}^{\frac{3}{2}}/\widetilde{t}_8,\widetilde{q}^{\frac{1}{2}}/\widetilde{t}_1,\ldots,
\widetilde{q}^{\frac{1}{2}}/\widetilde{t}_6,\widetilde{q}^{\frac{3}{2}}/\widetilde{t}_7;\widetilde{q}\bigr)
S_t\bigl(t_6/q,t_1,\ldots,t_5,t_7,t_8/q;q\bigr)\\
& \qquad +\ (u_6\leftrightarrow u_7),
\end{split}
\end{equation}
with the parameters $t_j$ and $\widetilde{t}_j$ given by
\eqref{parameterrel}.
\end{thm}
\begin{proof}
We put several additional conditions on the parameters, which can later
be removed by analytic continuity.
We assume that $\omega_1,-\omega_2\in\mathbb{H}_+$ and that $\Re(i\omega)<0$.
We furthermore choose parameters $u\in\mathcal{G}_{2i\omega}$ satisfying $\Re(u_j-i\omega)>0$
and $\Im(u_j-i\omega)<0$ for $j=1,\ldots,8$. Then
\begin{equation*}
\begin{split}
S_h(u)&=\int_{\mathbb{R}}\frac{G(i\omega\pm
2x)}{\prod_{j=1}^8G(u_j\pm
x)}dx\\
&=e\left(\frac{7}{24}\bigl(\frac{\omega_1}{\omega_2}+\frac{\omega_2}{\omega_1}\bigr)\right)
e\bigl((\omega^2+\sum_{j=1}^8u_j^2)/2\omega_1\omega_2\bigr)\int_{\mathbb{R}}W(x)\widetilde{W}(x)dx,
\end{split}
\end{equation*}
where
\begin{equation*}
\begin{split}
W(x)&=\frac{\bigl(e(\pm 2ix/\omega_2);q\bigr)_{\infty}}
{\prod_{j=1}^8\bigl(t_je(\pm
ix/\omega_2);q\bigr)_{\infty}},\\
\widetilde{W}(x)&=e\bigl(2x^2/\omega_1\omega_2\bigr)
\frac{\prod_{j=1}^8\bigl(\widetilde{t}_je(\pm
ix/\omega_1);\widetilde{q}\bigr)_{\infty}}
{\bigl(\widetilde{q}e(\pm
2ix/\omega_1);\widetilde{q}\bigr)_{\infty}}
\end{split}
\end{equation*}
by \eqref{Sh}. Using Cauchy's Theorem and elementary asymptotic
estimates of the integrand, we may rotate the integration contour $\mathbb{R}$
to $i\omega_2\mathbb{R}$. Since the factor $W(x)$ is
$i\omega_2$-periodic, we can fold the resulting integral,
interchange summation and integration by Fubini's Theorem, to obtain
the expression
\[S_h(u)=e\left(\frac{7}{24}\bigl(\frac{\omega_1}{\omega_2}+\frac{\omega_2}{\omega_1}\bigr)\right)
e\bigl((\omega^2+\sum_{j=1}^8u_j^2)/2\omega_1\omega_2\bigr)\int_{0}^{i\omega_2}W(x)F(x)dx,
\]
where
\begin{equation*}
\begin{split}
F(x)&=\sum_{n=-\infty}^{\infty}\widetilde{W}(x+mi\omega_2)\\
&=\widetilde{W}(x)\,
{}_{10}\psi_{10}\left(\begin{matrix}
\widetilde{q}e(ix/\omega_1), -\widetilde{q}e(ix/\omega_1),
\{\widetilde{q}\widetilde{t}_j^{-1}e(ix/\omega_1)\}_{j=1}^8\\
e(ix/\omega_1), -e(ix/\omega_1),
\{\widetilde{t}_je(ix/\omega_1)\}_{j=1}^8\end{matrix}\,;\widetilde{q},\widetilde{q}\right).
\end{split}
\end{equation*}
At this stage we have to resort to \cite[(5.6.3)]{GenR}, which
expresses a very-well-poised ${}_{10}\psi_{10}$ bilateral series
as a sum of three very-well-poised ${}_{10}\phi_9$ unilateral
series. This results in the formula
\begin{equation*}
\begin{split}
F(x)&=e\bigl(2x^2/\omega_1\omega_2\bigr)\theta\bigl(\widetilde{t}_6e(\pm
ix/\omega_1), \widetilde{t}_7e(\pm
ix/\omega_1);\widetilde{q}\bigr)\\
&\times\frac{\bigl(\widetilde{q};\widetilde{q}\bigr)_{\infty}\prod_{j=1}^5
\bigl(\widetilde{t}_j\widetilde{t}_8/\widetilde{q},
\widetilde{q}\widetilde{t}_j/\widetilde{t}_8;\widetilde{q}\bigr)_{\infty}}
{\bigl(\widetilde{q}^2/\widetilde{t}_6\widetilde{t}_8,
\widetilde{q}^2/\widetilde{t}_7\widetilde{t}_8,
\widetilde{t}_8/\widetilde{t}_6, \widetilde{t}_8/\widetilde{t}_7,
\widetilde{q}^3/\widetilde{t}_8^2;\widetilde{q}\bigr)_{\infty}}\,
{}_{10}W_9\bigl(\widetilde{q}^2/\widetilde{t}_8^2;\{\widetilde{q}^2/\widetilde{t}_j\widetilde{t}_8\}_{j=1}^7;
\widetilde{q},\widetilde{q}\bigr)\\
&+(u_8;u_6,u_7),
\end{split}
\end{equation*}
where $(u_8;u_6,u_7)$ means cyclic permutation of the parameters
$(u_8,u_6,u_7)$. Note that the ${}_{10}\phi_9$ series in the expression of $F(x)$
are independent of $x$. Combining Jacobi's inversion formula, the Jacobi triple product identity
and the modularity
\begin{equation}\label{eta}
\frac{\bigl(q;q\bigr)_{\infty}}{\bigl(\widetilde{q};\widetilde{q}\bigr)_{\infty}}
=\sqrt{\frac{\omega_2}{-i\omega_1}}\,e\left(-\frac{1}{24}\bigl(\frac{\omega_1}{\omega_2}+\frac{\omega_2}{\omega_1}
\bigr)\right)
\end{equation}
of Dedekind's eta function, we obtain
\begin{equation}\label{Jacobi}
\theta\bigl(e(u/\omega_1);\widetilde{q}\bigr)=
e\left(-\frac{1}{24}\bigl(\frac{\omega_1}{\omega_2}+\frac{\omega_2}{\omega_1}\bigr)\right)
e\bigl((u+\omega)^2/2\omega_1\omega_2\bigr)\theta\bigl(e(-u/\omega_2);q\bigr)
\end{equation}
for the rescaled Jacobi theta function $\theta(\cdot)$,
see e.g. \cite{GenR} or \cite{Stokman}. As a result, we can
rewrite the theta functions in the expression of $F(x)$ as
theta functions in base $q$,
\begin{equation*}
\begin{split}
e\bigl(2x^2/\omega_1\omega_2\bigr)\theta\bigl(\widetilde{t}_6e(\pm
ix/\omega_1),\widetilde{t}_7e(\pm
ix/\omega_1);\widetilde{q}\bigr)&=
e\left(-\frac{1}{6}\bigl(\frac{\omega_1}{\omega_2}+\frac{\omega_2}{\omega_1}\bigr)\right)
e\bigl(-(u_6^2+u_7^2)/\omega_1\omega_2\bigr)\\
&\times\theta\bigl(qe(\pm ix/\omega_2)/t_6, qe(\pm
ix/\omega_2)/t_7;q\bigr).
\end{split}
\end{equation*}
We thus obtain the expression
\begin{equation}\label{ShuHT}
\begin{split}
S_h(u)&=C(u_8;u_6,u_7)\int_{0}^{i\omega_2}\theta\bigl(qe(\pm ix/\omega_2)/t_6,
qe(\pm ix/\omega_2)/t_7;q\bigr)W(x)dx+ (u_8;u_6,u_7)\\
&=i\omega_2C(u_8;u_6,u_7)S_t\bigl(t_6/q,t_1,\ldots,t_5,t_8,t_7/q;q)+
(u_8;u_6,u_7),
\end{split}
\end{equation}
where we have used that $|t_j|<1$ for $j=1,\ldots,6$, with
\begin{equation*}
\begin{split}
C(u_8;u_6,u_7)&=
e\left(\frac{1}{8}\bigl(\frac{\omega_1}{\omega_2}+\frac{\omega_2}{\omega_1}\bigr)\right)
e\left((\omega^2+\sum_{j=1}^5u_j^2-u_6^2-u_7^2+u_8^2)/2\omega_1\omega_2\right)\\
&\times\frac{\bigl(\widetilde{q};\widetilde{q}\bigr)_{\infty}\prod_{j=1}^5
\bigl(\widetilde{t}_j\widetilde{t}_8/\widetilde{q},
\widetilde{q}\widetilde{t}_j/\widetilde{t}_8;\widetilde{q}\bigr)_{\infty}}
{\bigl(\widetilde{q}^2/\widetilde{t}_6\widetilde{t}_8,
\widetilde{q}^2/\widetilde{t}_7\widetilde{t}_8,
\widetilde{t}_8/\widetilde{t}_6, \widetilde{t}_8/\widetilde{t}_7,
\widetilde{q}^3/\widetilde{t}_8^2;\widetilde{q}\bigr)_{\infty}}\,
{}_{10}W_9\bigl(\widetilde{q}^2/\widetilde{t}_8^2;\{\widetilde{q}^2/\widetilde{t}_j\widetilde{t}_8\}_{j=1}^7;
\widetilde{q},\widetilde{q}\bigr).
\end{split}
\end{equation*}
We thus have obtained an expression of $S_h(u)$ as a sum of three
trigonometric hypergeometric functions $S_t$ in base $q$, with
coefficients expressed as very-well-poised ${}_{10}\phi_9$ series
in base $\widetilde{q}$. The next step is to use a three term
transformation for $S_t$ to write $S_h(u)$ as a sum of two
trigonometric hypergeometric functions $S_t$ in base $q$, with
coefficients now being a sum of two very-well-poised
${}_{10}\phi_9$ series.

Concretely, we consider the contiguous relation \eqref{eqellcont2}
for $S_e$ with $p\leftrightarrow q$ and with parameters
specialized to $(t_1,\ldots,t_5, pt_6, pt_7, t_8)$.
Taking the limit $p\rightarrow 0$ leads to
the three term transformation
\begin{equation}\label{diff3}
S_t\bigl(t_6/q,t_1,\ldots,t_5,t_8,t_7/q;q\bigr)=
\frac{\theta\bigl(t_7/t_8;q\bigr)}{\theta\bigl(t_7/t_6;q\bigr)}
\prod_{j=1}^5\frac{\theta\bigl(t_jt_6;q\bigr)}{\theta\bigl(t_jt_8;q\bigr)}\,
S_t\bigl(t_7/q,t_1,\ldots,t_6,t_8/q;q\bigr)+ (u_6\leftrightarrow
u_7)
\end{equation}
for $S_t$. Rewriting the coefficients in \eqref{diff3} in base $\widetilde{q}$
using the Jacobi inversion formula \eqref{Jacobi},
\[\frac{\theta\bigl(t_7/t_8;q\bigr)}{\theta\bigl(t_7/t_6;q\bigr)}
\prod_{j=1}^5\frac{\theta\bigl(t_jt_6;q\bigr)}{\theta\bigl(t_jt_8;q\bigr)}=
e\bigl((u_6^2-u_8^2)/\omega_1\omega_2\bigr)
\frac{\theta\bigl(\widetilde{t}_8/\widetilde{t}_7;\widetilde{q}\bigr)}
{\theta\bigl(\widetilde{t}_6/\widetilde{t}_7;\widetilde{q}\bigr)}
\prod_{j=1}^5\frac{\theta\bigl(\widetilde{q}^2/\widetilde{t}_j\widetilde{t}_6;\widetilde{q}\bigr)}
{\theta\bigl(\widetilde{q}^2/\widetilde{t}_j\widetilde{t}_8;\widetilde{q}\bigr)}
\]
and using the resulting three term transformation in
\eqref{ShuHT}, we obtain
\[
S_h(u)=D(u_6,u_7)S_t(t_7/q,t_1,\ldots,t_6,t_8/q;q)+
(u_6\leftrightarrow u_7)\]
with
\[D(u_6,u_7)=i\omega_2\Bigl(e\bigl((u_6^2-u_8^2)/\omega_1\omega_2\bigr)
\frac{\theta\bigl(\widetilde{t}_8/\widetilde{t}_7;\widetilde{q}\bigr)}
{\theta\bigl(\widetilde{t}_6/\widetilde{t}_7;\widetilde{q}\bigr)}
\prod_{j=1}^5\frac{\theta\bigl(\widetilde{q}^2/\widetilde{t}_j\widetilde{t}_6;\widetilde{q}\bigr)}
{\theta\bigl(\widetilde{q}^2/\widetilde{t}_j\widetilde{t}_8;\widetilde{q}\bigr)}C(u_8;u_6,u_7)
+C(u_6;u_7,u_8)\Bigr).
\]
The coefficient $D(u_6,u_7)$ is a sum of two very-well-poised ${}_{10}\phi_9$
series in base $\widetilde{q}$, which can be expressed in terms of the trigonometric
integral $U_t$ (in base $\widetilde{q}$) by direct computations using Lemma \ref{series}.
This yields the desired result.
\end{proof}
\begin{remark}\label{frem}
{\bf i)} Note that the $W(E_6)$-symmetry of the trigonometric
integrals $S_t$ and $U_t$ is upgraded to a $W(E_7)$-symmetry in Theorem \ref{hyptrigS} since the second term
in the right hand side of \eqref{hyptrigSfor} is
the first term with the role of $u_6$ and $u_7$ interchanged.

{\bf ii)} Specializing the parameters in Theorem \ref{hyptrigS} to $u\in\mathcal{G}_{2i\omega}$
with $u_1=-u_6$ (so that $t_1t_6=q$ and
$\widetilde{t}_1\widetilde{t}_6=\widetilde{q}$), the left hand
side of the identity can be evaluated by the hyperbolic
Nassrallah-Rahman integral evaluation \eqref{eqhypbetint}.
For the right hand of the identity, the
second term vanishes because $\theta\bigl(\widetilde{t}_1\widetilde{t}_6/\widetilde{q};\widetilde{q})=0$
under the particular parameter specialization.
The remaining product of two trigonometric integrals can be evaluated
by Corollary \ref{evaltrig}. The equality of both sides of the resulting identity
can be reconfirmed using \eqref{Sh} and \eqref{eta}. It follows from this argument that the evaluation of the
hyperbolic Nassrallah-Rahman integral is in fact a consequence of fusing
trigonometric identities, an approach to hyperbolic beta integrals which was analyzed in detail
in \cite{Stokman}.

{\bf (iii)} More generally, specializing \eqref{hyptrigSfor} at generic $u\in\mathcal{G}_{2i\omega}$ satisfying
$u_1+u_6= n i \omega_1 + mi \omega_2$ ($n,m \in \mathbb{Z}_{\geq 0}$), the
second term on the right hand side of \eqref{hyptrigSfor} still vanishes while the first term reduces to
the product of two terminating very-well-poised ${}_{10}\phi_9$ series, one in base $q$ and the other in base $\widetilde{q}$.
The terminating ${}_{10}\phi_9$ series is Rahman's \cite{Ra} biorthogonal rational ${}_{10}\phi_9$ function (cf. Remark \ref{biorthtrig}{\bf (ii)}),
while the resulting expression for $S_h$ is the corresponding two-index hyperbolic analogue of Rahman's biorthogonal rational function,
considered by Spiridonov \cite[\S 8.3]{Spir2} (cf. Remark \ref{biorthell} on the elliptic level).
\end{remark}

\begin{cor}\label{hyptrigE}
We have
\begin{equation*}
\begin{split}
E_h(u)&=\omega_2e\left(\frac{2\omega^2+\sum_{j=1}^4u_j^2-u_5^2+u_6^2}{2\omega_1\omega_2}\right)\\
&\times\frac{\bigl(\widetilde{q},\widetilde{q};\widetilde{q}\bigr)_{\infty}}{2}
\frac{\prod_{j=1}^4\theta\bigl(\widetilde{t}_j\widetilde{t}_6/\widetilde{q};\widetilde{q}\bigr)}
{\theta\bigl(\widetilde{t}_6/\widetilde{t}_5;\widetilde{q}\bigr)}\,
V_t\bigl(\frac{\widetilde{q}^{\frac{3}{2}}}{\widetilde{t}_6},\frac{\widetilde{q}^{\frac{1}{2}}}{\widetilde{t}_1},
\ldots,\frac{\widetilde{q}^{\frac{1}{2}}}{\widetilde{t}_5};\widetilde{q}\bigr)
E_t\bigl(\frac{t_5}{q},t_1,\ldots,t_4,t_6;q\bigr)\\
&+(u_5\leftrightarrow u_6)
\end{split}
\end{equation*}
as meromorphic functions in $u\in\mathbb{C}^6$, where
$t_j=e\bigl((iu_j+\omega)/\omega_2\bigr)$ and $\widetilde{t}_j=e\bigl((iu_j-\omega)/\omega_1\bigr)$
\textup{(}$j=1,\ldots,6$\textup{)} as before.
\end{cor}
\begin{proof}
For generic $u\in \mathcal{G}_{2i\omega}$ we have
\[E_h(u_2,\ldots,u_7)=\lim_{s\rightarrow\infty}S_h(u_1+s,u_2,\ldots,u_7,u_8-s)e\left(
\frac{(u_8-u_1-2s)(u_1+u_8)}{2\omega_1\omega_2}\right)
\]
under suitable parameter restraints by Proposition
\ref{prophyplimE}. By Theorem \ref{hyptrigS}, we alternatively
have
\begin{equation*}
\begin{split}
\lim_{s\rightarrow\infty}&S_h(u_1+s,u_2,\ldots,u_7,u_8-s)e\left(
\frac{(u_8-u_1-2s)(u_1+u_8)}{2\omega_1\omega_2}\right)=\\
&\lim_{s\rightarrow\infty}\omega_2e\left(\frac{2\omega^2+\sum_{j=2}^5u_j^2-u_6^2+u_7^2}{2\omega_1\omega_2}\right)
\frac{\bigl(\widetilde{q},\widetilde{q};\widetilde{q}\bigr)_{\infty}}{2}
\frac{\theta\bigl(\widetilde{t}_1\widetilde{t}_7\widetilde{u}/\widetilde{q};\widetilde{q}\bigr)
\prod_{j=2}^5\theta\bigl(\widetilde{t}_j\widetilde{t}_7/\widetilde{q};\widetilde{q}\bigr)}
{\theta\bigl(\widetilde{t}_7/\widetilde{t}_6;\widetilde{q}\bigr)}\\
&\qquad\times
U_t\bigl(\widetilde{q}^{\frac{3}{2}}/\widetilde{t}_7,\widetilde{q}^{\frac{1}{2}}/\widetilde{t}_2,\ldots,
\widetilde{q}^{\frac{1}{2}}/\widetilde{t}_6,\widetilde{q}^{\frac{1}{2}}/\widetilde{t}_1\widetilde{u},
\widetilde{q}^{\frac{3}{2}}\widetilde{u}/\widetilde{t}_8;\widetilde{q}\bigr)
S_t\bigl(t_6/q,t_2,\ldots,t_5,t_7,t_1u,t_8/qu;q\bigr)\\
&\qquad\qquad
+(u_6\leftrightarrow u_7)
\end{split}
\end{equation*}
where $u=e\bigl(is/\omega_2\bigr)$ and
$\widetilde{u}=e\bigl(is/\omega_1\bigr)$. We have $u,\widetilde{u}\rightarrow 0$
as $s\rightarrow \infty$ since $\Re(\omega_1),\Re(\omega_2)>0$, hence application of
Proposition \ref{trigellprop} gives the right hand side of the desired identity with respect
to the parameters $(u_2,\ldots,u_7)$.
\end{proof}

\begin{remark}
Alternatively Corollary \ref{hyptrigE} can be proved by repeating
the arguments of Theorem \ref{hyptrigS}. The argument
simplifies, since one now only needs the expression \cite[(5.6.1)]{GenR} of
a very-well-poised ${}_8\psi_8$ as a sum of two very-well-poised
${}_8\phi_7$ series, and one does not need to use three term
transformations for the trigonometric integrals.
\end{remark}

We conclude this section by relating the $R$-function to the Askey-Wilson function using
Corollary \ref{hyptrigE}. The answer deviates from Ruijsenaars' \cite[\S 6.6]{Ruijsproc} hunch
that $R$ is (up to an elliptic prefactor) the product of an Askey-Wilson function in base $q$
and an Askey-Wilson function in base $\widetilde{q}$:  it is the appearance below of two such terms
which upgrades the $W(D_3)$-symmetry of the Askey-Wilson functions to the $W(D_4)$-symmetry of $R$
(cf. Remark \ref{frem}{\,\bf i)}). For notational convenience, we
write $w_0=-1\in W(D_4)$ for the longest Weyl group element,
acting as $w_0\gamma=-\gamma$ on the Askey-Wilson parameters $\gamma$. We
define $\psi(\gamma;x,\lambda)=\psi(\gamma;x,\lambda;\omega_1,\omega_2)$ by
\begin{equation*}
\begin{split}
\psi(\gamma;x,\lambda;\omega_1,\omega_2)&=\frac{h(\hat{\gamma};\lambda,x;\omega_1,\omega_2)}
{h(-\hat{\gamma};\lambda,x;\omega_1,\omega_2)}\bigl(w_0\cdot\phi\bigr)(\gamma;x,\lambda;\omega_1,\omega_2)\\
&=\frac{h(\gamma;x,\lambda;\omega_1,\omega_2)h(\hat{\gamma};\lambda,x;\omega_1,\omega_2)}
{h(-\gamma;x,\lambda;\omega_1,\omega_2)h(-\hat{\gamma};\lambda,x;\omega_1,\omega_2)}\,
\phi(-\gamma;x,\lambda;\omega_1,\omega_2),
\end{split}
\end{equation*}
where the gauge factor $h$ is given by \eqref{h}.
Note that $\psi(\gamma;x,\lambda;\omega_1,\omega_2)$
is a self-dual solution of the Askey-Wilson difference equation \eqref{eigenvalueproblem}, see
Remark \ref{vb}. We furthermore define the multiplier
\[M(\gamma;x)=
\frac{\theta\bigl(e((\hat{\gamma}_0-\gamma_3-ix)/\omega_2),
e((\omega+\gamma_2+ix)/\omega_2),
e((\omega+\gamma_3-ix)/\omega_2);q\bigr)}
{\theta\bigl(e((\gamma_3-\hat{\gamma}_0-ix)/\omega_2),
e((\omega+\gamma_0-ix)/\omega_2),
e((\omega+\gamma_1-ix)/\omega_2);q\bigr)},
\]
which is elliptic in $x$ with respect to the period lattice $\mathbb{Z}i\omega_1+\mathbb{Z}i\omega_2$.

\begin{thm}
We have
\[
R(\gamma;x,\lambda;\omega_1,\omega_2)=K(\gamma)M(\gamma;x)
M(\hat{\gamma};\lambda)\phi(s_{23}\gamma;x,\lambda;\omega_1,\omega_2)
\psi(-\gamma;x,\lambda;-\omega_2,\omega_1)+ (\gamma_2\leftrightarrow\gamma_3),
\]
where $s_{23}\gamma=(\gamma_0,\gamma_1,\gamma_3,\gamma_2)$ and with
\begin{equation*}
\begin{split}
K(\gamma;\omega_1,\omega_2)&=\sqrt{-i}\,e\bigl(-\frac{1}{24}(\frac{\omega_1}{\omega_2}+
\frac{\omega_2}{\omega_1})\bigr)e\Bigl(-\frac{3\omega(\gamma_0+\gamma_1+\gamma_2-\gamma_3)}{2\omega_1\omega_2}\Bigr)
\frac{\prod_{j=1}^3G(i\omega+i\gamma_0+i\gamma_j)}
{\theta\bigl(e((\gamma_2-\gamma_3)/\omega_2);q\bigr)}\\
&\times e\Bigl(\frac{-\gamma_0^2-\gamma_1^2-\gamma_2^2+\gamma_3^2-2\gamma_0\gamma_1-2\gamma_0\gamma_2+
2\gamma_0\gamma_3}{4\omega_1\omega_2}\Bigr).
\end{split}
\end{equation*}
\end{thm}
\begin{proof}
Using the second hyperbolic Euler integral representation of $R$ from Theorem \ref{Eulerreph} and
subsequently applying Corollary \ref{hyptrigE}, we obtain an expression of $R(\gamma;x,\lambda;\omega_1,\omega_2)$
in terms of trigonometric integrals $E_t$ and $V_t$ with parameter specializations which allows us to
rewrite them as Askey-Wilson functions by Lemma \ref{phiEV}. This leads to the expression
\begin{equation}\label{RAW}
R(\gamma;x,\lambda;\omega_1,\omega_2)=C(\gamma;x,\lambda;\omega_1,\omega_2)
\phi(s_{23}\gamma;x,\lambda;\omega_1,\omega_2)
\psi(-\gamma;x,\lambda;-\omega_2,\omega_1)+ (\gamma_2\leftrightarrow\gamma_3),
\end{equation}
with the explicit prefactor
\begin{equation*}
\begin{split}
C(\gamma;\lambda,x;\omega_1,\omega_2)&=\sqrt{\frac{\omega_2}{\omega_1}}
e\left(\frac{8\omega^2+8x^2+2(i\omega-i\hat{\gamma}_0-\lambda)^2+(i\omega-2i\gamma_0+i\hat{\gamma}_0+\lambda)^2}
{8\omega_1\omega_2}\right.\\
&\qquad\left.+\frac{(i\omega-2i\gamma_1+i\hat{\gamma}_0+\lambda)^2-(i\omega-2i\gamma_2+i\hat{\gamma}_0+\lambda)^2+
(i\omega-2i\gamma_3+i\hat{\gamma}_0+\lambda)^2}{8\omega_1\omega_2}\right)\\
&\times\frac{\bigl(\widetilde{q};\widetilde{q}\bigr)_{\infty}}{\bigl(q;q\bigr)_{\infty}}\frac{\prod_{j=1}^3G(i\omega+i\gamma_0+i\gamma_j)}
{\theta\bigl(e((\gamma_3-\gamma_2)/\omega_1);\widetilde{q}\bigr)}\prod_{j=0}^2\frac{\bigl(e((-\omega+\gamma_j+ix)/\omega_1);\widetilde{q}\bigr)_{\infty}}
{\bigl(e((-\omega-\gamma_j+ix)/\omega_1);\widetilde{q}\bigr)_{\infty}}\\
&\times \frac{\theta\bigl(e((-\omega+\gamma_3\pm ix)/\omega_1);\widetilde{q}\bigr)\bigl(e((\omega-\gamma_2\pm ix)/\omega_2);q\bigr)_{\infty}
\prod_{j=0}^3G(-i\gamma_j\pm x)}
{\bigl(e((\omega+\gamma_0\pm ix)/\omega_2), e((\omega+\gamma_1\pm ix)/\omega_2), e((\omega+\gamma_3\pm ix)/\omega_2);q\bigr)_{\infty}}\\
&\times\frac{\theta\bigl(e((\gamma_3-\hat{\gamma}_0+ix)/\omega_1);\widetilde{q}\bigr)}
{\theta\bigl((\hat{\gamma}_0-\gamma_3+ix)/\omega_1);\widetilde{q}\bigr)}
\frac{\bigl(e((-\omega+\gamma_3-ix)/\omega_1);\widetilde{q}\bigr)_{\infty}}
{\bigl(e((-\omega-\gamma_3-ix)/\omega_1);q\bigr)_{\infty}}\\
&\times \frac{\theta\bigl(e((\hat{\gamma}_3-\gamma_0+i\lambda)/\omega_1), e((-\omega-\hat{\gamma}_2+i\lambda)/\omega_1),
e((-\omega+\hat{\gamma}_3+i\lambda)/\omega_1);\widetilde{q}\bigr)}
{\theta\bigl(e((\gamma_0-\hat{\gamma}_3+i\lambda)/\omega_1);\widetilde{q}\bigr)}\\
&\times\frac{G(\lambda-i\hat{\gamma}_1,\lambda-i\hat{\gamma}_2, \lambda-i\hat{\gamma}_3)}{G(\lambda+i\hat{\gamma}_0)}\\
&\times\frac{\bigl(e((-\omega+\hat{\gamma}_3-i\lambda)/\omega_1);\widetilde{q}\bigr)_{\infty}}
{\bigl(e((-\omega-\hat{\gamma}_3-i\lambda)/\omega_1);\widetilde{q}\bigr)_{\infty}}
\prod_{j=0}^2\frac{\bigl(e((-\omega+\hat{\gamma}_j+i\lambda)/\omega_1);\widetilde{q}\bigr)_{\infty}}
{\bigl(e((-\omega-\hat{\gamma}_j+i\lambda)/\omega_1);\widetilde{q}\bigr)_{\infty}}\\
&\times\frac{\bigl(e((-\omega-\hat{\gamma}_3+i\lambda)/\omega_1);\widetilde{q}\bigr)_{\infty}}
{\bigl(e((-\omega+\hat{\gamma}_0+i\lambda)/\omega_1), e((-\omega+\hat{\gamma}_1-i\lambda)/\omega_1), e((-\omega+\hat{\gamma}_2-i\lambda)/\omega_1);
\widetilde{q}\bigr)_{\infty}}\\
&\times\frac{\bigl(e((\omega-\hat{\gamma}_2-i\lambda)/\omega_2);q\bigr)_{\infty}}
{\bigl(e((\omega+\hat{\gamma}_0-i\lambda)/\omega_2), e((\omega+\hat{\gamma}_1+i\lambda)/\omega_2), e((\omega+\hat{\gamma}_3+i\lambda)/\omega_2);
q\bigr)_{\infty}}.
\end{split}
\end{equation*}
Elaborate but straightforward computations using \eqref{Sh}, \eqref{eta} and \eqref{Jacobi}
now yields the desired result.
\end{proof}


\begin{thebibliography}{99}
\bibitem{vdBult} F.J. van de Bult, \textit{Ruijsenaars' hypergeometric function and the modular double of
$\mathcal{U}_q(\mathfrak{sl}_2(\mathbb{C}))$}, arxiv:math.QA/0501405, Adv. Math. (to appear).


\bibitem{vD} J.F. van Diejen {\it Integrability of difference Calogero-Moser systems}, J. Math. Phys. {\bf 35} (1994), no. 6, 2983--3004.


\bibitem{Faddeev} L. Faddeev, \textit{Modular double of a Quantum Group}, in ``Conf\'erence Mosh\'e Flato'' 1999, Vol I (Dijon),
 Math. Phys. Stud., \textbf{21}, Kluwer Acad. Publ., Dordrecht, 2000, 149--156.

\bibitem{FrenkelTuraev} I.N. Frenkel, V.G. Turaev, \textit{Elliptic solutions of the Yang-Baxter equation and modular hypergeometric functions},
in: V.I. Arnold, I.M. Gelfand, V.S. Retakh, M. Smirnov (Eds.), The Arnold-Gelfand mathematical seminars, Birkh\"auser, Boston, 1997, 171--204.

\bibitem{GenR} G. Gasper, M. Rahman, Basic Hypergeometric Series, Encyclopedia of Mathematics and its
Applications, Vol. 96, 2nd ed., Cambridge (2004).

\bibitem{G} W. Groenevelt, {\it The Wilson function transform}, Int. Math. Res. Not. {\bf 2003}, no. 52, 2779--2817.

\bibitem{GuptaMasson} D.P. Gupta, D.R. Masson, \textit{Contiguous relations, continued fractions and orthogonality},
Trans. Amer. Math. Soc. \textbf{350} (1998), no. 2, 769--808.

\bibitem{HS} G. Heckman, H. Schlichtkrull, Harmonic Analysis and Special Functions on Symmetric Spaces, Perspectives in Mathematics, Vol. 16 (1994).

\bibitem{IR} M.E.H. Ismail, M. Rahman, {\textit Associated
Askey-Wilson polynomials}, Trans. Amer. Math. Soc. {\bf 328}
(1991), 201--239.

\bibitem{KenS0} E. Koelink, J.V. Stokman, {\it Fourier transforms on the quantum $\hbox{SU}(1,1)$ quantum group}
(with an appendix of M. Rahman), Publ. Res. Math. Sci. {\bf 37} (2001), no. 4, 621--715.

\bibitem{KenS} E. Koelink, J.V. Stokman, \textit{The Askey-Wilson function transform},
Int. Math. Res. Not. \textbf{2001} (22), 1203--1227.

\bibitem{KNR} E. Koelink, Y. van Norden, H. Rosengren, {\it Elliptic} $\hbox{U}(2)$ {\it quantum group and elliptic hypergeometric series},
Comm. Math. Phys {\bf 245} (2004), no. 3, 519--537.

\bibitem{KH} Y. Komori, K. Hikami, {\it Quantum integrability of the generalized elliptic Ruijsenaars models},
J. Phys. A: Math. Gen. {\bf 30} (1997), 4341--4364.

\bibitem{Let} G. Letzter, {\it Quantum zonal spherical functions and Macdonald polynomials}, Adv. Math. {\bf 189} (2004),
no. 1, 88--147.

\bibitem{LiJ} S. Lievens, J. Van der Jeugt, \textit{Symmetry
groups of Bailey's transformations for ${}_{10}\phi_9$-series},
preprint.

\bibitem{Macd} I.G. Macdonald, {\it Orthogonal polynomials associated to root systems}, S{\'e}m. Lothar. Combin. {\bf 45} (2000/01), Art. B45a.

\bibitem{NR} B. Nassrallah, M. Rahman, {\it Projection formulas, a reproducing kernel and a generating function for $q$-Wilson
polynomials}, SIAM J. Math. Anal. {\bf 16} (1985), 186--197.

\bibitem{N} M. Noumi, {\it Macdonald's symmetric polynomials and as zonal spherical functions on some quantum homogeneous spaces},
Adv. Math. {\bf 123} (1996), no. 1, 16--77.

\bibitem{Put} M. van der Put, M.F. Singer, Galois theory of difference equations, Lecture Notes in Mathematics,
1666. Springer-Verlag, Berlin, 1997.

\bibitem{Ra} M. Rahman, {\it An integral representation of a ${}_{10}\phi_9$ and continuous biorthogonal ${}_{10}\phi_9$
rational functions}, Canad. J. Math. {\bf 38} (1986), 605--618.

\bibitem{Rainstrafo} E.M. Rains, \textit{Tranformations of elliptic hypergeometric integrals}, arxiv:math.QA/0309252.

\bibitem{Rainsrecur} E.M. Rains, \textit{Recurrences for elliptic hypergeometric integrals}, arxiv:math.CA/0504285.

\bibitem{Rainslimits} E.M. Rains, \textit{Limits of elliptic hypergeometric integrals}, arxiv:math.CA/0607093.

\bibitem{Ruijs1} S.N.M. Ruijsenaars, \textit{First order analytic difference equations and integrable quantum systems},
J. Math. Phys. \textbf{38} (1997), 1069--1146.

\bibitem{Ruijsproced} S.N.M. Ruijsenaars, \textit{Systems of Calogero-Moser type}, in: Particles and Fields (Banff, AB, 1994),
251--352, CRM Ser. Math. Phys., Springer, New York, 1999.

\bibitem{Ruijs2} S.N.M. Ruijsenaars, \textit{A generalized hypergeometric function satisfying four analytic
difference equations of Askey-Wilson type}, Comm. Math. Phys. \textbf{206} (1999), no. 3, 639--690.

\bibitem{Ruijsproc} S.N.M. Ruijsenaars, \textit{Special functions defined by analytic difference equations},
 Special functions 2000: current perspective and future directions (Tempe, AZ),
281--333, NATO Sci. Ser. II Math. Phys. Chem., {\bf 30}, Kluwer Acad. Publ., Dordrecht, 2001.

\bibitem{Ruijs3} S.N.M. Ruijsenaars, \textit{A generalized hypergeometric function II. Asymptotics and $D_4$ symmetry},
Comm. Math. Phys. \textbf{243} (2003), no 3, 389--412.

\bibitem{Ruijs4} S.N.M. Ruijsenaars, \textit{A generalized hypergeometric function III.
Associated Hilbert space transform},  Comm. Math. Phys. {\bf 243}  (2003),
no. 3, 413--448.

\bibitem{Shintani} T. Shintani, \textit{On a Kronecker limit
formula for real quadratic fields}, J. Fac. Sci. Univ. Tokyo,
Sect. 1A {\bf 24} (1977), 167--199.

\bibitem{Spir} V.P. Spiridonov, \textit{On the elliptic beta function},
Russian Math. Surveys {\bf 56} (2001), no.1, 185--186.

\bibitem{Spiri} V.P. Spiridonov, \textit{Short proofs of the elliptic beta integrals}, arxiv:math.CA/0408369.

\bibitem{Spir2} V.P. Spiridonov, \textit{Classical elliptic hypergeometric functions and their applications},
Rokko Lect. in Math. Vol. 18, Dept. of Meth., Kobe Univ., 2005, 253--287.

\bibitem{Spi4} V.P. Spiridonov, \textit{Theta hypergeometric integrals}, Algebra i Analiz \textbf{15} (2003), 161--215
(St. Petersburg Math J. \textbf{15} (2004), 929--967).

\bibitem{SZ} V.P. Spiridonov, A.S. Zhedanov, {\it Spectral transformation chains and some new biorthogonal rational functions},
Comm. Math. Phys. {\bf 210} (2000), 49--83.

\bibitem{Stok} J.V. Stokman, \textit{Difference Fourier transforms
for nonreduced root systems}, Selecta Math. (N.S.) {\bf 9} (2003),
no. 3, 409--494.

\bibitem{Stokman} J. V. Stokman, \textit{Hyperbolic beta integrals}, Adv. Math. \textbf{190} (2005), 119--160.

\end{thebibliography}
\end{document}